\documentclass[nohyperref]{article}

\usepackage{microtype}
\usepackage{graphicx}
\usepackage{subfigure}
\usepackage{booktabs} 

\usepackage{hyperref}

\usepackage[accepted]{icml2022}

\usepackage{amsmath}
\usepackage{amssymb}
\usepackage{mathtools}
\usepackage{amsthm}

\usepackage[capitalize,noabbrev]{cleveref}

\theoremstyle{plain}
\newtheorem{theorem}{Theorem}[section]

\newtheorem{lemma}[theorem]{Lemma}
\newtheorem{corollary}[theorem]{Corollary}
\theoremstyle{definition}

\theoremstyle{remark}

\newcommand{\cut}[1]{{}}

\makeatletter
\newcommand*{\rom}[1]{\expandafter\@slowromancap\romannumeral #1@}
\makeatother

\newcommand{\vw}{{\mathbf{w}}}
\newcommand{\vx}{{\mathbf{x}}}

\newcommand{\cL}{{\mathcal{L}}}
\newcommand{\cM}{{\mathcal{M}}}
\newcommand{\cN}{{\mathcal{N}}}
\newcommand{\cO}{{\mathcal{O}}}
\newcommand{\cP}{{\mathcal{P}}}

\newcommand{\cR}{{\mathcal{R}}}

 \usepackage[bb=boondox]{mathalfa}

\usepackage{xparse}
\DeclareFontFamily{U}{ntxmia}{}
\DeclareFontShape{U}{ntxmia}{m}{it}{<-> ntxmia }{}
\DeclareFontShape{U}{ntxmia}{b}{it}{<-> ntxbmia }{}
\DeclareSymbolFont{lettersA}{U}{ntxmia}{m}{it}
\SetSymbolFont{lettersA}{bold}{U}{ntxmia}{b}{it}

\ExplSyntaxOn
\NewDocumentCommand{\varmathbb}{m}
 {
  \tl_map_inline:nn { #1 }
  {
    \use:c { varbb##1 }
  }
 }
\tl_map_inline:nn { ABCDEFGHIJKLMNOPQRSTUVWXYZ }
 {
  \exp_args:Nc \DeclareMathSymbol{varbb#1}{\mathord}{lettersA}{\int_eval:n { `#1+67 }}
 }
\exp_args:Nc \DeclareMathSymbol{varbbk}{\mathord}{lettersA}{169}
\ExplSyntaxOff

\newcommand*{\fix}{\mathrm{Fix}\,}
\newcommand*{\zer}{\mathrm{Zer}\,}
\newcommand*{\gra}{\mathrm{Gra}\,}

\newcommand{\prox}{\mathrm{Prox}}

\newcommand{\dom}{\mathrm{dom}\,} 

\newcommand{\reals}{\mathbb{R}}

\definecolor{lightgrey}{gray}{0.8}
\definecolor{medgrey}{gray}{0.6}
\definecolor{darkgrey}{gray}{0.4}

\newcommand{\opA}{{\varmathbb{A}}}
\newcommand{\opB}{{\varmathbb{B}}}

\newcommand{\opG}{{\varmathbb{G}}}

\newcommand{\opI}{{\varmathbb{I}}}
\newcommand{\opJ}{{\varmathbb{J}}}
\newcommand{\opM}{{\varmathbb{M}}}
\newcommand{\opN}{{\varmathbb{N}}}
\newcommand{\opR}{{\varmathbb{R}}}
\newcommand{\opS}{{\varmathbb{S}}}
\newcommand{\opT}{{\varmathbb{T}}}

\newcommand{\Span}{\mathrm{span}}

\newcommand{\resA}{\opJ_\opA}
\newcommand{\al}{\mathbf{A}}

\usepackage[textsize=tiny]{todonotes}

\icmltitlerunning{Exact Optimal Accelerated Complexity for Fixed-Point Iterations}

\begin{document}

\twocolumn[
\icmltitle{Exact Optimal Accelerated Complexity for Fixed-Point Iterations}




\begin{icmlauthorlist}
\icmlauthor{Jisun Park}{to}
\icmlauthor{Ernest K. Ryu}{to}
\end{icmlauthorlist}

\icmlaffiliation{to}{Department of Mathematical Sciences, Seoul National University}

\icmlcorrespondingauthor{Ernest K. Ryu}{ernestryu@snu.ac.kr}

\icmlkeywords{Machine Learning, ICML}

\vskip 0.3in
]



\printAffiliationsAndNotice{}  

\begin{abstract}
    Despite the broad use of fixed-point iterations throughout applied mathematics, the optimal convergence rate of general fixed-point problems with nonexpansive nonlinear operators has not been established.    This work presents an acceleration mechanism for fixed-point iterations with nonexpansive operators, contractive operators, and nonexpansive operators satisfying a H\"older-type growth condition. We then provide matching complexity lower bounds to establish the exact optimality of the acceleration mechanisms in the nonexpansive and contractive setups. Finally, we provide experiments with CT imaging, optimal transport, and decentralized optimization to demonstrate the practical effectiveness of the acceleration mechanism.
\end{abstract}

\section{Introduction}
\label{sec:intro}
The fixed-point iteration with $\opT \colon \reals^n \to \reals^n$ computes
\[
    x_{k+1} = \opT x_k
\]
for $k=0,1,\dots$ with some starting point $x_0\in \reals^n$.
The general rubric of formulating solutions of a problem at hand as fixed points of an operator and then performing the fixed-point iterations is ubiquitous throughout applied mathematics, science, engineering, and machine learning.

Surprisingly, however, the iteration complexity of the abstract fixed-point iteration has not been thoroughly studied. This stands in sharp contrast with the literature on convex optimization algorithms, where convergence rates and matching lower bounds are carefully studied.

In this paper, we establish the exact optimal complexity of fixed-point iterations by providing an accelerated method and a matching complexity lower bound. The acceleration is based on a Halpern mechanism, which follows the footsteps of \citet{Lieder2021halpern,kim2021accelerated,yoon2021accelerated}, and is distinct from Nesterov's acceleration.

\subsection{Preliminaries and notations}
\label{prelim}
We review standard definitions and set up the notation. 

\paragraph{Monotone and set-valued operators.}
We follow standard notation of \citet{bauschke2011convex,ryu2020large}.
For the underlying space, consider $\mathbb{R}^n$ with standard inner product $\langle \cdot,\cdot\rangle$ and norm $\| \cdot \|$, although our results can be extended to infinite-dimensional Hilbert spaces.

We say $\opA$ is an operator on $\mathbb{R}^n$ and write $\opA\colon\mathbb{R}^n\rightrightarrows\mathbb{R}^n$ if $\opA$ maps a point in $\mathbb{R}^n$ to a subset of $\mathbb{R}^n$.
For notational simplicity, also write $\opA x=\opA(x)$.
Write $\gra \opA = \{(x, u) \mid u \in \opA x\}$ for the graph of $\opA$.
Write $\opI\colon\mathbb{R}^n\rightarrow\mathbb{R}^n$ for the identity operator.
We say $\opA\colon\mathbb{R}^n\rightrightarrows\mathbb{R}^n$ is monotone if
\[
\langle \opA x-\opA y,x-y\rangle \ge 0,\qquad\forall x,y\in \mathbb{R}^n,
\]
i.e., if $\langle u-v,x-y\rangle\ge 0$ for all $u\in \opA x$ and $v\in \opA y$.
For $\mu\in (0,\infty)$, say $\opA \colon\mathbb{R}^n\rightrightarrows\mathbb{R}^n$ is $\mu$-strongly monotone if
\[
\langle \opA x-\opA y,x-y\rangle \geq \mu\|x-y\|^2,\qquad\forall x,y\in \mathbb{R}^n.
\]
An operator $\opA$ is maximally monotone if there is no other monotone $\opB$ such that $\gra\opA \subset \gra\opB$ properly, and is maximally $\mu$-strongly monotone if there is no other $\mu$-strongly monotone  $\opB$ such that $\gra\opA \subset \gra\opB$ properly.

For $L\in (0,\infty)$, single-valued operator $\opT \colon\mathbb{R}^n\to\mathbb{R}^n$ is $L$-Lipschitz if
\[
\|\opT x-\opT y\|\le L\|x-y\|,\qquad\forall x,y\in \mathbb{R}^n.
\]
$\opT$ is contractive if it is $L$-Lipschitz with $L <1$ and nonexpansive if it is $1$-Lipschitz.
For $\theta\in(0,1)$, an operator $\opS \colon \reals^n \to \reals^n$ is $\theta$-averaged if
$\opS = (1-\theta)\opI + \theta\opT$ and a nonexpansive operator $\opT$.

Write $\opJ_\opA = (\opI + \opA)^{-1}$ for the resolvent of $\opA$, and $\opR_\opA = 2\resA - \opI$ for the reflected resolvent of $\opA$.
When $\opA$ is maximal monotone, it is well known that $\resA$ is single-valued with $\dom\resA=\reals^n$, $\opR_\opA$ is a nonexpansive operator, and $\resA = \frac{1}{2}\opI + \frac{1}{2}\opR_\opA$ is $1/2$\nobreakdash-averaged.

We say $x_\star \in \reals^n$ is a zero of $\opA$ if $0 \in \opA x_\star$. We say $y_\star$ is a fixed-point of $\opT$ if $\opT y_\star = y_\star$. Write $\zer\opA$ for the set of zeros of $\opA$ and $\fix\opT$ for the set of all fixed-points of $\opT$. For any $x\in \mathbb{R}^n$ such that $x = \resA y$ for some $y\in \mathbb{R}^n$, define $\tilde{\opA} x = y - \resA y$ as the resolvent residual of $\opA$ at $x$.
Note that $\tilde{\opA} x \in \opA x$. 
For any $y\in \mathbb{R}^n$, define $y - \opT y$ as the fixed-point residual of $\opT$ at $y$.

\paragraph{Fixed-point iterations.}
There is a long and rich history of iterative methods for finding a fixed point of an operator $\opT \colon \reals^n \to \reals^n$ \cite{rhoades1991some,brezinski2000convergence,rhoades2001some,berinde2007iterative}.
In this work, we consider the following three:
the \emph{Picard iteration}
\[
    y_{k+1} = \opT y_k,
\]
the \emph{Krasnosel'ski{\u\i}--Mann iteration (KM iteration)}
\[
    y_{k+1} =  \lambda_{k+1} y_k + (1-\lambda_{k+1}) \opT y_k ,
\]
and the \emph{Halpern iteration}
\[
    y_{k+1} = \lambda_{k+1} y_0+(1-\lambda_{k+1}) \opT y_k ,
\]
where $y_0\in\reals^n$ is an initial point and $\{\lambda_k\}_{k\in\mathbb{N}}\subset(0,1)$.
Under suitable assumptions, the $\{y_k\}_{k\in\mathbb{N}}$ sequence of these iterations converges to a fixed point of $\opT$.

\subsection{Prior work}

\paragraph{Fixed-point iterations.}
Picard iteration's convergence with a contractive operator was established by Banach's fixed-point theorem \cite{banach1922operations}.
What we refer to as the Krasnosel'ski{\u\i}--Mann iteration is a generalization of the setups by \citet{krasnosel1955two} and \citet{mann1953mean}.
Its convergence with general nonexpansive operators is due to \citet{martinet1972algorithmes}.
The iteration of \citet{halpern1967fixed} converges for the wider choice of parameter $\lambda_k$ (including $\lambda_k = \frac{1}{k+1}$) due to  \citet{wittmann1992approximation}.
Halpern iteration is later generalized to the sequential averaging method \cite{xu2004viscosity}.
Ishikawa iteration \cite{ishikawa1976fixed} is an iteration with two sequences updated in an alternating manner.
Anderson acceleration \cite{anderson1965} is another acceleration scheme for fixed-point iterations, and it has recently attracted significant interest \cite{walker2011anderson,scieur2020regularized,barre2020convergence,zhang2020globally,bertrand2021anderson}.
A number of inertial fixed-point iterations have also been proposed to accelerate fixed-point iterations \cite{mainge2008convergence,dong2018modified,shehu2018convergence,reich2021inertial}.
Our presented method is optimal (in the sense made precise by the theorems) when compared these prior non-stochastic fixed-point iterations.

\paragraph{Convergence rates of fixed-point iterations.}
The squared fixed-point residual $\|y_k - \opT y_k\|^2$ is the error measure for fixed-point problems that we focus on.
Its convergence to $0$ (without a specified rate) is referred to as asymptotic regularity \cite{browder1966solution}, and it has been established for KM \cite{ishikawa1976fixed,borwein1992krasnoselski} and Halpern \cite{wittmann1992approximation,xu2002iterative}.

The convergence rate of the KM iteration in terms of $\|y_k - \opT y_k\|^2$ was shown to exhibit $\cO(1/k)$-rate \cite{cominetti2014rate,liang2016convergence,bravo2018sharp} and $o(1/k)$-rate \cite{baillon1992optimal,davis2016convergence,matsushita2017convergence} under various setups.
In addition, \citet{borwein2017convergence,lin2021convergence} studied the convergence rate of the distance to solution under additional bounded H\"older regularity assumption.

For the convergence rate of the Halpern iteration in terms of $\|y_k - \opT y_k\|^2$, \citet{leustean2007rates} proved a $\cO(1/(\log{k})^2)$-rate and later \citet{kohlenbach2011quantitative} improved this to a $\cO(1/k)$-rate.
\citet{sabach2017first} first proved the $\cO(1/k^2)$-rate of Halpern iteration, and this rate has been improved in its constant by a factor of $16$ by \citet{Lieder2021halpern}.

\paragraph{Monotone inclusions and splitting methods.}
As we soon establish in \cref{sec:equivalence}, monotone operators are intimately connected to fixed-point iterations.
Splitting methods such as forward-backward splitting (FBS) \cite{bruck1977weak,passty1979ergodic}, augmented Lagrangian method \cite{hestenes1969multiplier,powell1969method}, Douglas--Rachford splitting (DRS) \cite{peaceman1955numerical,douglas1956numerical,lions1979splitting}, alternating direction method of multiplier (ADMM) \cite{gabay1976dual}, Davis--Yin splitting (DYS) \cite{davis2017three}, (PDHG) \cite{chambolle2011first}, and Condat--V\~u \cite{condat2013primal,vu2013splitting} are all fixed-point iterations with respect to specific nonexpansive operators.
Therefore, an acceleration of the abstract fixed-point iteration is applicable to the broad range of splitting methods for monotone inclusions.

\paragraph{Acceleration.}
Since the seminal work by \citet{nesterov1983method} on accelerating gradient methods convex minimization problems, much work as been dedicated to algorithms with faster accelerated rates.
Gradient descent \cite{cauchy1847methode} can be accelerated in terms of function value suboptimality for smooth convex minimization problems \cite{nesterov1983method,KF2016ogm}, smooth strongly convex minimization problems \cite{Nesterov2004convex,van2017fastest,park2021factor,taylor2021optimal,salim2021optimal}, and convex composite minimization problems \cite{guler1992new,beck2009fast}.
Recently, accelerated methods for reducing the squared gradient magnitude for smooth convex minimization \cite{kim2021optimizing,lee2021geometric} and smooth convex-concave minimax optimization \cite{diakonikolas2021potential,yoon2021accelerated}
were presented.

Recently, it was discovered that acceleration is also possible in solving monotone inclusions.
The accelerated proximal point method (APPM) \cite{kim2021accelerated} provides an accelerated $\cO(1/k^2)$-rate of $\|\tilde{\opA}x_k\|^2$ compared to the $\cO(1/k)$-rate of proximal point method (PPM) \cite{Martinet1970ppm,gu2020tight} for monotone inclusions.
\citet{mainge2021accelerated} improved this rate to $o(1/k^2)$ rate with another accelerated variant of proximal point method called CRIPA-S.

\paragraph{Complexity lower bound.}
Under the information-based complexity framework \cite{nemirovski1992information}, complexity lower bound on first-order methods for convex optimization has been thoroughly studied \cite{Nesterov2004convex,drori2017exact,drori2020stoc-complexity,carmon2020stationary1,carmon2021stationary2,drori2022oracle}.
When a complexity lower bound matches an algorithms' guarantee, it establishes optimality of the algorithm \cite{nemirovski1992information,drori2016optimal-kelley,KF2016ogm,taylor2021optimal,yoon2021accelerated,salim2021optimal}.
In the fixed-point theory literature, \citet{Diako2020halpern} provided the lower bound result for the rate of $\langle \opA x_k, x_k-x_\star \rangle$ for variational inequalities with Lipschitz, monotone operator.
\citet{colao2021rate} showed $\Omega\left(1/k^{2-\frac{2}{q}}\right)$ lower bound on  $\|y_k-y_\star\|^2$ for Halpern iterations in $q$-uniformly smooth Banach spaces.
Recently, there has been work establishing complexity lower bounds for the more restrictive ``1-SCLI'' class of algorithms \cite{arjevani2016lower}. The class of 1-SCLI fixed-point iterations includes the KM iteration but not Halpern. Up-to-constant optimality of the KM iteration among 1-SCLI algorithms was proved with the $\Omega(1/k)$ lower bound by \citet{diakonikolas2021potential}.

There also has been recent work on lower bounds for the general class of algorithms (not just 1-SCLI) for fixed-point problems.
\citet{contreras2021optimal} established a $\Omega(1/k^2)$ lower bound on the fixed-point residual for the general Mann iteration, which includes the KM and Halpern iterations, in Banach spaces.
Our $\Omega(1/k^2)$ lower bound of \cref{sec:lower-bound} is more general than the result of \citet{contreras2021optimal} as it applies to all deterministic algorithms, not just Mann iterations.
\citet{diakonikolas2021potential} established a $\Omega(1/k^2)$ lower bound on the squared operator norm for algorithms finding zeros of cocoercive operators, which are equivalent to methods finding fixed points of nonexpansive operators.
Our lower bound  of \cref{sec:lower-bound} improves upon this result (by a constant of about $80$) and establishes exact optimality of the methods in \cref{sec:method}.



\paragraph{Acceleration with restart.}
Restarting is a technique that allows one to render a standard accelerated method to be adaptive to the local structure  \cite{nemirovski1985optimal,nesterov2013gradient,lin2014adaptive,o2015adaptive,kim2018adaptive,fercoq2019adaptive,roulet2020sharpness,ito2021nearly}.
Our method of \cref{sec:restart} was inspired specifically by the restarting scheme of \citet{roulet2020sharpness}.

\paragraph{Performance estimation problem.}
The discovery of the main algorithm of \cref{sec:method} heavily relied on the use of the performance estimation problem (PEP) technique \cite{drori2014performance}. 
Loosely speaking, the PEP is a computer-assisted methodology for finding optimal methods by numerically solving semidefinite programs \cite{drori2014performance,KF2016ogm,taylor2018exactpgm,drori2020efficient,kim2021optimizing}.
We discuss the details of our use of the PEP in Section~\ref{detail:pep} of the appendix.

\subsection{Contributions}
We summarize the contribution of this work as follows. First, we present novel accelerated fixed-point iteration \eqref{algo:oc-halpern} and its equivalent form \eqref{algo:os-ppm} for monotone inclusions. Second, we present exact matching complexity lower bounds and thereby establish the exact optimality of our presented methods. Third, using a restarting mechanism, we extend the acceleration to a broader setup with operators satisfying a H\"older-type growth condition. Finally, we demonstrate the effectiveness of the proposed acceleration mechanism through extensive experiments.

\section{Equivalence of nonexpansive operators and monotone operators}
\label{sec:equivalence}
Before presenting the main content, we quickly establish the equivalence between the fixed-point problem
\[
    \begin{array}{ll}
    \underset{y\in \reals^n}{\mbox{find}} & y= \opT y
    \end{array}
\]
and the monotone inclusion
\[
    \begin{array}{ll}
    \underset{x\in \reals^n}{\mbox{find}} & 0\in \opA x,
    \end{array}
\]
where $\opT\colon\reals^n\to\reals^n$ is $1/\gamma$-Lipschitz with $\gamma\ge1$ and $\opA\colon\reals^n\rightrightarrows\reals^n$ is maximal $\mu$-strongly monotone with $\mu\ge0$.

\begin{lemma} \label{equiv:sm-lips}
    Let $\opT \colon \reals^n\to\reals^n$ and $\opA \colon \reals^n\rightrightarrows\reals^n$.
    If $\opT$ is $1/\gamma$-Lipschitz with $\gamma \ge 1$, then 
    \[
        \opA = \left(\opT + \frac{1}{\gamma} \opI \right)^{-1} \left( 1+ \frac{1}{\gamma} \right) - \opI
    \]
    is maximal $\frac{\gamma-1}{2}$-strongly monotone.
    Likewise, If $\opA$ is maximal $\mu$-strongly monotone with $\mu\ge 0$, then
    \[
        \opT = \left( 1 + \frac{1}{1+2\mu} \right) \resA - \frac{1}{1+2\mu} \opI
    \]
    is $\frac{1}{1+2\mu}$-Lipschitz.
    Under these transformations, $x_\star$ is a zero of $\opA$ if and only if it is a fixed point of $\opT$, i.e., $\zer \opA = \fix \opT$.
\end{lemma}

The equivalence in case of $\gamma=1$ and $\mu=0$ is well known in optimization literature \citep[Theorem~23.8]{bauschke2011convex} \cite{bauschke2012firmly,combettes2018monotone}. This lemma generalizes the equivalence to $\gamma\ge1$ and $\mu\ge0$.
As we see in \cref{proof:equivalence} of the appendix, the equivalence is straightforwardly established using the scaled relative graph (SRG) \cite{ryu2021scaled}, but we also provide a classical proof based on inequalities without using the SRG.

\emph{Remark.}    Since $\opI-\opT = (1+\tfrac{1}{1+2\mu})(\opI-\resA)$, finding an algorithm that effectively reduces $\|y_{N-1} - \opT y_{N-1}\|^2$ for fixed-point problem is equivalent to finding an algorithm that effectively reduces $\| \tilde{\opA} x_N \|^2$ for monotone inclusions.

\section{Exact optimal methods}
\label{sec:method}
We now present our methods and their accelerated rates.

For a $1/\gamma$-contractive operator $\opT \colon \reals^n\to\reals^n$, the \emph{Optimal Contractive Halpern} \eqref{algo:oc-halpern} is
\begin{align*}
    y_k &= \left( 1 - \frac{1}{\varphi_k} \right) \opT y_{k-1} + \frac{1}{\varphi_k} y_0
    \tag{OC-Halpern}
    \label{algo:oc-halpern}
\end{align*}
for $k=1,2,\dots$, where $\varphi_k=\sum_{i=0}^{k} \gamma^{2i}$ and $y_0\in \mathbb{R}^n$ is a starting point.
For a maximal $\mu$-strongly monotone operator $\opA \colon \reals^n\rightrightarrows\reals^n$, the \emph{Optimal Strongly-monotone Proximal Point Method} \eqref{algo:os-ppm} is
\begin{align*}
    x_k &= \resA y_{k-1}
    \tag{OS-PPM}
    \label{algo:os-ppm}\\
    y_k 
    &= x_k + \frac{\varphi_{k-1}-1}{\varphi_k} (x_k-x_{k-1}) - \frac{2\mu \varphi_{k-1}}{\varphi_k} (y_{k-1}-x_k) \\
    &~\qquad + \frac{(1+2\mu) \varphi_{k-2}}{\varphi_k} (y_{k-2}-x_{k-1})
\end{align*}
for $k=1,2,\dots$,
where $\varphi_k = \sum_{i=0}^k (1+2\mu)^{2i}$, $\varphi_{-1} = 0$, and $x_0 = y_0 = y_{-1} \in \mathbb{R}^n$ is a starting point.
These two methods are equivalent.

\begin{lemma} \label{equiv:method}
Suppose $\gamma = 1+2\mu$. 
Let $\opA = \left( \opT + \frac{1}{\gamma} \opI \right)^{-1} \left( 1 + \frac{1}{\gamma} \right) - \opI$ given $\opT$, or equivalently let $\opT = \left( 1+ \frac{1}{1+2\mu} \right) \resA - \frac{1}{1+2\mu} \opI$ given $\opA$. Then the $y_k$-iterates of \eqref{algo:oc-halpern} and \eqref{algo:os-ppm} are identical provided they start from the same initial point $y_0 = \tilde{y}_0$.
\end{lemma}

We now state the convergence rates.

\begin{theorem}
\label{lyapunov:os-ppm}
    Let $\opA\colon\reals^n\rightrightarrows\reals^n$ be maximal $\mu$-strongly monotone with $\mu \ge 0$.
    Assume $\opA$ has a zero and let $x_\star\in\zer\opA$.
    For $N=1,2,\dots$, \eqref{algo:os-ppm} exhibits the rate
    \begin{align*}
        \| \tilde{\opA} x_N \|^2 
        &\le
        \left( \frac{1}{\sum_{k=0}^{N-1} (1+2\mu)^k} \right)^2 \| y_0-x_\star \|^2.
    \end{align*}
\end{theorem}
\begin{corollary}
\label{lyapunov:oc-halpern}
    Let $\opT\colon\reals\to\reals$ be $\gamma^{-1}$-contractive with $\gamma\ge1$.
    Assume $\opT$ has a fixed point and let $y_\star\in\fix\opT$.
    For $N=0,1,\dots$, \eqref{algo:oc-halpern} exhibits the rate
    \[
        \| y_{N} - \opT y_{N} \|^2
        \le
        \left( 1+ \frac{1}{\gamma} \right)^2
        \left( \frac{1}{\sum_{k=0}^{N} \gamma^k} \right)^2 \| y_0-y_\star \|^2.
    \]
\end{corollary}

When $\opA$ is strongly monotone ($\mu>0$), \eqref{algo:os-ppm} exhibits an accelerated $\cO(e^{-4\mu N})$-rate compared to the $\cO(e^{-2\mu N})$-rate of the proximal point method (PPM) \cite{rockafellar1976monotone,bauschke2011convex}.
When $\opT$ is contractive ($\gamma<1$), both \eqref{algo:oc-halpern} and the Picard iteration exhibit $\cO(\gamma^{-2N})$-rates on the squared fixed-point residual.
In fact, the Picard iteration with the $\opT$ of Lemma~\ref{equiv:sm-lips} instead of $\resA $ is faster than the regular PPM and achieves a $\cO(e^{-4\mu N})$ rate.
\eqref{algo:oc-halpern} is exactly optimal and is faster than Picard in higher order terms hidden in the big-$\mathcal{O}$ notation. To clarify, the $\mathcal{O}$ considers the regime $\mu\to0$.

When $\opA$ is not strongly monotone ($\mu=0$) or $\opT$ is not contractive ($\gamma=1$), \eqref{algo:os-ppm} and \eqref{algo:oc-halpern} respectively reduces to accelerated PPM (APPM) of \citet{kim2021accelerated} and Halpern iteration of \citet{Lieder2021halpern}, sharing the same $\cO(1/N^2)$-rate.
In this paper, we refer to the method of \citet{Lieder2021halpern} as the optimized Halpern method (OHM).

The discovery of \eqref{algo:oc-halpern} and \eqref{algo:os-ppm} was assisted by the performance estimation problem \cite{drori2014performance,kim2016optimized,taylor2017smooth,drori2020efficient,ryu2020operator,kim2021optimizing,park2021optimal}
The details are discussed in Section~\ref{detail:pep} of the appendix.

\subsection{Proof outline of \cref{lyapunov:os-ppm}}
Here, we quickly outline the proof of \cref{lyapunov:os-ppm} while deferring the full proof to Section~\ref{proof:method} of the appendix.

Define the Lyapunov function
\begin{align*}
    &V^k = (1+\gamma^{-k}) \Bigg[ 
        \left(\sum_{n=0}^{k-1} \gamma^n\right)^2 \|\tilde{\opA}x_k\|^2 \\
    &
        + 2\left(\sum_{n=0}^{k-1}\gamma^n\right) \langle \tilde{\opA}x_k - \mu(x_k-x_\star), x_k-x_\star\rangle \\
    &
        + \gamma^{-k} \left\| \left(\sum_{n=0}^{k-1}\gamma^n\right) \tilde{\opA}x_k - \gamma^k (x_k-x_\star) + (x_k-y_0) \right\|^2
    \Bigg] \\
    &\qquad
        + (1-\gamma^{-k})\|y_0-x_\star\|^2
    \tag{OS-PPM-Lyapunov}
    \label{lyapunov:potential}
\end{align*}
for $k=0,1,\dots$, where $\gamma=1+2\mu$ and $\tilde{\opA}x_k = y_{k-1} - x_k \in \opA x_k$.
After some calculations (deferred to the appendix), we use $\mu$-strong monotonicity of $\opA$ to conclude
\begin{align*}
    V^{k+1} - V^{k} &= 
    - 2\gamma^{-2k}(1+\gamma)\varphi_k \varphi_{k-1} \\ &\langle \tilde{\opA}x_{k+1} - \tilde{\opA}x_k - \mu(x_{k+1}-x_k), x_{k+1}-x_k \rangle \\
    &\le 0.
\end{align*}
Therefore, 
\[
    V^N \le V^{N-1} \le \dots \le V^0 = 2\|y_0-x_\star\|^2
\]
and we conclude
\[
    \| \tilde{\opA} x_N \|^2 \le \left( \frac{1}{\sum_{k=0}^{N-1} \gamma^k} \right)^2 \| y_0 - x_\star \|^2.
\]

\section{Complexity lower bound}
\label{sec:lower-bound}
We now establish \emph{exact} optimality of  \eqref{algo:oc-halpern} and \eqref{algo:os-ppm} through matching complexity lower bound.
By exact, we mean that the lower bound is exactly equal to upper bounds of \cref{lyapunov:os-ppm} and \cref{lyapunov:oc-halpern}.

\begin{theorem}
    \label{lower-bound:oc-halpern}
    For $n \ge N+1$ and any initial point $y_0\in\reals^n$, there exists an $1/\gamma$-Lipschitz operator $\opT\colon\reals^n\rightarrow \reals^n$ with a fixed point $y_\star\in\fix\opT$ such that 
    \[
        \| y_{N} - \opT y_{N} \|^2 \ge \left( 1 + \frac{1}{\gamma} \right)^2 \left( \frac{1}{\sum_{k=0}^{N} \gamma^k} \right)^2 \| y_0 - y_\star \|^2
    \]
    for any iterates $\{y_k\}_{k=0}^N$ satisfying 
    {    \color{black}
    \[
        y_k \in y_0 + \Span \{y_0 - \opT y_0, y_1 - \opT y_1, \dots, y_{k-1} - \opT y_{k-1}\}
    \]}
    for $k=1,\dots,N$.
\end{theorem}
The following corollary translates \cref{lower-bound:oc-halpern} to an equivalent complexity lower bound for proximal point methods in monotone inclusions.
\begin{corollary}
    \label{lower-bound:os-ppm}
    For $n \ge N$ and any initial point $x_0 = y_0\in\reals^n$, there exists a maximal $\mu$-strongly monotone operator $\opA\colon\reals^n\rightrightarrows\reals^n$ with a zero $x_\star\in\zer\opA$ such that
    \[
        \| \tilde{\opA} x_N \|^2 \ge \left( \frac{1}{\sum_{k=0}^{N-1} (1+2\mu)^k} \right)^2 \| y_0 - x_\star \|^2
    \]
    for any iterates $\{x_k\}_{k=0,1,\dots}$ and $\{y_k\}_{k=0,1,\dots}$ satisfying
    \begin{align*}
        x_k &= \resA y_{k-1} \\
        y_k &\in y_0 + \Span \{\tilde{\opA} x_1, \tilde{\opA} x_2, \dots, \tilde{\opA} x_{k}\}
    \end{align*}
    for $k=1,\dots,N$, where $\tilde{\opA}x_k = y_{k-1}-x_k$.
\end{corollary}

\eqref{algo:oc-halpern} and \eqref{algo:os-ppm} satisfy the span assumptions stated in \cref{lower-bound:oc-halpern} and \cref{lower-bound:os-ppm}, respectively.
Therefore, the rates of \eqref{algo:oc-halpern} and \eqref{algo:os-ppm} are exactly optimal.
The lower bounds in the cases where $\gamma=1$ and $\mu=0$ establish that the prior rates of OHM \cite{Lieder2021halpern} and APPM \cite{kim2021accelerated} are exactly optimal.
To clarify, the lower bound is novel even for the case $\gamma=1$ and $\mu=0$.

\subsection{Construction of the worst-case operator}
We now describe the construction of the worst-case operator, while deferring the proofs to Section~\ref{proof:lower-bound} of the appendix.
Let $e_k$ be the canonical basis vector with $1$ at the $k$\nobreakdash-th entry and $0$ at remaining entries.

\begin{lemma} \label{lem:equiv-ops}
    $\opT$ is $\frac{1}{\gamma}$-contractive if and only if $\opG = \frac{\gamma}{1+\gamma}(\opI-\opT)$ is $\frac{1}{1+\gamma}$-averaged.
\end{lemma}

By \cref{lem:equiv-ops}, finding the worst-case $\tfrac{1}{\gamma}$-contractive operator $\opT$ is equivalent to finding the worst-case $\tfrac{1}{1+\gamma}$-averaged operator $\opG$, which we define in the following lemma.

\begingroup
\allowdisplaybreaks
\begin{lemma} \label{lem:T-aver}
    Let $R > 0$.
    Define  $\opN, \opG \colon \reals^{N+1} \to \reals^{N+1}$ as
    \begin{align*}
        \opN(x_1, x_2, \dots, x_{N}, x_{N+1}) 
        &= 
        (x_{N+1}, -x_1, -x_2, \dots, -x_{N}) \\
        &\quad - \frac{1+\gamma^{N+1}}{\sqrt{1+\gamma^2+\dots+\gamma^{2N}}} R e_1
    \end{align*}
    and
    \[
        \opG = \frac{1}{1+\gamma} \opN + \frac{\gamma}{1+\gamma} \opI.
    \]
    That is,
    \begin{align*}
        \opG x
        &=
        \underbrace{
            \frac{1}{1+\gamma}
            \begin{bmatrix}
                \gamma & 0 & \cdots & 0 & 1 \\
                -1 & \gamma & \cdots & 0 & 0 \\
                \vdots & \vdots & \ddots & \vdots & \vdots \\
                0 & 0 & \cdots & \gamma & 0 \\
                0 & 0 & \cdots & -1 & \gamma 
            \end{bmatrix}
        }_{=: H}
        x \\
        &\quad - 
        \underbrace{
            \frac{1}{1+\gamma} \frac{1+\gamma^{N+1}}{\sqrt{1+\gamma^2+\dots+\gamma^{2N}}} Re_1
        }_{=: b}.
    \end{align*}
    Then $\opN$ is nonexpansive, and $\opG$ is $\frac{1}{1+\gamma}$-averaged.
\end{lemma}
\endgroup

Following lemma states the property of iterations $\{y_k\}_{k=0}^N$ with respect to $\opG$, that proper span condition results in gradually expanding support of $y_k$.

\begin{lemma} \label{lem:span}
    Let $\opG \colon \reals^{N+1} \to \reals^{N+1}$ be defined as in \cref{lem:T-aver}.
    For any $\{y_k\}_{k=0}^N$ with $y_0 = 0$ satisfying
    \[
        y_k \in y_0 + \Span \{ \opG y_0, \opG y_1, \dots, \opG y_{k-1} \}, \quad k=1,\dots,N,
    \]
    we have
    \begin{align*}
        y_k 
        &\in \Span \left\{ e_1, e_2, \dots, e_{k} \right\} \\
        \opG y_k
        &\in \Span \left\{ e_1, e_2, \dots, e_{k+1} \right\},
        \quad k=0,\dots,N.
    \end{align*}
\end{lemma}

\subsection{Proof outline of \cref{lower-bound:oc-halpern}}
    \label{subsec:lower-bound}
    Let $\opT_0\colon\reals^n\to\reals^n$ be the worst-case $\frac{1}{\gamma}$-contraction for initial point $0$.
    For any given $y_0\in\reals^n$, we show in section \ref{proof:lower-bound} of the appendix that $\opT\colon\reals^n\to\reals^n$ defined as $\opT (\cdot) = \opT_0(\cdot-y_0) + y_0$ becomes the worst-case $\frac{1}{\gamma}$-contraction with initial point $y_0\in\reals^n$.
    Therefore, it suffices to consider the case $y_0 = 0$.
    
    Define $\opG$, $H$, and $b$ as in \cref{lem:T-aver}. By \cref{lem:equiv-ops}, $\opT = \opI - \frac{1+\gamma}{\gamma} \opG$ is a $1/\gamma$-contraction.
    Note that $H$ is invertible, as we can use Gaussian elimination on $H$ to obtain an upper triangular matrix with nonzero diagonals.
    This makes $\opG$ an invertible affine operator with the unique zero
    \[
        y_\star = 
        \displaystyle \frac{R}{\sqrt{1+\gamma^2 + \dots + \gamma^{2N}}} \begin{bmatrix} \gamma^{N} & \gamma^{N-1} & \cdots & \gamma & 1 \end{bmatrix}^\intercal.
    \]
So $\fix\opT = \zer\opG = \{y_\star\}$ and $\| y_0 - y_\star \|=\| y_\star\| = R$.
    
    Let the iterates $\{y_k\}_{k=0}^N$ satisfy the span condition of \cref{lower-bound:oc-halpern}, which is equivalent to
    \[
        y_k \in y_0 + \Span \{ \opG y_0, \opG y_1, \dots, \opG y_{k-1} \},
        \quad k=1,\dots,N.
    \]
    By \cref{lem:span}, $y_{N} \in \Span \{ e_1, \dots, e_{N} \}$.
    Therefoere
    \[
        \opG y_{N} = H y_{N} - b \in \Span \{ H e_1, \dots, H e_{N} \} - b.
    \]
    and
    \[
        \left\| \opG y_{N} \right\|^2
        \ge 
        \left\| \cP_{{\Span \{ H e_1, \dots, H e_{N} \}}^\perp} (b) \right\|^2,
    \]
    where $\cP_V$ is the orthogonal projection onto the subspace $V$.
    As ${\Span \{ H e_1, \dots, H e_{N} \}}^\perp = \Span \{ v \}$ with
    \[
        v=\begin{bmatrix} 1 & \gamma & \cdots & \gamma^{N-1} & \gamma^{N} \end{bmatrix}^\intercal,
    \]
    we get
    \begin{align*}\!\!
        \| \opG y_{N} \|^2
        &\!\ge\!
        \left\| \cP_{\Span\{v\}} (b) \right\|^2\!
        \!=\!
        \left\|
            \frac{\langle b, v \rangle}{\langle v, v \rangle} v
        \right\|^2 \!\!\stackrel{(*)}{=}\!
        \left( \frac{1}{\sum_{k=0}^{N} \gamma^k} \right)^2 \!\!\!R^2,
    \end{align*}
    where $(*)$ is established in the Section~\ref{proof:lower-bound} of the appendix.
    Finally,
    \begin{align*}
        \| y_{N} - \opT y_{N} \|^2
        &=
        \left\| \left(1+\frac{1}{\gamma}\right) \opG y_{N} \right\|^2 \\
        &\ge
        \left(1+\frac{1}{\gamma}\right)^2 \left( \frac{1}{\sum_{k=0}^{N} \gamma^k} \right)^2 R^2 \\
        &=
        \left(1+\frac{1}{\gamma}\right)^2 \left( \frac{1}{\sum_{k=0}^{N} \gamma^k} \right)^2 \|y_0-y_\star\|^2.
    \end{align*}

\subsection{Generalized complexity lower bound result}
In order to extend the lower bound results of \cref{lower-bound:oc-halpern} and \cref{lower-bound:os-ppm} to general deterministic fixed-point iterations and proximal point methods (which do not necessarily satisfy the span condition), we use the resisting oracle technique of \citet{nemirovski1983problem}.
Here, we quickly state the result while deferring the proofs to the Section~\ref{proof:lower-bound} of the appendix.

\begin{theorem}
    \label{lower-bound:general-halpern}
    Let $n\ge 2N$ for $N\in\mathbb{N}$.
    For any deterministic fixed-point iteration $\al$
    and any initial point $y_0\in\reals^n$, there exists a $\frac{1}{\gamma}$-Lipschitz operator $\opT \colon \reals^n \to \reals^n$ with a fixed point $y_\star\in\fix\opT$ such that
    \[
        \|y_N - \opT y_N\|^2 \ge \left( 1 + \frac{1}{\gamma} \right)^2
        \left( \frac{1}{\sum_{k=0}^N \gamma^k} \right)^2 \|y_0-y_\star\|^2
    \]
    where $\{y_t\}_{t\in\mathbb{N}} = \al[y_0; \opT]$.
\end{theorem}

\begin{figure*}[ht]
    \centering
    \subfigure[Fixed-point residual of $\opT_\theta$]{
        \centering
        \includegraphics[width=0.8\columnwidth]{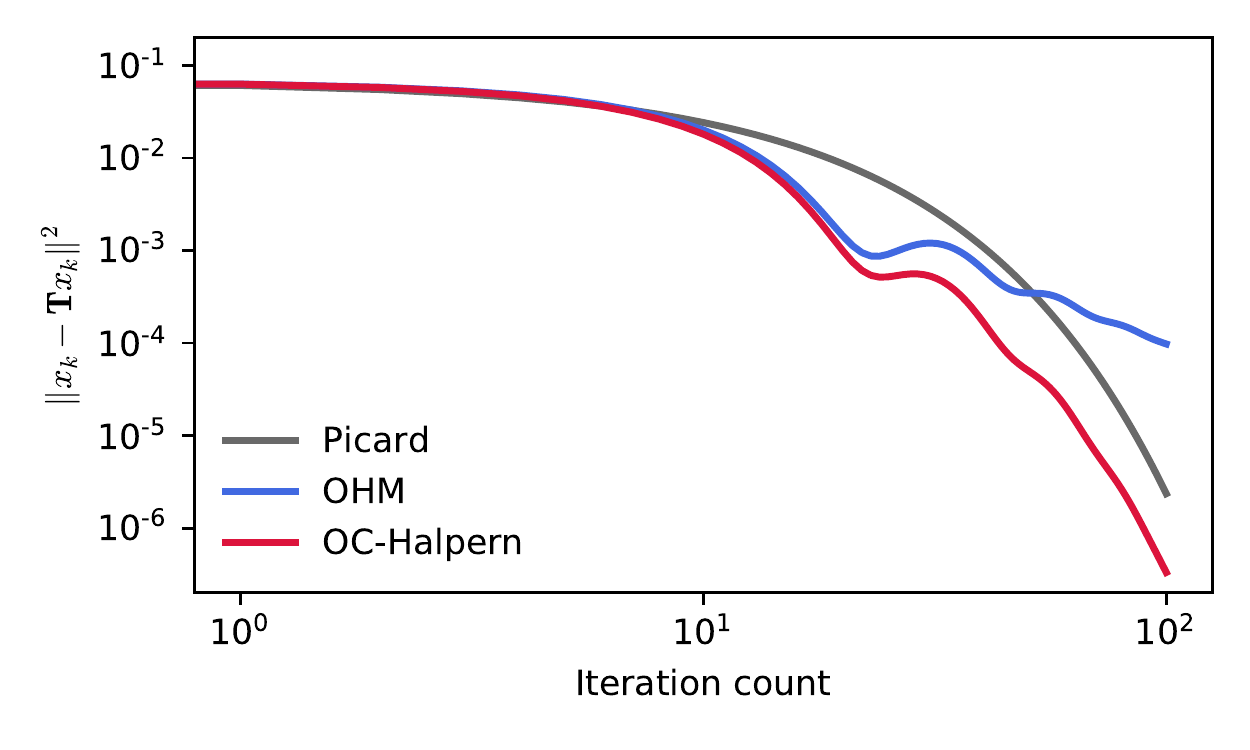}
        \label{fig:oc-halpern-fpr}
    }
    \hspace{.05\textwidth}
    \subfigure[Resolvent residual norm of $\opM$]{
        \centering
        \includegraphics[width=0.8\columnwidth]{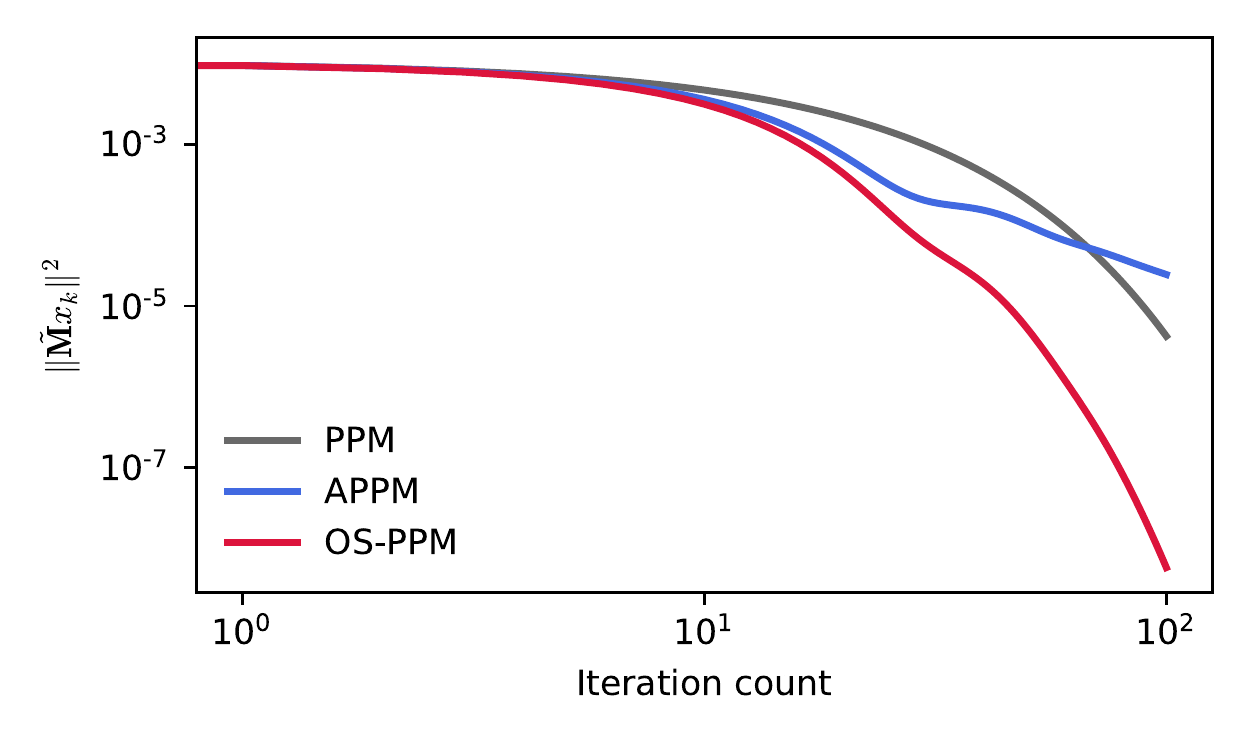}
        \label{fig:os-ppm-fpr}
    }
    \caption{Fixed-point and resolvent residuals versus iteration count for the 2D toy example of \cref{exper:toy}.
    Here, $\gamma=1/0.95=1.0526$, $\mu=0.035$, $\theta=15^\circ$ and $N=101$.
    Indeed, \eqref{algo:oc-halpern} and \eqref{algo:os-ppm} exhibit the fastest rates.}
        \label{fig:fig1}
\end{figure*}

\begin{figure*}[ht]
    \centering
    \subfigure[Trajectory of $\opT_\theta$]{
        \centering
        \includegraphics[width=0.75\columnwidth]{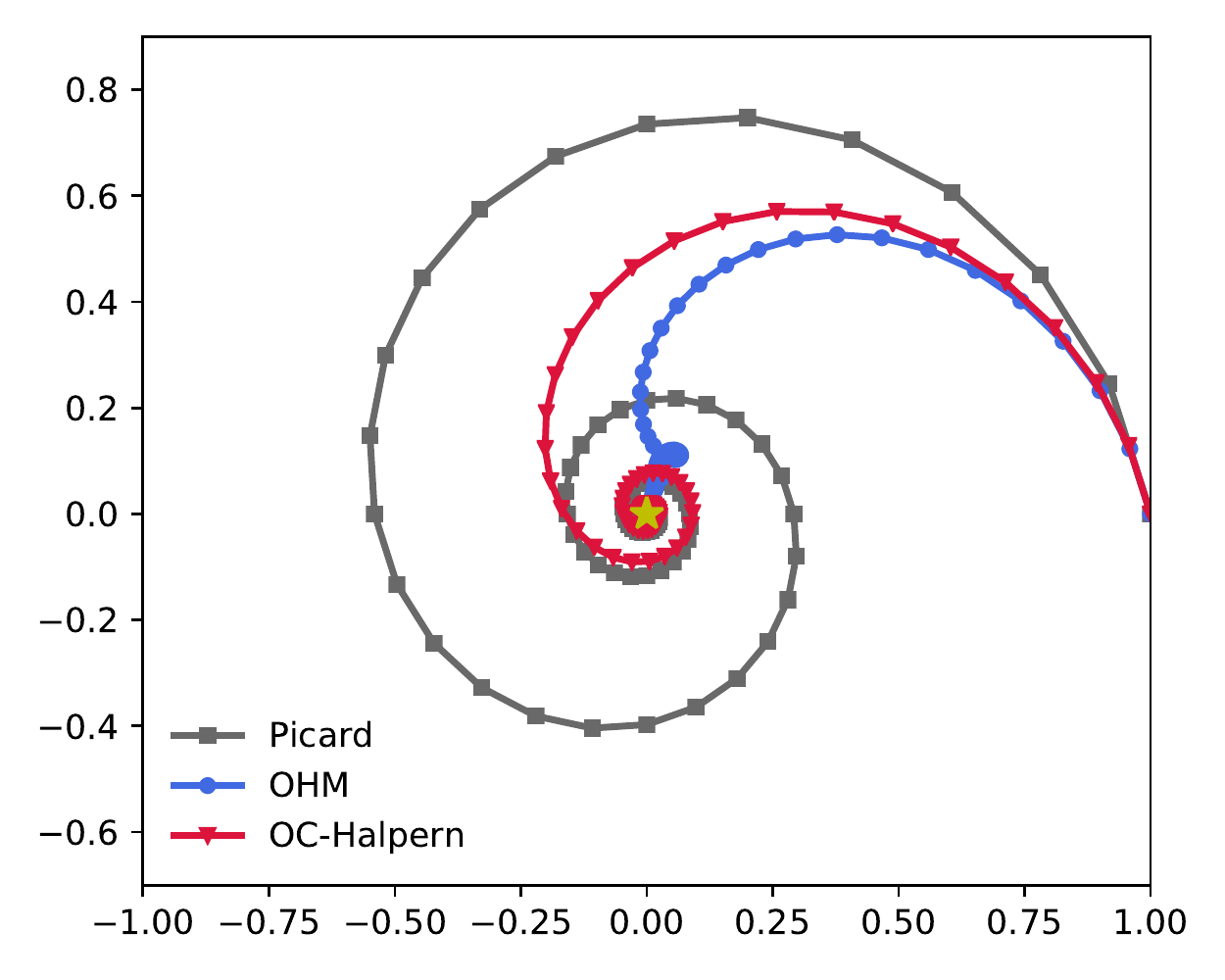}
        \label{fig:oc-halpern-trajectory}
    }
    \hspace{.1\textwidth}
    \subfigure[Trajectory of $\opM$]{
        \centering
        \includegraphics[width=0.75\columnwidth]{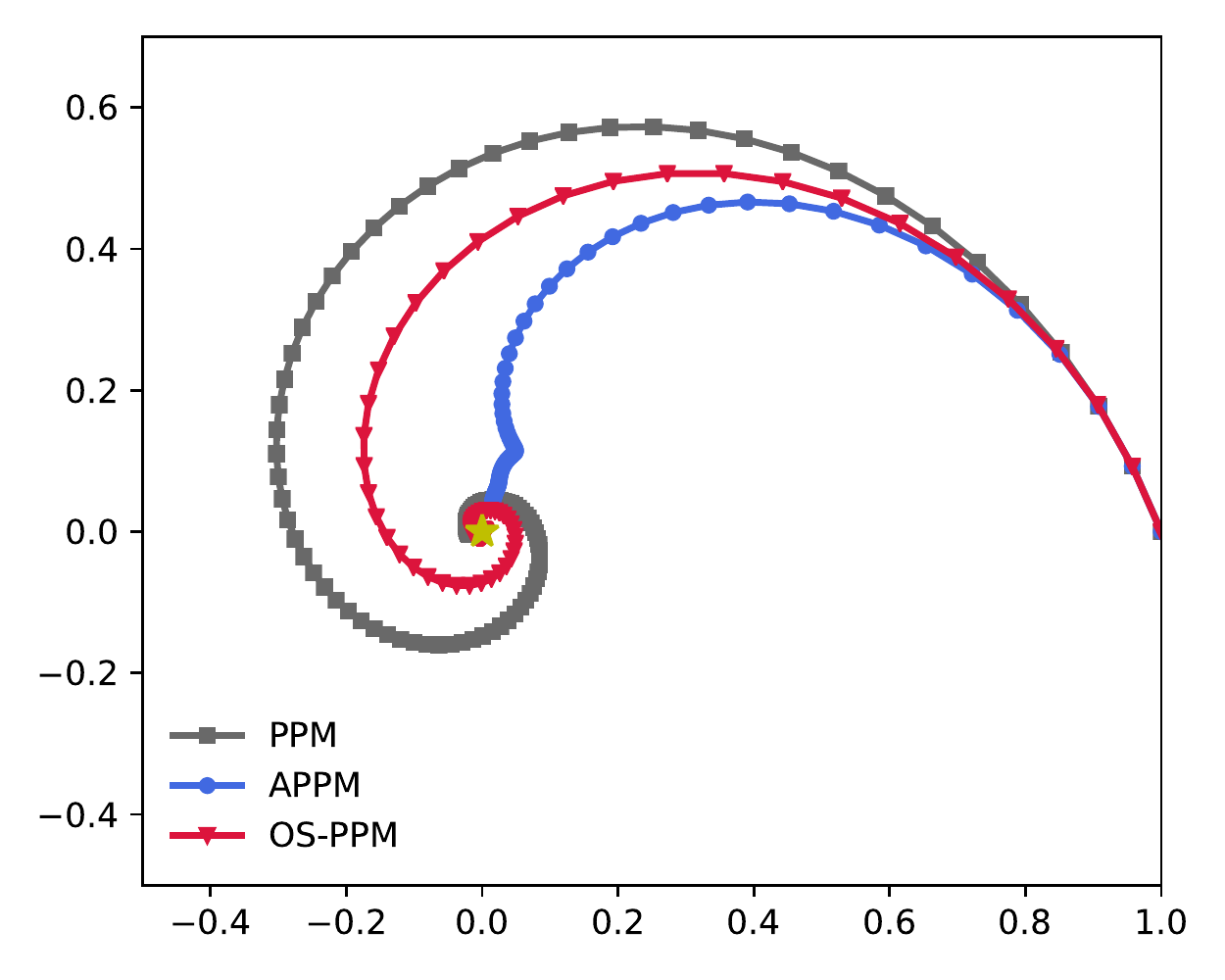}
        \label{fig:os-ppm-trajectory}
    }
    \caption{Trajectories of iterates for the 2D toy example of \cref{exper:toy}. 
    Here, $\gamma=1/0.95=1.0526$, $\mu=0.035$, $\theta=15^\circ$ and $N=101$.
    A marker is placed at every iterate.
    Picard and PPM are slowed down by the cyclic behavior. Halpern and APPM dampens the cycling behavior, but does so too aggressively.
    The fastest rate is achieved by \eqref{algo:oc-halpern} and \eqref{algo:os-ppm}, which appears to be due to the cycling behavior being optimally dampened.
    }
\end{figure*}

\section{Acceleration under H\"older-type growth condition} \label{sec:restart}

While \eqref{algo:os-ppm} provides an accelerated rate when the underlying operator is monotone or strongly monotone, many operators encountered in practice have a structure lying between these two assumptions.
For \eqref{algo:oc-halpern}, this corresponds to a fixed-point operator that is not strictly contractive but has structure stronger than nonexpansiveness.
In this section, we accelerate the proximal point method when the underlying operator is uniformly monotone, an assumption weaker than strong monotonicity but stronger than monotonicity.

We say an operator $\opA \colon \reals^n \rightrightarrows \reals^n$ is \emph{uniformly monotone with parameters $\mu > 0$ and $\alpha > 1$} if it is monotone and
\[
    \langle \opA x, x - x_\star \rangle \ge \mu \|x - x_\star\|^{\alpha+1}
\]
for any $x \in \reals^n$ and $x_\star\in\zer\opA$.
This is a special case of uniform monotonicity in \citet[Definition~22.1]{bauschke2011convex}.
We also refer to this as a H\"older-type growth condition, as it resembles the H\"olderian error bound condition with function-value suboptimality replaced by $\langle \opA x, x - x_\star \rangle$ \cite{lojasiewicz1963propriete,bolte2017error}.


The following theorem establishes a convergence rate of the (unaccelerated) proximal point method.
This rate serves as a baseline to improve upon with acceleration.
\begin{theorem} \label{sharp:ppm}
    Let $\opA \colon \reals^n \rightrightarrows \reals^n$ be uniformly monotone with parameters $\mu > 0$ and $\alpha > 1$. Let $x_\star\in \zer\opA$. 
    Then the iterates $\{x_k\}_{k=0}^N$ generated by the proximal point method $x_{k+1} = \resA x_k$ starting from $x_0 \in \reals^n$ satisfy
    \begin{align*}
        \|\tilde{\opA} x_N\|^2 
        &\le 
        \frac{ 2^\frac{\alpha+3}{\alpha-1} \max\left\{ \left(\frac{2^\frac{\alpha}{\alpha-1}-2}{\mu}\right)^\frac{2}{\alpha-1}, \|x_0-x_\star\|^2 \right\} }{ N^\frac{\alpha+1}{\alpha-1} }
    \end{align*}
    for $N \in \mathbb{N}$ where $\tilde{\opA}x_N = x_{N-1} - x_N$.
\end{theorem}

We now present an accelerated method based on \eqref{algo:os-ppm} and restarting \cite{nesterov2013gradient,roulet2020sharpness}.
Given a uniformly monotone operator $\opA\colon\reals^n\rightrightarrows\reals^n$ with $\mu>0$ and $\alpha>1$, $x_\star\in\zer\opA$, and an initial point $x_0\in\reals^n$, 
\emph{Restarted OS-PPM} is:
\begin{align*}
    \tilde{x}_0 &= \resA x_0
    \tag{OS-PPM$_0^\mathrm{res}$}
    \label{algo:os-ppm-res} \\
    \tilde{x}_k &\gets \textbf{OS-PPM}_0 (\tilde{x}_{k-1}, t_k), \qquad k=1,\dots,R,
\end{align*}
where $\textbf{OS-PPM}_0(\tilde{x}_{k-1}, t_k)$ is the execution of $t_k$ iterations of \eqref{algo:os-ppm} with $\mu=0$ starting from $\tilde{x}_{k-1}$.
The following theorem provides a restarting schedule, i.e., specified values of $t_1,\dots,t_R$, and an accelerated rate.

\newpage
\begin{theorem} \label{sharp:appm}
    Let $\opA \colon \reals^n \rightrightarrows \reals^n$ be uniformly monotone with parameters $\mu > 0$ and $\alpha > 1$, $x_\star\in\zer\opA$, and $N$ be the total number of iterations.
    Define
    \[
        \lambda 
        = 
        \left(\frac{e}{\mu}\right)^{\frac{1}{\alpha}} \|x_0-x_\star\|^{1-\frac{1}{\alpha}},
        \qquad
        \beta 
        = 
        1-{\frac{1}{\alpha}}.
    \]
    Let $R\in \mathbb{N}$ be an integer satisfying
    \[
        \sum_{k=1}^{R} \lceil \lambda e^{\beta k} \rceil \le N-1 < \sum_{k=1}^{R+1} \lceil \lambda e^{\beta k} \rceil,
    \]
    and let $t_k$ be defined as
    \[
        t_k = \begin{cases}
            \left\lceil \lambda e^{\beta k} \right\rceil & \text{ for }k=1,\dots,R-1\\
            N - 1 - \sum_{k=1}^{R-1} t_k. & \text{ for }k=R.
        \end{cases}
    \]
    Then \eqref{algo:os-ppm-res} exhibits the rate
    \begin{align*}
        &\|\tilde{\opA} x_N\|^2 \\
        &\le 
        \left\{
            {\scriptstyle\frac{e^\beta-1}{\lambda e^{2\beta}}} \left(N-2-{\scriptstyle\frac{1}{\beta}} 
            \log\left({\scriptstyle \frac{e^\beta-1}{\lambda e^\beta}}(N-1) + 1 \right) \right) + {\scriptstyle\frac{1}{e^\beta}}
        \right\}^{-\frac{2\alpha}{\alpha-1}} \\
        &\qquad \times \|x_0-x_\star\|^2 \\
        &=
        \mathcal{O} \left( N^{- \frac{2\alpha}{\alpha-1}} \right).
    \end{align*}
\end{theorem}
The proofs of Theorems~\ref{sharp:ppm} and \ref{sharp:appm} are presented in Section~\ref{proof:restart} of the appendix.
When the values of $\alpha$, $\mu$, and $\|x_0-x_\star\|^2$ are unknown, as in the case in most practical setups, one can use a grid search as in \citet{roulet2020sharpness} and retain the $\cO\left(N^{-\frac{2\alpha}{\alpha-1}}(\log N)^2\right)$-rate.
Using \cref{lem:equiv-ops}, \eqref{algo:os-ppm-res} can be translated into a restarted OC-Halpern method. 
The experiments of Section \ref{sec:experiment} indicate that \eqref{algo:os-ppm-res} does provide an acceleration in cases where \eqref{algo:os-ppm} by itself does not.

\section{Experiments} \label{sec:experiment}

\begin{figure*}[ht]
    \centering
    \subfigure[CT imaging]{
        \includegraphics[width=.65\columnwidth]{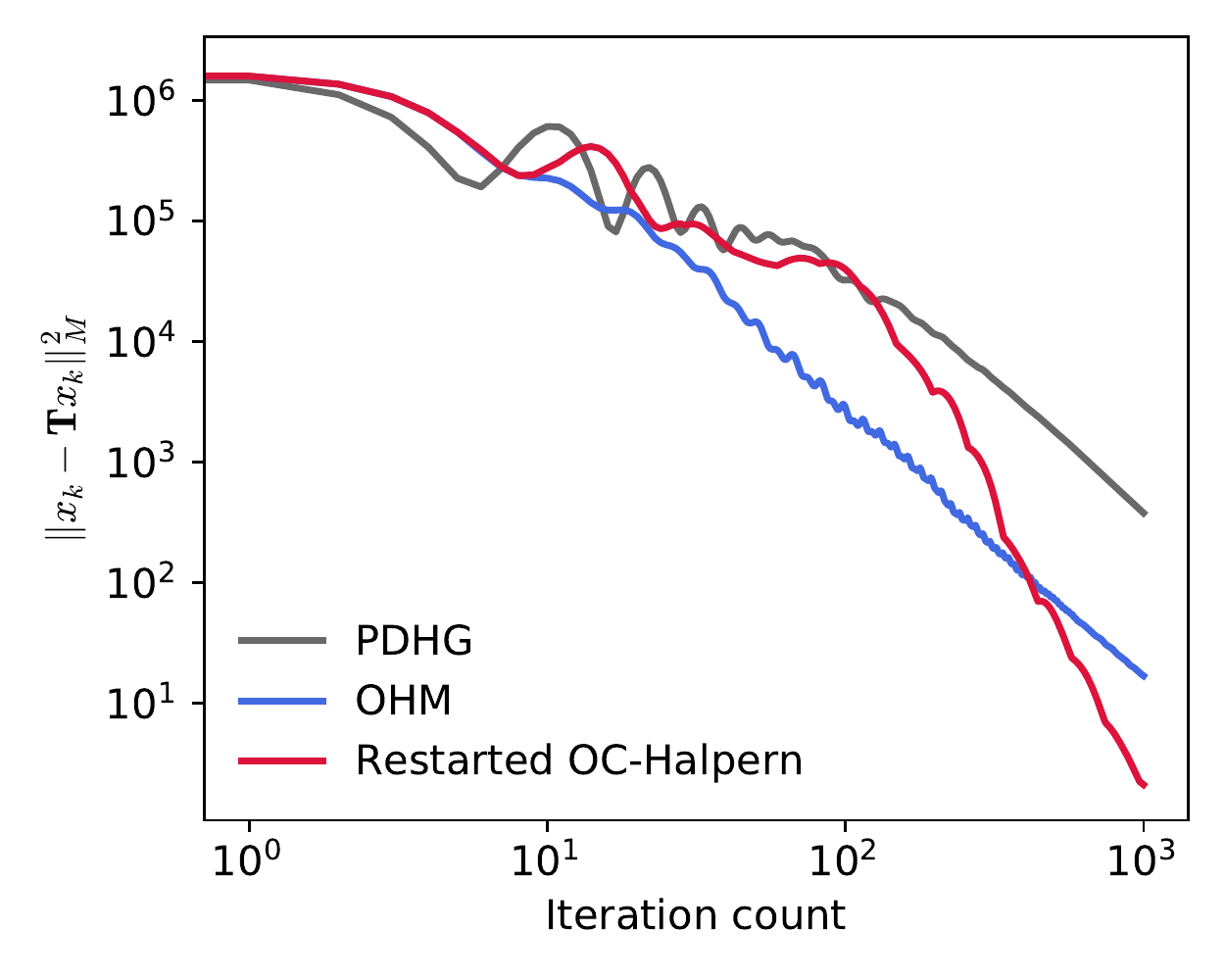}
        \label{fig:ct-fpr}
    }
    \hfill
    \subfigure[Earth mover's distance]{
        \includegraphics[width=.65\columnwidth]{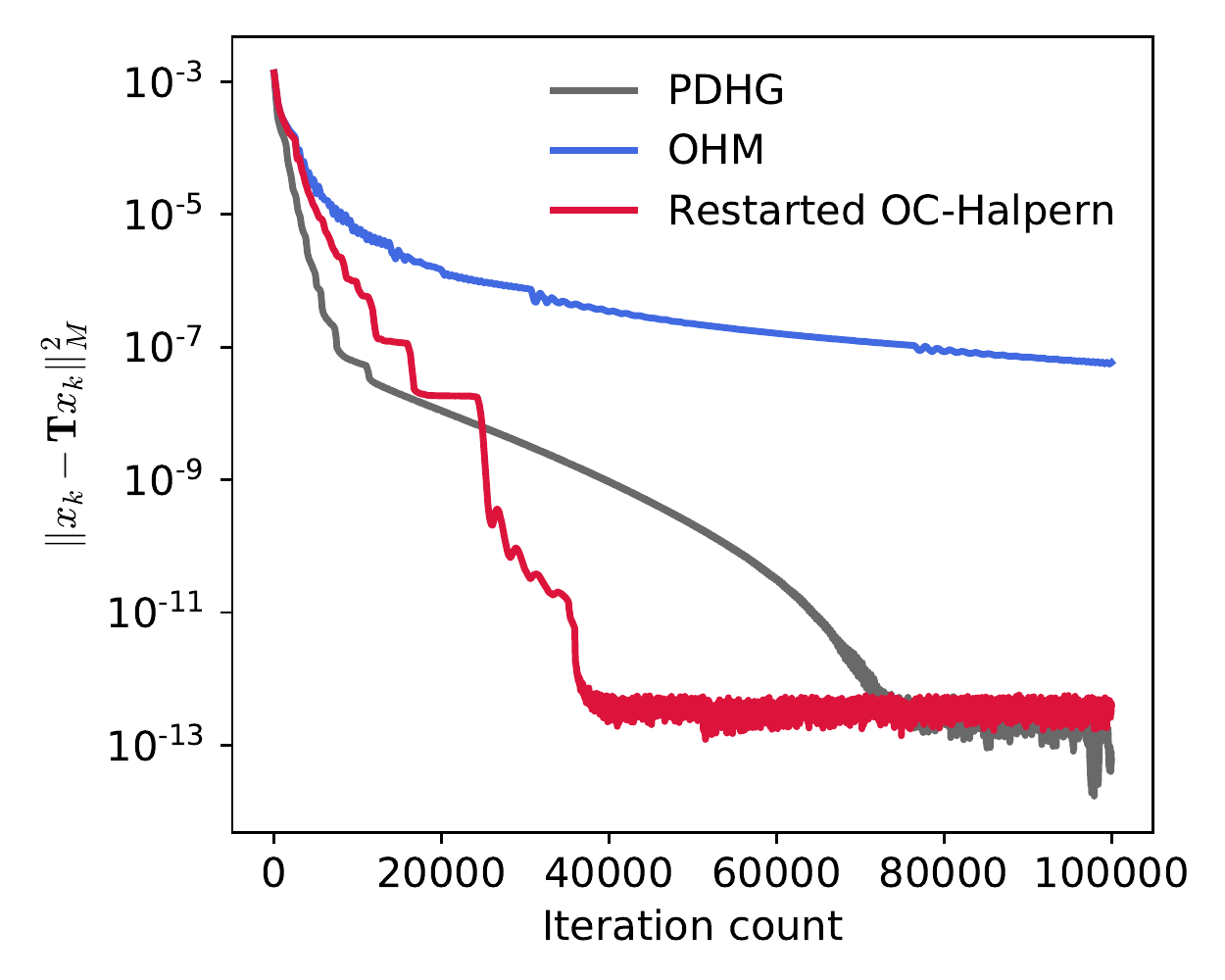}
        \label{fig:emd-fpr}
    }
    \hfill
    \subfigure[Decentralized compressed sensing]{
        \includegraphics[width=.65\columnwidth]{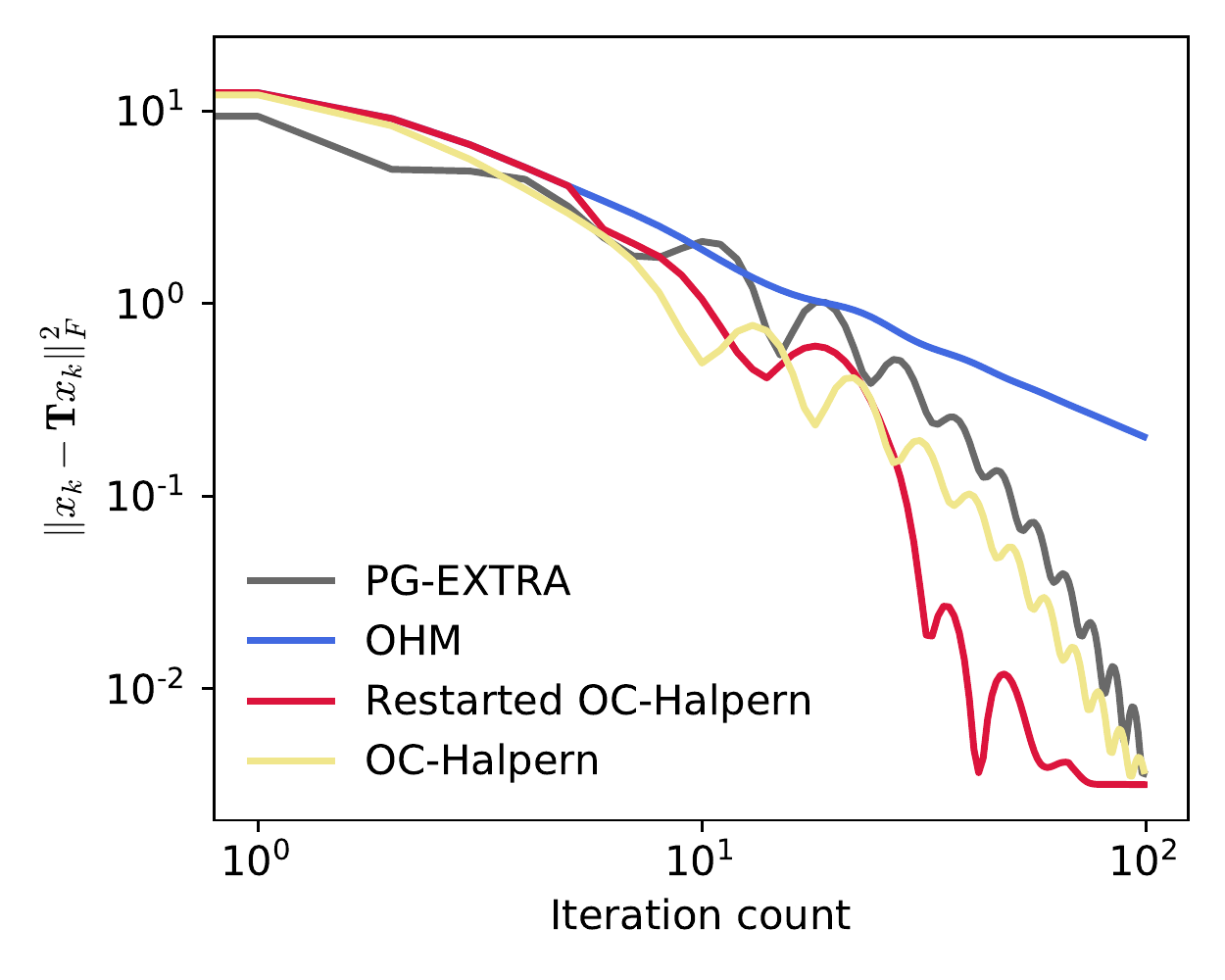}
        \label{fig:decen-fpr}
    }
    \vspace{-0.05in}
    \caption{Reduction of fixed-point residuals. The norms are the norms under which the fixed-point operator $\opT$ is nonexpansive.}
    \vspace{-0.1in}
\end{figure*}

We now present experiments with illustrative toy examples and real-world problems in medical imaging, optimal transport, and decentralized compressed sensing.
Further experimental details are provided in Section~\ref{sec:expr-details} of the appendix.

\subsection{Illustrative 2D toy examples} \label{exper:toy}

Consider a $\frac{1}{\gamma}$-contractive operator $\opT_\theta \colon \reals^2 \to \reals^2$
\[
    \opT_\theta \begin{bmatrix} x_1 \\ x_2 \end{bmatrix} 
    =
    \frac{1}{\gamma}
    \begin{bmatrix} \cos \theta & - \sin \theta \\ \sin \theta & \cos \theta \end{bmatrix} \begin{bmatrix} x_1 \\ x_2 \end{bmatrix}
\]

and a maximal $\mu$-strongly monotone operator $\opM \colon \reals^2 \to \reals^2$
\[
    \opM \begin{bmatrix} x_1 \\ x_2 \end{bmatrix}
    =
    \left( \frac{1}{N-1} \begin{bmatrix} 0 & 1 \\ -1 & 0 \end{bmatrix} + \begin{bmatrix} \mu & 0 \\ 0 & \mu \end{bmatrix} \right)
    \begin{bmatrix} x_1 \\ x_2 \end{bmatrix}.
\]
$\opT_\theta$ is a counterclockwise $\theta$-rotation followed by $\frac{1}{\gamma}$-scaling on 2D plane, and $\opM$ is a linear combination of the worst-case instances of the proximal point method applied to monotone operators \cite{gu2020tight} and $\mu$-strongly monotone operators \cite{rockafellar1976monotone}.
The results of \cref{fig:fig1} indicate that \eqref{algo:oc-halpern} and \eqref{algo:os-ppm} indeed provide acceleration.

\subsection{Computed tomography (CT) imaging} \label{exper:ct}
Consider the medical imaging application of total variation regularized computed tomography (CT), which solves 
\[
    \begin{array}{ll}
        \underset{x\in \reals^n}{\mbox{minimize}} & \displaystyle \frac{1}{2}\|Ex-b\|^2 + \lambda\|Dx\|_1,
    \end{array}
\]
where $x\in\reals^n$ is a vectorized image, $E\in\reals^{m\times n}$ is the discrete Radon transform, $b=Ex$ is the measurement, and $D$ is the finite difference operator.
We use primal-dual hybrid gradient (PDHG) \cite{zhu2008efficient,pock2009algorithm,esser2010general,chambolle2011first}, an instance of a nonexpansive fixed-point iteration via variable metric PPM \cite{he2012convergence}.
The results of \cref{fig:ct-fpr} indicate that restarted OC-Halpern \eqref{algo:os-ppm-res} provides an acceleration.

\subsection{Earth mover's distance} \label{exper:emd}
Consider the earth mover's distance between two probability measures, also referred to as the Wasserstein distance or the optimal transport problem.
The distance is defined through the discretized optimization problem
\[
    \begin{array}{ll}
        \underset{m_x, m_y}{\mbox{minimize}} & \displaystyle \|\mathbf{m}\|_{1,1} = \sum_{i=1}^n \sum_{j=1}^n |m_{x,ij}| + |m_{y,ij}| \\
        \mbox{subject to} & \mathrm{div}(\mathbf{m}) + \rho_1 - \rho_0 = 0,
    \end{array}
\]
where $\rho_0, \rho_1$ are probability measures on $\reals^{n\times n}$, $\mathrm{div}$ is a discrete divergence operator, and $\mathbf{m} = (m_x, m_y) \in \reals^{(n-1)\times n} \times \reals^{n\times(n-1)}$ is the optimization variable.
We use the algorithm of \citet{li2018parallel}, an instance of a nonexpansive fixed-point iteration via PDHG.
The results of \cref{fig:emd-fpr} indicate that restarted OC-Halpern \eqref{algo:os-ppm-res} provides an acceleration.

\subsection{Decentralized optimization with PG-EXTRA} \label{exper:decentralized}
Consider a decentralized optimization setting where each agent $i\in\{1,2,\dots,n\}$ has access to the sensing matrix $A_{(i)}\in\reals^{m_i\times n}$ and the noisy measurement $ b_{(i)}\approx A_{(i)}x$.
The goal is to recover the sparse signal $x\in\reals^n$ by solving the following compressed sensing problem:
\[
    \begin{array}{ll}
        \underset{x\in \reals^n}{\mbox{minimize}} & \displaystyle \frac{1}{n} \sum_{i=1}^n \|A_{(i)} x - b_{(i)}\|^2 + \lambda \|x\|_1.
    \end{array}
\]
We use PG-EXTRA \cite{shi2015proximal},
which is an instance of a nonexpansive fixed-point iteration via the Condat--V\~u \cite{condat2013primal,vu2013splitting} splitting method \cite{wu2017decentralized}.
The results of \cref{fig:decen-fpr} indicate that restarted OC-Halpern \eqref{algo:os-ppm-res} provides an acceleration.

\section{Conclusion}
This work presents an acceleration mechanism for fixed-point iterations and provides an exact matching complexity lower bound. The acceleration mechanism is an instance of Halpern's method, also referred to as anchoring, and the complexity lower bound is based on an explicit construction satisfying the zero-chain condition.

In this work, we measure the suboptimality of iterates with the fixed-point residual.
However, the fixed-point iteration is a \emph{meta}-algorithm, and almost all instances of it have further specific structure and suboptimality measures that are better suited for the particular problem of interest, such as function-value suboptimality, infeasibility for constrained problems, and primal-dual gap for minimax problems. Therefore, the fact that our proposed method accelerates the reduction of the fixed-point residual does not necessarily imply that it accelerates the reduction of the problem-specific suboptimality measure of practical interest.

Interestingly, the experimental results of Sections~\ref{sec:experiment} and \ref{sec:expr-details} indicate that our proposed acceleration does indeed provide a benefit in practice. This raises the following question: Under what setups can we expect anchoring-based acceleration to theoretically provide a benefit in terms of other suboptimality measures? Investigating this question would be an interesting direction of future work.

\newpage
\section*{Acknowledgements}
JP and EKR were supported by the National Research Foundation of Korea (NRF) Grant funded by the Korean Government (MSIP) [No. 2022R1C1C1010010], the National Research Foundation of Korea (NRF) Grant funded by the Korean Government (MSIP) [No. 2022R1A5A6000840], and the Samsung Science and Technology Foundation (Project Number SSTF-BA2101-02).
We thank TaeHo Yoon for providing careful reviews and valuable feedback. We thank Jelena Diakonikolas for the discussion on the prior work on complexity lower bounds of the fixed-point iterations. Finally, we thank the anonymous reviewers for their thoughtful comments.

\bibliography{main-camera-ready}
\bibliographystyle{icml2022}


\newpage
\onecolumn
\appendix

\section{Omitted proofs of Section \ref{sec:equivalence}} \label{proof:equivalence}
\begin{proof}[Proof of \cref{equiv:sm-lips} with inequalities]
    Suppose $\opT\colon\reals^n\to\reals^n$ is $\frac{1}{\gamma}$-Lipschitz for $\gamma\ge1$.
    Define $\opA\colon\reals^n\rightrightarrows\reals^n$ as
    \[
        \opA = \left( \opT + \frac{1}{\gamma} \opI \right)^{-1} \left(1 + \frac{1}{\gamma} \right) - \opI.
    \]
    For any $x, y \in \reals^n$, let $u \in \opA x$ and $v \in \opA y$.
    Then
    \begin{align*}
        u \in \opA x
        &\implies
        u \in \left( \opT + \frac{1}{\gamma} \opI \right)^{-1} \left(1 + \frac{1}{\gamma} \right) x - x \\
        &\iff
        x+u \in \left( \opT + \frac{1}{\gamma} \opI \right)^{-1} \left(1 + \frac{1}{\gamma} \right) x \\
        &\iff
        \left( \opT + \frac{1}{\gamma} \opI \right) (x+u) = x + \frac{1}{\gamma} x \\
        &\iff
        \opT (x+u) = x - \frac{1}{\gamma} u
    \end{align*}
    Likewise,
    \[
        \opT (y+v) = y - \frac{1}{\gamma} v.
    \]
    From the $\frac{1}{\gamma}$-Lipschitzness of $\opT$,
    \begin{align*}
        \|\opT (x+u) - \opT (y+v)\| \le \frac{1}{\gamma} \|(x+u)-(y+v)\|
        &\iff
        \left\|\left(x-\frac{1}{\gamma}u\right) - \left(y-\frac{1}{\gamma}v\right) \right\|
        \le
        \frac{1}{\gamma} \|(x+u)-(y+v)\| \\
        &\iff
        \left\|(x-y)-\frac{1}{\gamma}(u-v)\right\|^2
        \le
        \frac{1}{\gamma^2} \|(x-y)+(u-v)\|^2 \\
        &\iff
        \left(1 - \frac{1}{\gamma^2}\right) \|x-y\|^2 \le \left( \frac{2}{\gamma^2} + \frac{2}{\gamma} \right) \langle u-v, x-y \rangle \\
        &\iff
        \langle u-v, x-y \rangle \ge \frac{\gamma-1}{2} \|x-y\|^2.
    \end{align*}
    This holds for any $u\in\opA x$ and $v\in\opA y$ for any $x,y\in\reals^n$, so $\opA$ is $\frac{\gamma-1}{2}$-strongly monotone.
    
    We can further prove that
    \begin{align*}
        x_\star\in\zer\opA
        &\iff
        0 \in \opA x_\star = \left(\opT + \frac{1}{\gamma}\opI\right)^{-1} \left(x_\star+\frac{1}{\gamma}x_\star \right) - x_\star \\
        &\iff
        x_\star \in \left(\opT + \frac{1}{\gamma}\opI\right)^{-1} \left(x_\star+\frac{1}{\gamma}x_\star \right) \\
        &\iff
        \opT x_\star + \frac{1}{\gamma} x_\star = x_\star + \frac{1}{\gamma} x_\star \\
        &\iff
        x_\star = \opT x_\star \\
        &\iff
        x_\star \in \fix\opT.
    \end{align*}

    Suppose $\opA\colon\reals^n\rightrightarrows\reals^n$ is $\mu$-strongly monotone for $\mu\ge0$.
    Define $\opT\colon\reals^n\to\reals^n$ as
    \[
        \opT = \left(1 + \frac{1}{1+2\mu} \right) \resA - \frac{1}{1+2\mu} \opI.
    \]
    For any $x,y\in\reals^n$, let $u = \opT x$ and $v = \opT y$.
    Then
    \begin{align*}
        &u = \opT x \\
        &\iff
        u = \left(1 + \frac{1}{1+2\mu} \right) \resA x - \frac{1}{1+2\mu} x \\
        &\iff
        \frac{1}{1+2\mu} x + u = \left(1 + \frac{1}{1+2\mu}\right) \resA x \\
        &\iff
        \frac{1+2\mu}{2+2\mu}\left( \frac{1}{1+2\mu} x + u \right) = \frac{1}{2+2\mu}x + \frac{1+2\mu}{2+2\mu}u = \resA x \\
        &\iff
        x \in (\opI+\opA) \left( \frac{1}{2+2\mu}x + \frac{1+2\mu}{2+2\mu}u \right) \\
        &\iff
        \frac{1+2\mu}{2+2\mu} (x-u) \in \opA \left( \frac{1}{2+2\mu}x + \frac{1+2\mu}{2+2\mu}u \right).
    \end{align*}
    Likewise,
    \[
        \frac{1+2\mu}{2+2\mu} (y-v) \in \opA \left( \frac{1}{2+2\mu}y + \frac{1+2\mu}{2+2\mu}v \right).
    \]
    From the $\mu$-strong monotonicity of $\opA$,
    \begin{align*}
        &\left\langle \opA \left( \frac{1}{2+2\mu}x + \frac{1+2\mu}{2+2\mu}u \right) - \opA \left( \frac{1}{2+2\mu}y + \frac{1+2\mu}{2+2\mu}v \right), \left( \frac{1}{2+2\mu}x + \frac{1+2\mu}{2+2\mu}u \right) - \left( \frac{1}{2+2\mu}y + \frac{1+2\mu}{2+2\mu}v \right) \right\rangle \\
        &\ge
        \mu \left\| \left( \frac{1}{2+2\mu}x + \frac{1+2\mu}{2+2\mu}u \right) - \left( \frac{1}{2+2\mu}y + \frac{1+2\mu}{2+2\mu}v \right) \right\|^2 \\
        &\implies
        \left\langle \frac{1+2\mu}{2+2\mu} (x-u) - \frac{1+2\mu}{2+2\mu} (y-v), \left( \frac{1}{2+2\mu}x + \frac{1+2\mu}{2+2\mu}u \right) - \left( \frac{1}{2+2\mu}y + \frac{1+2\mu}{2+2\mu}v \right) \right\rangle \\
        &\qquad\ge
        \mu \left\| \left( \frac{1}{2+2\mu}x + \frac{1+2\mu}{2+2\mu}u \right) - \left( \frac{1}{2+2\mu}y + \frac{1+2\mu}{2+2\mu}v \right) \right\|^2 \\
        &\iff
        \left\langle \frac{1+2\mu}{2+2\mu} (x-y) - \frac{1+2\mu}{2+2\mu} (u-v), \frac{1}{2+2\mu}(x-y) + \frac{1+2\mu}{2+2\mu}(u-v) \right\rangle \\
        &\qquad\ge
        \mu \left\| \frac{1}{2+2\mu}(x-y) + \frac{1+2\mu}{2+2\mu}(u-v) \right\|^2 \\
        &\iff
        \langle (1+2\mu)(x-y) - (1+2\mu)(u-v), (x-y) + (1+2\mu)(u-v) \rangle
        \ge
        \mu \|(x-y) + (1+2\mu)(u-v)\|^2 \\
        &\iff
        (1+\mu) \|x-y\|^2 \ge (1+\mu)(1+2\mu)^2\|u-v\|^2 \\
        &\iff
        \|u-v\|^2 \le \frac{1}{(1+2\mu)^2} \|x-y\|^2.
    \end{align*}
    This holds for any $u=\opT x$ and $v = \opT y$ for any $x,y\in\reals^n$, so $\opT$ is $\frac{1}{1+2\mu}$-Lipschitz.
    
    Finally, we can also prove that
    \begin{align*}
        x_\star \in \fix\opT
        &\iff
        x_\star = \opT x_\star = \left(1+\frac{1}{1+2\mu}\right)\resA x_\star - \frac{1}{1+2\mu} x_\star \\
        &\iff
        \frac{2+2\mu}{1+2\mu} x_\star = \frac{2+2\mu}{1+2\mu} \resA x_\star \\
        &\iff
        x_\star = \resA x_\star = (\opI+\opA)^{-1} x_\star \\
        &\iff
        x_\star \in x_\star + \opA x_\star \\
        &\iff
        0\in\opA x_\star \\
        &\iff
        x_\star \in \zer\opA.
    \end{align*}
\end{proof}

\begin{proof}[Proof of \cref{equiv:sm-lips} with scaled relative graph]
    In this proof, we use the notations of \citet{ryu2021scaled} for the operator classes, which we list below.
    Consider a class of operators $\cM_\mu$ of $\mu$-strongly monotone operators and $\cL_{1/\gamma}$ of $\frac{1}{\gamma}$-contractions.
    As $\cM_\mu$, $\cL_{1/\gamma}$ are SRG-full classes, which means that the inclusion of the SRG of some operator to the SRG of an operator class is equivalent to membership of that operator to the given operator class \citep[Section~3.3]{ryu2021scaled}. Instead of showing that the operators satisfy the equivalent inequality condition to the membership, we show the membership in terms of the SRGs.
    
    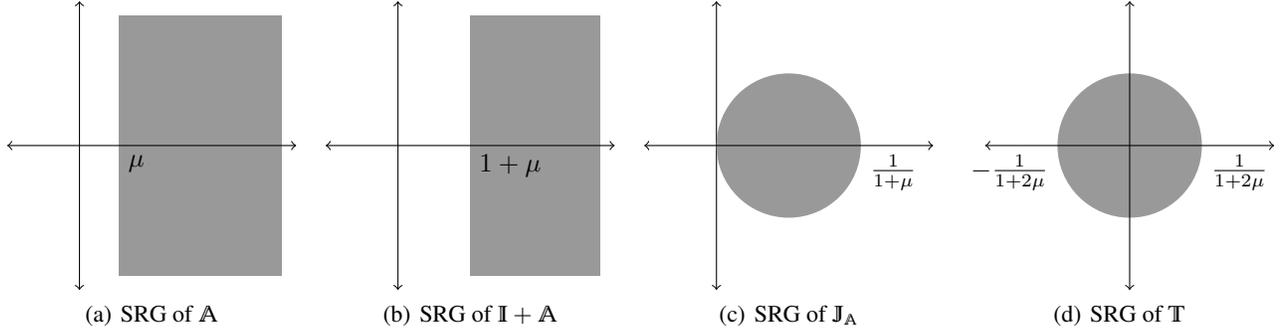
\begin{figure}[ht]
    \vskip 0.1in
        \centering
        \subfigure[SRG of $\opA$]{
            \begin{tikzpicture}[scale=0.48]
                \draw [<->] (-2,4) -- (-2,-4);
                \def\m{.9}
                \fill [fill=medgrey]
                (-\m,-4*\m) -- (-\m,4*\m) -- (4*\m,4*\m) -- (4*\m,-4*\m) -- cycle;
                \draw [<->] (-4,0) -- (4,0);
                \draw (-\m,0) node [below right] {$\mu$};
            \end{tikzpicture}
        } \hfill
        \subfigure[SRG of $\opI + \opA$]{
            \begin{tikzpicture}[scale=0.48]
                \draw [<->] (-2,4) -- (-2,-4);
                \def\m{.9}
                \fill [fill=medgrey]
                (0,-4*\m) -- (0,4*\m) -- (4*\m,4*\m) -- (4*\m,-4*\m) -- cycle;
                \draw [<->] (-4,0) -- (4,0);
                \draw (0,0) node [below right] {$1+\mu$};
            \end{tikzpicture}
        } \hfill
        \subfigure[SRG of $\resA$]{
            \begin{tikzpicture}[scale=0.48]
                \draw [<->] (-2,4) -- (-2,-4);
                \fill [fill=medgrey] (0,0) circle (2);
                \draw [<->] (-4,0) -- (4,0);
                \draw (2,0) node [below right] {$\frac{1}{1+\mu}$};
            \end{tikzpicture}
        } \hfill
        \subfigure[SRG of $\opT$]{
            \begin{tikzpicture}[scale=0.48]
                \fill [fill=medgrey] (0,0) circle (2);
                \draw [<->] (-4,0) -- (4,0);
                \draw [<->] (0,4) -- (0,-4);
                \draw (2,0) node [below right] {$\frac{1}{1+2\mu}$};
                \draw (-2,0) node [below left] {$-\frac{1}{1+2\mu}$};
            \end{tikzpicture}
        }
        \caption{SRG changing with invertible transformation $F$.}
    \vskip 0.1in
    \end{figure}
    
    Consider an invertible transformation $F : \mathbb{C} \cup \{\infty\} \to \mathbb{C} \cup \{\infty\}$ defined as
    \[
        F(z) = \left( 1 + \frac{1}{1+2\mu} \right) (1+z)^{-1} - \frac{1}{1+2\mu}.
    \]
    $F$ is a composition of only scalar addition/subtraction/multiplication and inversion, therefore preserves the SRG of $\cM_\mu$ and $\cL_{1/\gamma}$. 
    SRG of $F(\cM_\mu)$ and $\cL_{1/\gamma}$ match, and the SRG of $F^{-1}(\cL_{1/\gamma})$ and $\cM_\mu$ match.
\end{proof}

\section{Omitted proofs of Section \ref{sec:method}} \label{proof:method}

\subsection{Proof of \cref{equiv:method}}

\begin{lemma} \label{os-ppm-anchor}
    The $y_k$-update in algorithm \eqref{algo:os-ppm} is equivalent to
    \begin{align*}
        x_k &= \resA y_{k-1} \\
        y_k &= \left( 1 - \frac{1}{\varphi_k} \right) \left\{ \left(1 + \frac{1}{\gamma} \right) x_k - \frac{1}{\gamma} y_{k-1} \right\} + \frac{1}{\varphi_k} y_0
    \end{align*}
    where $\gamma=1+2\mu$.
\end{lemma}

\begin{proof}
    It suffices to show the equivalence of $y_k$-iterates.
    For $k=1$, from \eqref{algo:os-ppm} update,
    \begin{align*}
        y_1
        &=
        x_1 + \frac{\varphi_0-1}{\varphi_1} (x_1-y_0)
        - \frac{(\gamma-1)\varphi_0}{\varphi_1} (y_0-x_1) \\
        &=
        x_1 - \frac{\gamma-1}{\gamma^2+1} (y_0-x_1) 
        \tag{$\varphi_1 = 1 + \gamma^2$} \\
        &=
        \left(1 - \frac{1}{\varphi_1}\right) \left\{ \left(1+\frac{1}{\gamma}\right) x_1 - \frac{1}{\gamma} y_0 \right\} + \frac{1}{\varphi_1} y_0.
    \end{align*}
    
    Assume that the equivalence of the iterates holds for $k=1,2,\dots,l$.
    From the \eqref{algo:os-ppm} update,
    \begin{align*}
        y_{l+1}
        &=
        x_{l+1} + \frac{\varphi_l-1}{\varphi_{l+1}}(x_{l+1}-x_l)
        - \frac{(\gamma-1)\varphi_l}{\varphi_{l+1}}(y_l-x_{l+1}) + \frac{\gamma \varphi_{l-1}}{\varphi_{l+1}} (y_{l-1}-x_{l}) \\
        &=
        \left\{ 1 + \frac{\varphi_l-1}{\varphi_{l+1}} + \frac{(\gamma-1)\varphi_l}{\varphi_{l+1}} \right\} x_{l+1}
        - \left( \frac{\varphi_l-1}{\varphi_{l+1}} + \frac{\gamma \varphi_{l-1}}{\varphi_{l+1}} \right) x_l
        - \frac{(\gamma-1)\varphi_l}{\varphi_{l+1}} y_l + \frac{\gamma \varphi_{l-1}}{\varphi_{l+1}} y_{l-1} \\
        &= \gamma(\gamma+1)\frac{\varphi_l}{\varphi_{l+1}} x_{l+1} - \gamma(\gamma+1) \frac{\varphi_{l-1}}{\varphi_{l+1}} x_l
        - \frac{(\gamma-1)\varphi_l}{\varphi_{l+1}} y_l + \frac{\gamma \varphi_{l-1}}{\varphi_{l+1}} y_{l-1}.
    \end{align*}
    From the inductive hypothesis, we have 
    \[
        y_l = \left(1 - \frac{1}{\varphi_l}\right) \left\{ \left(1+\frac{1}{\gamma}\right) x_l - \frac{1}{\gamma} y_{l-1} \right\} + \frac{1}{\varphi_l} y_0,
    \]
    or
    \[
        \gamma \varphi_{l-1} y_{l-1} = \gamma (\gamma+1) \varphi_{l-1} x_l - \varphi_l y_l + y_0.
    \]
    Plugging this into the $\gamma \varphi_{l-1} y_{l-1}$-term in $y_{l+1}$, we get
    \begin{align*}
        y_{l+1}
        &=
        \gamma(\gamma+1) \frac{\varphi_l}{\varphi_{l+1}} x_{l+1} - \gamma(\gamma+1)\frac{\varphi_{l-1}}{\varphi_{l+1}} x_l
        - \frac{(\gamma-1)\varphi_l}{\varphi_{l+1}} y_l
        + \frac{1}{\varphi_{l+1}} \left\{ \gamma (\gamma+1) \varphi_{l-1} x_l - \varphi_l y_l + y_0 \right\} \\
        &=
        \gamma(\gamma+1)\frac{\varphi_l}{\varphi_{l+1}} x_{l+1} - \gamma \frac{\varphi_l}{\varphi_{l+1}}y_l + \frac{1}{\varphi_{l+1}} y_0 \\
        &=
        \frac{\gamma^2 \varphi_l}{\varphi_{l+1}} \left\{ \left(1 + \frac{1}{\gamma}\right) x_{l+1} - \frac{1}{\gamma} y_l \right\} + \frac{1}{\varphi_{l+1}} y_0 \\
        &=
        \left( 1 - \frac{1}{\varphi_{l+1}} \right) \left\{ \left(1 + \frac{1}{\gamma}\right) x_{l+1} - \frac{1}{\gamma} y_l \right\} + \frac{1}{\varphi_{l+1}} y_0.
    \end{align*}
    The same equivalence holds for $y_{l+1}$, so we are done.
\end{proof}

\begin{proof}[Proof of \cref{equiv:method}]
    Start from the same initial iterate $y_0 = \tilde{y}_0$.
    Suppose $y_k = \tilde{y}_k$ for some $k \ge 0$.
    Then,
    \begin{align*}
        y_{k+1}
        &=
        \left( 1 - \frac{1}{\varphi_{k+1}}\right) \left\{ \left(1+\frac{1}{1+2\mu} \right) x_{k+1} - \frac{1}{1+2\mu} y_k \right\}
        + \frac{1}{\varphi_{k+1}} y_0
        \tag{\cref{os-ppm-anchor}} \\
        &=
        \left( 1 - \frac{1}{\varphi_{k+1}} \right) \left\{ \left(1 + \frac{1}{1+2\mu} \right) \opJ_\opA - \frac{1}{1+2\mu} \opI \right\} y_{k}
        + \frac{1}{\varphi_{k+1}} y_0 \\
        &=
        \left( 1 - \frac{1}{\varphi_{k+1}} \right) \opT y_{k} + \frac{1}{\varphi_{k+1}} y_0 \\
        &=
        \left( 1 - \frac{1}{\varphi_{k+1}} \right) \opT \tilde{y}_{k} + \frac{1}{\varphi_{k+1}} y_0
        \tag{$y_k = \tilde{y}_k$ by induction hypothesis}
        =
        \tilde{y}_{k+1}.
    \end{align*}
\end{proof}

\subsection{Proof of \cref{lyapunov:os-ppm}}
Recall that
\begin{align*}
    V^k
    &= (1+\gamma^{-k}) \Bigg[ 
        \left(\sum_{n=0}^{k-1} \gamma^n\right)^2 \|\tilde{\opA}x_k\|^2
        + 2\left(\sum_{n=0}^{k-1}\gamma^n\right) \langle \tilde{\opA}x_k - \mu(x_k-x_\star), x_k-x_\star\rangle \\
        &\quad
        + \gamma^{-k} \left\| \left(\sum_{n=0}^{k-1}\gamma^n\right) \tilde{\opA}x_k - \gamma^k (x_k-x_\star) + (x_k-y_0) \right\|^2
    \Bigg]
    + (1-\gamma^{-k})\|y_0-x_\star\|^2.
    \tag{OS-PPM-Lyapunov}
    \label{lyapunov:potential-repeat}
\end{align*}

for $k=1,2,\dots,N$ and $V^0 = 2\|y_0-x_\star\|^2$, where $\gamma = 1+2\mu$, $\varphi_{k} = \sum_{n=0}^{k} \gamma^{2n}$ and $\tilde{\opA} x_k = y_{k-1} - x_k \in \opA x_k$.
We will often use the following identity.
\[
    (1+\gamma) \varphi_k 
    = (1+\gamma) \sum_{n=0}^{k} \gamma^{2n} 
    = (1+\gamma^{k+1}) \sum_{n=0}^{k} \gamma^n.
    \label{identity:step}
\]

First, we show that $V^k$ has an alternate form as below.
This form is useful in proving the monotone decreasing property of $V^k$ in $k$.
\begin{lemma} \label{lem:lyap-form}
    $V^k$ defined in \eqref{lyapunov:potential-repeat} can be equivalently written as
    
    \[
        V^k
        =
        \gamma^{-2k} (1+\gamma)^2 \varphi_{k-1}^2 \|\tilde{\opA} x_k\|^2
        + 2\gamma^{-2k} (1+\gamma) \varphi_{k-1} \langle \tilde{\opA} x_k - \mu(x_k-y_0), x_k-y_0 \rangle
        + 2\|y_0-x_\star\|^2.
    \]
\end{lemma}
\begin{proof}
    Expanding the square term,
    \begin{align*}
        &\left\| \left(\sum_{n=0}^{k-1}\gamma^n\right) \tilde{\opA}x_k - \gamma^k (x_k-x_\star) + (x_k-y_0) \right\|^2 \\
        &=
        \left\| \left(\sum_{n=0}^{k-1}\gamma^n\right) \tilde{\opA}x_k - (\gamma^k-1) (x_k-y_0) - \gamma^k(y_0-x_\star) \right\|^2 \\
        &=
        \left(\sum_{n=0}^{k-1}\gamma^n\right)^2 \|\tilde{\opA}x_k\|^2
        - 2\left(\sum_{n=0}^{k-1}\gamma^n\right)(\gamma^k-1) \langle \tilde{\opA}x_k, x_k-y_0 \rangle
        - 2\left(\sum_{n=0}^{k-1}\gamma^n\right)\gamma^k \langle \tilde{\opA}x_k, y_0-x_\star \rangle \\
        &\quad
        + (\gamma^k-1)^2 \|x_k-y_0\|^2 - 2\gamma^k(\gamma^k-1) \langle x_k-y_0, y_0-x_\star\rangle
        + \gamma^{2k}\|y_0-x_\star\|^2.
    \end{align*}
    Also, we have
    \begin{align*}
        &\langle \tilde{\opA}x_k - \mu(x_k-x_\star), x_k-x_\star \rangle \\
        &=
        \langle \tilde{\opA}x_k - \mu(x_k-y_0) - \mu(y_0-x_\star), (x_k-y_0)+(y_0-x_\star) \rangle \\
        &=
        \langle \tilde{\opA}x_k-\mu(x_k-y_0), x_k-y_0 \rangle + \langle \tilde{\opA}x_k, y_0-x_\star \rangle - 2\mu \langle x_k-y_0, y_0-x_\star \rangle.
    \end{align*}
    Then $V^k$ is expressed as
    \begin{align*}
        V^k
        &=
        2(1+\gamma^{-k}) \left( \sum_{n=0}^{k-1} \gamma^n \right) \langle \left\{
            \langle \tilde{\opA}x_k-\mu(x_k-y_0), x_k-y_0 \rangle + \langle \tilde{\opA}x_k, y_0-x_\star \rangle - 2\mu \langle x_k-y_0, y_0-x_\star \rangle
        \right\} \\
        &\qquad
        + (1+\gamma^{-k})\gamma^{-k} \Bigg\{
            \left(\sum_{n=0}^{k-1}\gamma^n\right)^2 \|\tilde{\opA}x_k\|^2
            - 2\left(\sum_{n=0}^{k-1}\gamma^n\right)(\gamma^k-1) \langle \tilde{\opA}x_k, x_k-y_0 \rangle
            + (\gamma^k-1)^2 \|x_k-y_0\|^2 \\
            &\qquad\qquad\qquad\qquad
            - 2\left(\sum_{n=0}^{k-1}\gamma^n\right)\gamma^k \langle \tilde{\opA}x_k, y_0-x_\star \rangle 
            - 2\gamma^k(\gamma^k-1) \langle x_k-y_0, y_0-x_\star\rangle
            + \gamma^{2k}\|y_0-x_\star\|^2
        \Bigg\} \\
        &\qquad
        + (1+\gamma^{-k}) \left( \sum_{n=0}^{k-1} \gamma^n \right)^2 \|\tilde{\opA}x_k\|^2 
        + (1-\gamma^{-k}) \|y_0-x_\star\|^2 \\
        &=
        (1+\gamma^{-k})^2 \left( \sum_{n=0}^{k-1} \gamma^n \right)^2 \|\tilde{\opA}x_k\|^2
        + 2(1+\gamma^{-k}) \left( \sum_{n=0}^{k-1} \gamma^n \right) \langle \tilde{\opA}x_k - \mu(x_k-y_0), x_k-y_0 \rangle \\
        &\qquad
        - 2\gamma^{-k}(1+\gamma^{-k})(\gamma^k-1) \left( \sum_{n=0}^{k-1} \gamma^n \right) \langle \tilde{\opA}x_k, x_k-y_0 \rangle
        + \gamma^{-k} (1+\gamma^{-k}) (\gamma^k-1)^2 \|x_k-y_0\|^2 
        + 2\|y_0-x_\star\|^2 \\
        &=
        (1+\gamma^{-k})^2 \left( \sum_{n=0}^{k-1} \gamma^n \right)^2 \|\tilde{\opA}x_k\|^2
        + 2\gamma^{-k}(1+\gamma^{-k}) \left( \sum_{n=0}^{k-1} \gamma^n \right) \langle \tilde{\opA}x_k-\mu(x_k-y_0), x_k-y_0 \rangle
        + 2\|y_0-x_\star\|^2.
    \end{align*}
    As
    \[
        (1+\gamma^{-k}) \left( \sum_{n=0}^{k-1} \gamma^n \right)
        =
        \frac{1+\gamma^k}{\gamma^k} \left( \sum_{n=0}^{k-1} \gamma^n \right)
        =
        \frac{1+\gamma}{\gamma^k}\varphi_{k-1},
    \]
    we have
    \begin{align*}
        V^k
        &=
        \gamma^{-2k}(1+\gamma)^2 \varphi_{k-1}^2 \|\tilde{\opA}x_k\|^2
        + 2\gamma^{-2k}(1+\gamma) \varphi_{k-1} \langle \tilde{\opA}x_k - \mu(x_k-y_0), x_k-y_0 \rangle
        + 2\|y_0-x_\star\|^2.
    \end{align*}
\end{proof}

Next, we prove that $\{V^k\}_{k=0}^N$ is monotonically decreasing in $k$.

\begin{lemma} \label{lem:lyap-mono}
    For $k=0,1,\dots,N$ with $V^k$ defined as \eqref{lyapunov:potential-repeat}, we have
    \[
        V^N \le V^{N-1} \le \dots \le V^1 \le V^0.
    \]
\end{lemma}
\begin{proof}
    We use the form of $V^k$ as in \cref{lem:lyap-form}.
    \begin{align*}
        V^{1} - V^0
        &=
        \gamma^{-2} (1+\gamma)^2 \|\tilde{\opA}x_1\|^2 + 2\gamma^{-2} (1+\gamma) \langle \tilde{\opA}x_1 - \mu(x_1-y_0), x_1-y_0 \rangle \\
        &=
        \gamma^{-2} (1+\gamma) \left\{ (1+\gamma) \|\tilde{\opA}x_1\|^2 + 2 \langle \tilde{\opA}x_1 - \mu(x_1-y_0), x_1-y_0 \rangle \right\} \\
        &=
        \gamma^{-2} (1+\gamma) \left\{ (1+\gamma) \|\tilde{\opA} x_1\|^2 - 2 (1+\mu) \|\tilde{\opA}x_1\|^2 \right\}
        \tag{$x_1-y_0 = - \tilde{\opA}x_1$} \\
        &=
        0.
        \tag{$1 + \gamma = 2(1+\mu)$}
    \end{align*}
    Now, consider $k \ge 1$. 
    Then,
    \begin{align*}
        V^{k+1} - V^k
        &=
        \gamma^{-2(k+1)} (1+\gamma)^2 \varphi_k^2 \|\tilde{\opA}x_{k+1}\|^2 - \gamma^{-2k} (1+\gamma)^2 \varphi_{k-1}^2 \|\tilde{\opA} x_k\|^2 \\
        &\quad
        + 2\gamma^{-2(k+1)} (1+\gamma) \varphi_k \langle \tilde{\opA}x_{k+1} - \mu(x_{k+1}-y_0), x_{k+1}-y_0 \rangle \\
        &\quad
        - 2\gamma^{-2k} (1+\gamma) \varphi_{k-1} \langle \tilde{\opA}x_k - \mu(x_k-y_0), x_k-y_0 \rangle .
    \end{align*}
    Now, we claim that
    \[
        V^{k+1} - V^k + 2\gamma^{-2k}(1+\gamma) \varphi_k \varphi_{k-1} \langle \tilde{\opA}x_{k+1} - \tilde{\opA}x_k - \mu(x_{k+1}-x_k), x_{k+1}-x_k \rangle = 0.
    \]
    First,
    \begin{align*}
        &V^{k+1} - V^k + 2\gamma^{-2k}(1+\gamma) \varphi_k \varphi_{k-1} \langle \tilde{\opA}x_{k+1} - \tilde{\opA}x_k - \mu(x_{k+1}-x_k), x_{k+1}-x_k \rangle \\
        &=
        V^{k+1} - V^k + 2\gamma^{-2k} (1+\gamma) \varphi_k \varphi_{k-1} \langle \tilde{\opA}x_{k+1} - \mu(x_{k+1}-y_0), x_{k+1}-x_k \rangle \\
        &\quad
        - 2\gamma^{-2k} (1+\gamma) \varphi_k \varphi_{k-1} \langle \tilde{\opA}x_k - \mu(x_k-y_0), x_{k+1}-x_k \rangle \\
        &=
        \gamma^{-2(k+1)} (1+\gamma)^2 \langle \varphi_k \tilde{\opA}x_{k+1} - \gamma \varphi_{k-1} \tilde{\opA}x_k, \varphi_k \tilde{\opA}x_{k+1} + \gamma \varphi_{k-1} \tilde{\opA}x_k \rangle \\
        &\quad
        + 2\gamma^{-2(k+1)} (1+\gamma) \varphi_k \langle \tilde{\opA}x_{k+1} - \mu(x_{k+1}-y_0), \gamma^2 \varphi_{k-1} (x_{k+1}-x_k) + (x_{k+1}-y_0) \rangle \\
        &\quad
        - 2\gamma^{-2k} (1+\gamma) \varphi_{k-1} \langle \tilde{\opA}x_k - \mu(x_k-y_0), \varphi_k (x_{k+1}-x_k) + (x_k-y_0) \rangle.
    \end{align*}
    From \cref{os-ppm-anchor}, we have
    \[
        y_k = \left( 1 - \frac{1}{\varphi_k} \right) \left\{ \left(1 + \frac{1}{\gamma} \right) x_k - \frac{1}{\gamma} y_{k-1} \right\} + \frac{1}{\varphi_k} y_0.
    \]
    Using the fact that $y_{k-1} = x_k + \tilde{\opA}x_k$, $y_k = x_{k+1} + \tilde{\opA} x_{k+1}$, and $\varphi_k = \gamma^2 \varphi_{k-1} + 1$, we obtain
    \[
        \varphi_k(x_{k+1}-y_0) + \varphi_k\tilde{\opA}x_{k+1}
        =
        \gamma^2 \varphi_{k-1} (x_k-y_0) - \gamma \varphi_{k-1} \tilde{\opA}x_k.
    \]
    Letting $U^k = \varphi_k (x_{k+1}-y_0) - \gamma^2 \varphi_{k-1} (x_k-y_0) = - \varphi_k \tilde{\opA}x_{k+1} - \gamma \varphi_{k-1} \tilde{\opA}x_k$, above formula is simplified as
    \begin{align*}
        &V^{k+1} - V^k + 2\gamma^{-2k}(1+\gamma) \varphi_k\varphi_{k-1} \langle \tilde{\opA}x_{k+1} - \tilde{\opA}x_k - \mu(x_{k+1}-x_k), x_{k+1}-x_k \rangle \\
        &=
        - \gamma^{-2(k+1)} (1+\gamma)^2 \langle \varphi_k \tilde{\opA}x_{k+1} - \gamma \varphi_{k-1} \tilde{\opA}x_k, U_k \rangle \\
        &\quad
        + 2\gamma^{-2(k+1)}(1+\gamma)\varphi_k \langle \tilde{\opA}x_{k+1} - \mu(x_{k+1}-y_0), U_k \rangle \\
        &\quad
        - 2\gamma^{-2k}(1+\gamma)\varphi_{k-1} \langle \tilde{\opA}x_k - \mu(x_k-y_0), U_k \rangle \\
        &=
        \gamma^{-2(k+1)}(1+\gamma) \big\langle
            - (1+\gamma)(\varphi_k\tilde{\opA}x_{k+1} - \gamma\varphi_{k-1}\tilde{\opA}x_k) + 2\varphi_k \{\tilde{\opA}x_{k+1} - \mu(x_{k+1}-y_0)\} \\
        &\qquad\qquad\qquad\qquad\qquad
            - 2\gamma^2\varphi_{k-1} \{\tilde{\opA}x_{k} - \mu(x_{k}-y_0)\}, U_k
        \big\rangle \\
        &=
        \gamma^{-2(k+1)}(1+\gamma) \big\langle (1-\gamma) (\varphi_k \tilde{\opA}x_{k+1} + \gamma \varphi_{k-1} \tilde{\opA}x_k) - 2\mu\{\varphi_k(x_{k+1}-y_0) - \gamma^2\varphi_{k-1} (x_k-y_0)\}, U_k \big\rangle \\
        &=
        \gamma^{-2(k+1)}(1+\gamma) \langle (\gamma-1) U_k - 2\mu U_k, U_k \rangle
        =
        0
        \tag{$\gamma-1=2\mu$}
    \end{align*}
\end{proof}

We now prove \cref{lyapunov:os-ppm}.
\begin{proof}[Proof of \cref{lyapunov:os-ppm}]
    According to \cref{lem:lyap-mono}, we have $V^N \le V^{N-1} \le \dots \le V^0 = 2\|y_0-x_\star\|^2$.
    Therefore,
    \begin{align*}
        2\|y_0-x_\star\|^2 &\ge V^N \\
        &=
        (1+\gamma^{-N}) \left(\sum_{n=0}^{N-1} \gamma^n\right)^2 \|\tilde{\opA}x_N\|^2 
        + 2(1+\gamma^{-N}) \left(\sum_{n=0}^{N-1} \gamma^n\right) \langle \tilde{\opA}x_N - \mu(x_N-x_\star), x_N-x_\star \rangle \\
        &\quad
        + \gamma^{-N}(1+\gamma^{-N}) \left\| \left(\sum_{n=0}^{N-1} \gamma^n\right) \tilde{\opA}x_N - \gamma^N (x_N-x_\star) + (x_N-y_0) \right\|^2 + (1-\gamma^{-N}) \|y_0-x_\star\|^2 \\
        &\ge
        (1+\gamma^{-N}) \left(\sum_{n=0}^{N-1} \gamma^n\right)^2 \|\tilde{\opA}x_N\|^2 + (1-\gamma^{-N}) \|y_0-x_\star\|^2,
    \end{align*}
    which can be simplified as
    \[
        (1+\gamma^{-N})\|y_0-x_\star\|^2 \ge (1+\gamma^{-N}) \left(\sum_{n=0}^{N-1} \gamma^n\right)^2 \|\tilde{\opA}x_N\|^2,
    \]
    or equivalently,
    \[
        \|\tilde{\opA}x_N\|^2 \le \left(\frac{1}{\sum_{n=0}^{N-1} \gamma^n}\right)^2 \|y_0-x_\star\|^2.
    \]
\end{proof}

\begin{proof}[Proof of \cref{lyapunov:oc-halpern}]
    This immediately follows from \cref{lyapunov:os-ppm} and \cref{equiv:method} by
    \[
        \tilde{\opA} x_N 
        = 
        y_{N-1} - x_N
        = 
        \left( 1 + \frac{1}{\gamma} \right)^{-1} (y_{N-1} - \opT y_{N-1}) 
        \in \opA x_N.
    \]
\end{proof}

\section{Details on the formulation of performance estimation problem for \eqref{algo:os-ppm}}
\label{detail:pep}

In order to obtain an estimate on the worst-case complexity of the algorithm, performance estimation problem (PEP) technique solves a certain form of semidefinite problem (SDP).
This SDP holds positive semidefinite matrix as an optimization variable, and solves the problem under constraints formulated from the interpolation condition of an operator in hand.

When discovering \eqref{algo:os-ppm}, we used maximal monotonicity as our interpolation condition, just as in \citet{ryu2020operator,kim2021accelerated}.
We further extended this to cover the case of maximal strongly-monotone operators, in a slightly different way with \citet{taylor2021optimal} who considered strongly convex interpolation.
The optimization variable is a positive semidefinite matrix, and this is of a Gram matrix form which stores information on the iterates of algorithms.
Usual choice of basis vectors for the gram matrix in PEP is usually $\nabla f(x)$ for convex minimization setup \cite{kim2016optimized,taylor2018exactpgm,taylor2021optimal}, or $\tilde{\opA} x$ for operator setup \cite{kim2021accelerated}.
Here, we used $x$-iterates to form the gram matrix of SDP.

This basic SDP is a primal problem of the PEP (Primal-PEP), and solving this returns an estimate to the worst-case complexity of given algorithm.
If we form a dual problem (dual-PEP) and minimize the optimal value of dual-PEP over possible choices of stepsizes as in \citet{kim2016optimized,taylor2018exactpgm,kim2021accelerated,taylor2021optimal}, this provides possibly the fastest rate, and solution to this minimization problem gives possibly optimal algorithms.
We considered a class of algorithms satisfying the span assumption in \cref{lower-bound:os-ppm}, and obtained \eqref{algo:os-ppm}.

\section{Omitted proofs of Section \ref{sec:lower-bound}} \label{proof:lower-bound}

\subsection{Proving complexity lower bound with span condition}

\begin{proof}[Proof of \cref{lem:equiv-ops} with inequalities]
    From \citep[Proposition~4.35]{bauschke2011convex}, $\opG$ is $\frac{1}{1+\gamma}$-averaged if and only if
    \[
        \|\opG x - \opG y\|^2 + \frac{\gamma-1}{\gamma+1}\|x-y\|^2 \le \frac{2\gamma}{1+\gamma} \langle \opG x - \opG y, x - y \rangle, \qquad \forall x, y \in \reals^n.
    \]
    Then for any $x,y\in\reals^n$, we get the chain of equivalences as follows.
    \begin{align*}
        \|\opT x - \opT y\|^2 \le \frac{1}{\gamma^2} \|x-y\|^2
        &\iff
        \| \gamma \opT x - \gamma \opT y \|^2 \le \|x-y\|^2 \\
        &\iff
        \| \{(1+\gamma)\opG x - \gamma x\} - \{(1+\gamma)\opG y - \gamma y\}\|^2 \le \|x-y\|^2 \\
        &\iff
        (1+\gamma)^2 \|\opG x - \opG y\|^2 - 2\gamma(1+\gamma) \langle \opG x - \opG y, x - y \rangle + \gamma^2 \|x-y\|^2 \le \|x-y\|^2 \\
        &\iff
        (1+\gamma)^2 \|\opG x - \opG y\|^2 + (\gamma^2-1)\|x-y\|^2 \le 2\gamma(1+\gamma)\langle \opG x - \opG y, x - y \rangle \\
        &\iff
        \|\opG x - \opG y \|^2 + \frac{\gamma-1}{\gamma+1}\|x-y\|^2 \le \frac{2\gamma}{\gamma+1}\langle \opG x - \opG y, x - y \rangle.
        \tag{$\because 1+\gamma>0$}
    \end{align*}
    Therefore, $\opT$ is $\frac{1}{\gamma}$-contractive if and only if $\opG$ is $\frac{1}{1+\gamma}$-averaged.
\end{proof}

\begin{proof}[Proof of \cref{lem:equiv-ops} with scaled relative graph]
    Using the notion of SRG \cite{ryu2021scaled}, we get the following equivalence of SRGs.
    Here, $\cN_{\frac{1}{1+\gamma}}$ is a class of $\frac{1}{1+\gamma}$-averaged operators.
    \begin{figure}[ht]
        \centering
        \mbox{
            \begin{tikzpicture}[scale=0.5]
                \fill [fill=medgrey] (0,0) circle (2);
                \draw [<->] (-3,0) -- (3,0);
                \draw [<->] (0,3) -- (0,-3);
                \draw (2,0) node [below right] {$\frac{1}{\gamma}$};
                \draw (-2,0) node [below left] {$-\frac{1}{\gamma}$};
                \draw (-3,5) node [below right] {SRG of $\opT \in \cL_\frac{1}{\gamma}$};
                \draw (-3,-5) node [above right] {SRG of $\opI - \left( 1+\frac{1}{\gamma} \right) \opG$};
            \end{tikzpicture}
            \hspace{1cm}
            \begin{tikzpicture}[scale=0.5]
                \fill [fill=medgrey] (0,0) circle (2);
                \draw [<->] (-3,0) -- (3,0);
                \draw [<->] (-1,3) -- (-1,-3);
                \draw (2,0) node [below right] {1};
                \draw (-2,0) node [below left] {$1-\frac{2}{1+\gamma}$};
                \draw (-3,5) node [below right] {SRG of $\frac{\gamma}{1+\gamma} (\opI-\opT)$};
                \draw (-3,-5) node [above right] {SRG of $\opG \in \cN_\frac{1}{1+\gamma}$};
            \end{tikzpicture}
        }
        \caption{SRG of $\opT$ and $\opG$}
    \end{figure}
    Therefore, we get the chain of equivalences
    \begin{align*}
        \opT \in \cL_{1/\gamma}
        &\iff
        \gamma \opT \in \cL_{1}
        \iff
        - \gamma \opT \in \cL_{1} \\
        &\iff
        \opG = \frac{\gamma}{1+\gamma} \opI + \frac{1}{1+\gamma} (- \gamma \opT) \in \cN_{\frac{1}{1+\gamma}},
    \end{align*}
    and conclude that $\opT$ is $\frac{1}{\gamma}$-Lipschitz if and only if $\opG$ is $\frac{1}{1+\gamma}$-averaged.
\end{proof}

\begin{proof}[Proof of \cref{lem:T-aver}]
    We restate the definition of $\opN \colon \reals^{N+1} \to \reals^{N+1}$.
    \[
        \opN x = \opN (x_1, x_2, \dots, x_{N+1}) = (x_{N+1}, -x_1, \dots, -x_{N}) - \frac{1+\gamma^{N+1}}{\sqrt{1+\gamma^2+\dots+\gamma^{2N}}} Re_1, \qquad x \in \reals^{N+1}.
    \]
    For any $x,y \in \reals^{N+1}$ such that
    \[
        x = (x_1, x_2, \dots, x_{N+1}), \qquad y = (y_1, y_2, \dots, y_{N+1}),
    \]
    we have
    \begin{align*}
        \| \opN x - \opN y \|^2
        &=
        \| (x_{N+1}, -x_1, \dots, -x_N) - (y_{N+1}, -y_1, \dots, -y_N) \|^2 \\
        &=
        (x_{N+1}-y_{N+1})^2 + (x_1-y_1)^2 + \dots + (x_N-y_N)^2 \\
        &=
        \| x - y \|^2.
    \end{align*}
    Then $\opN$ is nonexpansive, and by definition, $\opG = \frac{1}{1+\gamma} \opN + \frac{\gamma}{1+\gamma} \opI$ is a $\frac{1}{1+\gamma}$-averaged operator.
\end{proof}

\begin{proof}[Proof of \cref{lem:span}]
    By the definition of $\opG\colon\reals^{N+1}\to\reals^{N+1}$, for any $x\in\reals^{N+1}$,
    \begin{align*}
        \opG x
        &= 
        \frac{1}{1+\gamma} \opN x + \frac{\gamma}{1+\gamma} x \\
        &=
        \underbrace{
        \frac{1}{1+\gamma} \begin{bmatrix}
            \gamma & 0 & 0 & \hdots & 0 & 1 \\
            -1 & \gamma & 0 & \hdots & 0 & 0 \\
            0 & -1 & \gamma & \hdots & 0 & 0 \\
            \vdots & \vdots & \vdots & \ddots & \vdots & \vdots \\
            0 & 0 & 0 & \hdots & \gamma & 0 \\
            0 & 0 & 0 & \hdots & -1 & \gamma 
        \end{bmatrix}}_{= H} x
        - \underbrace{\frac{1}{1+\gamma} \frac{1+\gamma^{N+1}}{\sqrt{1+\gamma^2 + \dots + \gamma^{2N}}} Re_1}_{= b}
    \end{align*}
    where $\gamma=1+2\mu$.
    Observe that $\opG e_k \in \Span\{e_1, e_k, e_{k+1}\}$ for $k=1,\dots,N$.
    
    We use induction on $k$ to prove the Lemma.
    The claim holds for $k=0$ from
    \[
        \opG y_0 
        = \opG \mathbf{0}
        = - \frac{1}{1+\gamma} \frac{1+\gamma^{N+1}}{\sqrt{1+\gamma^2 + \dots + \gamma^{2N}}} Re_1
        \in \Span\{e_1\}.
    \]
    Now, suppose that the claim holds for $k < N$, i.e.,
    \begin{align*}
        y_k 
        &\in \Span\{ e_1, e_2, \dots, e_{k} \} \\
        \opG y_k 
        &\in \Span\{ e_1, e_2, \dots, e_{k+1} \}.
    \end{align*}
    Then
    \begin{align*}
        y_{k+1} 
        &\in 
        y_0 + \Span\{ \opG y_0, \opG y_1, \dots, \opG y_k \} \\
        &\subseteq
        \Span\{ e_1, e_2, \dots, e_{k+1} \} \\
        \opG y_{k+1}
        &=
        H y_{k+1} - b \\
        &\in
        H \Span\{ e_1, e_2, \dots, e_{k+1} \} - b \\
        &\subseteq
        \Span\{ e_1, e_2, \dots, e_{k+2} \}.
    \end{align*}
\end{proof}

\begin{proof}[Proof of \cref{lower-bound:oc-halpern}]
    The proof outline of \cref{lower-bound:oc-halpern} in \cref{subsec:lower-bound} is complete except for the part that the identity $(*)$ holds, and that \cref{lower-bound:oc-halpern} holds for any initial point $y_0\in\reals^n$ which is not necessarily zero.

    First, we show that for any initial point $y_0\in\reals^n$, there exists an worst-case operator $\opT\colon\reals^n\to\reals^n$ which cannot exhibit better than the desired rate.
    Denote by $\opT_0\colon\reals^n\to\reals^n$ the worst-case operator constructed in the proof of \cref{lower-bound:oc-halpern} for $y_0 = 0$.
    Define $\opT\colon\reals^n\to\reals^n$ as
    \[
        \opT y = \opT_0 (y-y_0) + y_0
    \]
    given $y_0\in\reals^n$.
    Then, first of all, the fixed point of $\opT$ is $y_\star = \tilde{y}_\star + y_0$ where $\tilde{y}_\star$ is the unique solution of $\opT_0$.
    Also, if $\{y_k\}_{k=0}^N$ satisfies the span condition
    \[
        y_k \in y_0 + \Span\left\{ y_0 - \opT y_0, \dots, y_{k-1} - \opT y_{k-1} \right\}, \qquad k=1,\dots,N,
    \]
    then $\tilde{y}_k = y_k - y_0$ forms a sequence satisfying
    \[
        \tilde{y}_k \in \underbrace{\tilde{y}_0}_{=0} + \Span \left\{ \tilde{y}_0 - \opT_0 \tilde{y}_0, \dots, \tilde{y}_{k-1} - \opT_0 \tilde{y}_{k-1} \right\},
        \qquad k=1,\dots,N,
    \]
    which is the same span condition in \cref{lower-bound:oc-halpern} with respect to $\opT_0$.
    This is true from the fact that
    \[
        y_k - \opT y_k 
        = \underbrace{y_k - y_0}_{= \tilde{y}_k} + \opT_0 (\underbrace{y_k - y_0}_{\tilde{y}_k})
        = \tilde{y}_k - \opT_0 \tilde{y}_k
    \]
    for $k=1,\dots,N$.
    
    Now, $\{\tilde{y}_k\}_{k=0}^N$ is a sequence starting from $\tilde{y}_0=0$ satisfying the span condition for $\opT_0$.
    This implies that, 
    \begin{align*}
        \|y_N - \opT y_N\|^2 
        &= 
        \|\tilde{y}_N - \opT_0 \tilde{y}_N\|^2 \\
        &\ge
        \left(1 + \frac{1}{\gamma}\right)^2 \left( \frac{1}{\sum_{k=0}^N \gamma^k} \right)^2 \|\tilde{y}_0 - \tilde{y}_\star\|^2 \\
        &=
        \left(1 + \frac{1}{\gamma}\right)^2 \left( \frac{1}{\sum_{k=0}^N \gamma^k} \right)^2 \|y_0 - y_\star\|^2.
    \end{align*}
    $\opT$ is our desired worst-case $\frac{1}{\gamma}$-contraction on $\reals^n$.

    It remains to show that
    \[
        \| \opG y_{N} \|^2
        \ge
        \left\| \cP_{\Span\{v\}} (b) \right\|^2\!
        =
        \left\|
            \frac{\langle b, v \rangle}{\langle v, v \rangle} v
        \right\|^2 
        \stackrel{(*)}{=}
        \left( \frac{1}{\sum_{k=0}^{N} \gamma^k} \right)^2 \!\!\!R^2
    \]
    where
    \[
        v = \begin{bmatrix}
            1 & \gamma & \gamma^2 & \dots & \gamma^N
        \end{bmatrix}^\intercal,
    \]
    especially the identity $(*)$.
    \begin{align*}
        \left\| \frac{\langle b, v \rangle}{\langle v, v \rangle}v \right\|^2
        &=
        \frac{|\langle b, v \rangle|^2}{\|v\|^2} \\
        &=
        \left( \frac{R}{1+\gamma} \times \frac{1+\gamma^{N+1}}{\sqrt{1+\gamma^2+\gamma^4+\dots+\gamma^{2N}}} \right)^2 \times \frac{1}{1+\gamma^2+\gamma^4+\dots+\gamma^{2N}} \\
        &=
        \left( \frac{R}{1+\gamma} \times \frac{1+\gamma^{N+1}}{1+\gamma^2+\gamma^4+\dots+\gamma^{2N}} \right)^2 \\
        &=
        \left( \frac{R}{1+\gamma+\gamma^2+\dots+\gamma^N} \right)^2 \\
        &=
        \left( \frac{1}{\sum_{k=0}^N \gamma^k} \right)^2 R^2.
    \end{align*}
\end{proof}

\begin{proof}[Proof of \cref{lower-bound:os-ppm}]
    According to \cref{equiv:sm-lips}, $\opT$ is $1/\gamma$-contractive if and only if $\opA = (\opT + 1/\gamma \opI)^{-1} - \opI$ is $\frac{\gamma-1}{2}$-strongly monotone.
    For any $y\in\reals^n$, if $x=\resA y$, then
    \[
        y - \opT y
        =
        y - \left\{ \left(1 + \frac{1}{\gamma}\right) \resA y - \frac{1}{\gamma} y \right\}
        =
        \left(1+\frac{1}{\gamma}\right) (y - x)
        =
        \left(1+\frac{1}{\gamma}\right) \tilde{\opA} x.
    \]
    This implies that
    \[
        y_k \in y_0 + \Span\{y_0-\opT y_0, y_1-\opT y_1, \dots, y_{k-1}-\opT y_{k-1}\},
        \qquad k=1,\dots,N,
    \]
    if and only if
    \begin{align*}
        x_k &= \resA y_{k-1} \\
        y_k &\in y_0 + \Span \{ \tilde{\opA}x_1, \dots, \tilde{\opA}x_k \}, \qquad k=1,\dots,N
    \end{align*}
    where $x_k = \resA y_{k-1}$.
    Span conditions in the statements of \cref{lower-bound:oc-halpern} and \cref{lower-bound:os-ppm} are equivalent under the transformation $\opA = (\opT+1/\gamma\opI)^{-1}-\opI$.
    Therefore, the lower bound result of this corollary can be derived from the lower bound result of \cref{lower-bound:oc-halpern}.
\end{proof}

\subsection{Deterministic algorithm classes}

In this section, we provide basic terminologies and necessary concepts in proving the complexity lower bound result for general algorithms.
We follow the \textit{information-based complexity} framework developed by \citet{nemirovski1983problem}, and use the resisting oracle technique to extend the results of \cref{lower-bound:oc-halpern} and \cref{lower-bound:os-ppm} to general fixed-point iterations and general proximal point methods.
The proof itself is motivated by the works of \citet{carmon2020stationary1,carmon2021stationary2}, and large portion of the definitions and notations are due to their work.

In the information-based complexity framework, every iterate $\{y_k\}_{k\in\mathbb{N}}$ is a query from an \textit{information oracle}, which returns restrictive information on a given function or operator.
Then, assumptions on the algorithm, such as linear span condition, illustrates how it uses such information.
For instance, provided with a gradient oracle $\cO_f(x) = \nabla f(x)$ of convex function $f$ to be minimized, usually the first-order algorithms search within the span of previous gradients to reach the next iterate.

A \textit{deterministic fixed-point iteration} $\al$ is a mapping of an initial point $y_0$ and an operator $\opT$ to a sequence of iterates $\{y_t\}_{t\in\mathbb{N}}$ and $\{\bar{y}_t\}_{t\in\mathbb{N}}$, such that the output depends on $\opT$ only through the \textit{fixed-point residual oracle} $\cO_\opT(y) = y - \opT y$.
Here, `deterministic' means that given the same initial point $y_0$
and the sequence of oracle evaluations $\{\cO_\opT(y_t)\}_{t\in \mathbb{N}}$, the algorithm yields the same sequence of iterates $\{(y_t, \bar{y}_t)\}_{t\in\mathbb{N}}$.
More precisely, we define $\al$ per iteration by setting $\al = \{\al_t\}_{t\in\mathbb{N}}$ 
with
\[
    (y_t, \bar{y}_t) = \al_t [y_0; \opT] = \al_t [y_0, \cO_\opT(y_0), \dots, \cO_\opT(y_{t-1})],
\]
where $y_t$ is the $t$\nobreakdash-th \textit{query point} and $\bar{y}_t$ is the $t$\nobreakdash-th \textit{approximate solution} produced by $\al_t$.
Here, we consider the algorithms whose query points and approximate solutions are identical ($y_t = \bar{y}_t)$.

Even though the $\al$ is defined to produce infinitely many $y_t$- and $\bar{y}_t$-iterates, the definition includes the case where algorithm terminates at a predetermined  total iteration count $N$, i.e., the algorithm may have a predetermined iteration count $N$ and the behavior may depend on the specified value of $N$.
In such cases, $y_N=\bar{y}_N= y_{N+1}=\bar{y}_{N+1} = \cdots$.

Similarly, a \emph{deterministic proximal point method} $\al$ is a mapping of an initial point $y_0$ and a maximal monotone operator $\opA$ to a sequence of query points $\{y_t\}_{t\in\mathbb{N}}$ and approximate solutions $\{\bar{y}_t\}_{t\in\mathbb{N}}$, such that the output depends on $\opA$ only through the resolvent residual oracle $\cO_\opA(y) = y - \resA y = \tilde{\opA} x \in \opA x$ where $x = \resA y$.
Indeed, this method $\al$ yields the same sequence of iterates given the same initial point $y_0$ and oracle evaluations $\{\cO_\opA(y_t)\}_{t\in \mathbb{N}}$.

\subsection{Generalized complexity lower bound}
As mentioned earlier, the general deterministic fixed-point iterations have no accounts for the span condition.
We use the resisting oracle technique \cite{nemirovski1983problem} to prove the lower bound result for general deterministic fixed-point iterations.
Recall that \cref{lower-bound:general-halpern} is
\renewcommand{\thesection}{\arabic{section}}
\setcounter{theorem}{5}
\begin{theorem}[Complexity lower bound of general deterministic fixed-point iterations]
    \label{lower-bound:general-halpern-repeat}
    Let $n\ge 2N$ for $N\in\mathbb{N}$.
    For any deterministic fixed-point iteration $\al$
    and any initial point $y_0\in\reals^n$, there exists a $\frac{1}{\gamma}$-Lipschitz operator $\opT \colon \reals^n \to \reals^n$ with a fixed point $y_\star\in\fix\opT$ such that
    \[
        \|y_N - \opT y_N\|^2 \ge \left( 1 + \frac{1}{\gamma} \right)^2
        \left( \frac{1}{\sum_{k=0}^N \gamma^k} \right)^2 \|y_0-y_\star\|^2
    \]
    where $\{y_t\}_{t\in\mathbb{N}} = \al[y_0; \opT]$.
\end{theorem}
\setcounter{theorem}{0}
\renewcommand{\thesection}{\Alph{section}}
By the equivalency of the optimization problems and algorithms stated in \cref{equiv:sm-lips} and \cref{equiv:method}, \cref{lower-bound:general-halpern} also generalizes \cref{lower-bound:os-ppm} to general proximal point methods.
\begin{corollary}[Complexity lower bound of general proximal point methods]
    \label{lower-bound:general-ppm}
    Let $n\ge 2N-2$ for $N\in\mathbb{N}$.
    For any deterministic proximal point method $\al$ and arbitrary initial point $y_0\in\reals^n$, there exists a $\mu$-strongly monotone operator $\opA\colon\reals^n\to\reals^n$ with a zero $x_\star\in\zer\opA$ such that
    \[
        \|\tilde{\opA}x_N\|^2 \ge \left( \frac{1}{1+\gamma+\dots+\gamma^{N-1}} \right)^2 \|y_0-x_\star\|^2
    \]
    where $\{y_t\}_{t\in\mathbb{N}} = \al[y_0;\opT]$.
\end{corollary}
\setcounter{theorem}{1}

\subsection{Proof of \cref{lower-bound:general-halpern-repeat}}
In order to prove \cref{lower-bound:general-halpern-repeat}, we first extend the result of \cref{lower-bound:oc-halpern} to the \emph{zero-respecting} sequences, which is a requirement slightly more general than the span assumption.
The worst-case operator of \cref{lower-bound:oc-halpern} covers the case of zero-respecting sequences, and this result will be successfully extended to general deterministic fixed-point iterations.

We say that a sequence $\{z_t\}_{t\in\mathbb{N}\cup\{0\}} \subseteq \reals^d$ is \emph{zero-respecting with respect to $\opT$} if
\[
    \mathrm{supp} \{ z_t \} \subseteq \cup_{s<t} \mathrm{supp} \{ z_s - \opT z_s \}
\]
for every $t\in\mathbb{N}\cup\{0\}$, where $\mathrm{supp}\{z\} := \{i\in[d] \mid \langle z, e_i \rangle \neq 0\}$.
An deterministic fixed-point iteration $\al$ is called \emph{zero-respecting} if $\al$ generates a sequence $\{z_t\}_{t\in\mathbb{N}\cup\{0\}}$ which is zero-respecting with respect to $\opT$ for any nonexpansive $\opT\colon\reals^d\to\reals^d$.
Note that by definition, $z_0 = 0$.
And for notational simplicity, define $\mathrm{supp} V = \bigcup_{z\in V} \mathrm{supp} \{z\}$.

This property serves as an important intermediate step to the generalization of \cref{lower-bound:oc-halpern}, where its similar form called `zero-chain' has numerously appeared on the relevant references in convex optimization \cite{Nesterov2004convex,drori2017exact,carmon2020stationary1,drori2022oracle}.
The worst-case operator found in the proof of \cref{lower-bound:oc-halpern} still performs the best among all the zero-respecting query points with respect to $\opT$, according to the following lemma.

\begin{lemma}
    \label{lem:worst-case}
    Let $\opT \colon \reals^{N+1} \to \reals^{N+1}$ be the worst-case operator defined in the proof of \cref{lower-bound:oc-halpern}.
    If the iterates $\{ z_t \}_{t=0}^N$ are zero-respecting with respect to $\opT$,
    \[
        \| z_N - \opT z_N \|^2 \ge \left( 1 + \frac{1}{\gamma} \right)^2 \left( \frac{1}{1+\gamma+\dots+\gamma^N} \right)^2 \|z_0-z_\star\|^2
    \]
    for $z_\star\in\fix\opT$.
\end{lemma}

\begin{proof}
    Let $\opG$ be defined as in the proof of \cref{lower-bound:oc-halpern}.
    Then we have
    \[
        z \in \Span\{e_1,e_2,\dots,e_k\} \implies \opG z \in \Span\{e_1,e_2,\dots,e_{k+1}\}.
    \]
    
    We claim that any zero-respecting sequence $\{z_k\}_{k=0,1,\dots,N}$ satisfies
    \begin{align*}
        z_k &\in \Span \big\{e_1,e_2,\dots,e_k \big\} \\
        \opG z_k = \frac{\gamma}{1+\gamma}(z_k - \opT z_k) &\in \Span \big\{e_1,e_2,\dots,e_{k+1} \big\}
    \end{align*}
    for $k=0,1,\dots,N$, so that the lower bound result of \cref{lower-bound:oc-halpern} is applicable.
    
    If $k=0$, then $y_0 = 0$ and from this, $\opG 0 \in \Span \{e_1\}$. So the case of $k=0$ holds. 
    Now, suppose that $0 < k \le N$ and the claim holds for all $n < k$.
    Then $\opG z_n \in \Span \{e_1,\dots,e_{n+1}\} \subseteq \Span \{e_1,\dots,e_{k}\}$ for $0\le k<n$.
    $\{z_k\}_{k=0}^N$ is zero-respecting with respect to $\opT$, so
    \begin{align*}
        \mathrm{supp} \{z_k\} 
        &\subseteq
        \bigcup_{n<k} \mathrm{supp} \{z_n - \opT z_n\} \\
        &=
        \mathrm{supp} \big\{\opG z_0, \opG z_1, \dots, \opG z_{k-1}\big\} \\
        &\subseteq \mathrm{supp} \{e_1,e_2,\dots,e_k\}.
    \end{align*}
    Therefore, $z_k \in \Span\{e_1,e_2,\dots,e_k\}$, and $\opG z_k \in \Span\{e_1,e_2,\dots,e_{k+1}\}$.
    The claim holds for $k=1,\dots,N$.
    
    According to the proof of \cref{lower-bound:oc-halpern},
    \[
        \|z_N - \opT z_N\|^2
        \ge
        \left( 1 + \frac{1}{\gamma} \right)^2 \left( \frac{1}{1+\gamma+\dots+\gamma^N} \right)^2 \|z_0-z_\star\|^2
    \]
    for any zero-respecting iterates $\{z_k\}_{k=0}^N$ with respect to $\opT$.
\end{proof}

We say that a matrix $U\in\reals^{m\times n}$ with $m \ge n$ is \emph{orthogonal}, if each columns $\{u_i\}_{i=1}^n\subseteq\reals^m$ of $U$ as in
\[
    U = 
    \begin{bmatrix}
        \vert & \dots & \vert \\
        u_1 & \dots & u_n \\
        \vert & \dots & \vert
    \end{bmatrix}
\]
are orthonormal to each other, or in other words, $U^\intercal U = I_n$.
It directly follows that $UU^\intercal$ is an orthogonal projection from $\reals^m$ to the range $\cR(U)$ of $U$.
\begin{lemma}
    \label{lem:invariant}
    For any orthogonal matrix $U \in \reals^{m \times n}$ with $m \ge n$ and any arbitrary vector $y_0 \in \reals^m$, if $\opT\colon\reals^n\to\reals^n$ is a $\frac{1}{\gamma}$-contractive operator with $\gamma \ge 1$, then $\opT_U \colon \reals^m \to \reals^m$ defined as
    \[
        \opT_U(y) := U \opT U^\intercal(y-y_0) + y_0, \qquad \forall y \in \reals^m
    \]
    is also a $\frac{1}{\gamma}$-contractive operator.
    Furthermore, $z_\star\in\fix\opT$ if and only if $y_\star=y_0 + Uz_\star\in\fix\opT_U$.
\end{lemma}

\begin{proof}
    For any $x, z \in \reals^m$,
    \begin{align*}
        \| \opT_U x - \opT_U z \|
        &=
        \| U \opT U^\intercal (x-y_0) - U \opT U^\intercal (z-y_0) \| \\
        &=
        \| \opT U^\intercal (x-y_0) - \opT U^\intercal (z-y_0) \| 
        \tag{$U$ is an orthogonal matrix} \\
        &\le
        \frac{1}{\gamma} \| U^\intercal (x-y_0) - U^\intercal (z-y_0) \|
        \tag{$\opT$ is $\frac{1}{\gamma}$-contractive} \\
        &=
        \frac{1}{\gamma} \| UU^\intercal (x - z) \| \\
        &\le
        \frac{1}{\gamma} \| x - z \|.
        \tag{$UU^\intercal$ is an orthogonal projection onto $\cR(U)$}
    \end{align*}
    Now, suppose $z_\star$ is a fixed point of $\opT$. Then
    \begin{align*}
        \opT_U (y_\star)
        &=
        U\opT U^\intercal U z_\star + y_0
        =
        U\opT z_\star + y_0 \\
        &=
        Uz_\star + y_0 = y_\star
    \end{align*}
    so $y_\star$ is a fixed point of $\opT_U$.
    On the other hand, if $y_\star$ is a fixed point of $\opT_U$, then $z_\star = U^\intercal (y_\star - y_0)$ satisfies
    \begin{align*}
        \opT (z_\star)
        &=
        \opT U^\intercal (y_\star - y_0) \\
        &=
        U^\intercal U \opT U^\intercal (y_\star - y_0)
        \tag{$U^\intercal U = I_n$} \\
        &=
        U^\intercal (\opT_U y_\star - y_0) \\
        &=
        U^\intercal (y_\star - y_0) = z_\star
        \tag{$y_\star\in\fix\opT_U$}
    \end{align*}
    so it is a fixed point of $\opT$.
\end{proof}

\begin{lemma}
    \label{lem:span2}
    Let $\al$ be a general deterministic fixed-point iteration, and $\opT \colon \reals^n \to \reals^n$ be a ${1/\gamma}$-contractive operator.
    For $m \ge n+N-1$ and any arbitrary point $y_0 \in \reals^m$, there exists an orthogonal matrix $U \in \reals^{m\times n}$ and the iterates $\{y_t\}_{t=1}^{N} = \al[y_0; \opT_U]$ with the following properties.
    \begin{itemize}
        \item[(i)]
        Let $z^{(t)} := U^\intercal (y_t-y_0)$ for $t = 0,1,\dots,N$.
        Then $\{z^{(t)}\}_{t=0}^{N}$ is zero-respecting with respect to $\opT$.
        
        \item[(ii)] $\{z^{(t)}\}_{t=0}^{N}$ satisfies
        \[
            \| z^{(t)} - \opT z^{(t)} \| \le \| y_t - \opT_U y_t \|, \qquad t=0,\dots,N.
        \]
    \end{itemize}
\end{lemma}

\begin{proof}
    We first show that $(i)$ implies $(ii)$.
    From $(i)$, we know that $z^{(t)} = U^\intercal (y_t-y_0)$ for $t = 0,1,\dots,N$.
    Therefore,
    \begin{align*}
        \| z^{(t)} - \opT z^{(t)} \|
        &=
        \| U^\intercal (y_t-y_0) - \opT U^\intercal (y_t-y_0) \| \\
        &=
        \| UU^\intercal \{(y_t-y_0) - U \opT U^\intercal (y_t-y_0) \} \|
        \tag{$U$ is orthogonal} \\
        &=
        \| UU^\intercal (y_t-y_0) - UU^\intercal U \opT U^\intercal (y_t-y_0) \|
        \tag{$U^\intercal U = I_n$} \\
        &=
        \| UU^\intercal \{(y_t-y_0) - U \opT U^\intercal (y_t-y_0) \} \| \\
        &\le
        \| (y_t-y_0) - U \opT U^\intercal (y_t-y_0) \|
        \tag{$UU^\intercal$ is an orthogonal projection} \\
        &=
        \| y_t - \opT_U y_t \|.
        \tag{Definition of $\opT_U$}
    \end{align*}
    
    Now we prove the existence of orthogonal $U\in\reals^{m\times n}$ with $\{y_t\}_{t=0}^N = \al[y_0; \opT_U]$ and $(i)$ holds.
    In order to show the existence of such orthogonal matrix $U$ as in $(i)$, we provide the inductive scheme that finds the columns of $U$ at each iteration.
    Before describing the actual scheme, we first provide some observations useful to deriving the necessary conditions for the columns $\{u_i\}_{i=1}^n$ of $U$ to satisfy.
    
    Let $t \in \{1,\dots,N\}$, and define the set of indices $S_t$ as
    \[
        S_t = \cup_{s<t} \mathrm{supp} \{ z^{(s)} - \opT z^{(s)} \}.
    \]
    For $\{z^{(t)}\}_{t=0}^N$ to satisfy the zero-respecting property with respect to $\opT$, $z^{(t)}$ is required to satisfy
    \[
        \mathrm{supp} \{z^{(t)}\} \subseteq S_t
    \]
    for $t=1,\dots,N$.
    This requirement is fulfilled when
    \[
        y_t-y_0 \in \Span \{ u_i \}_{i\in S_t}
    \]
    or equivalently,
    \[
        \langle u_i, y_t-y_0 \rangle = 0
    \]
    for every $i \notin S_t$.
    Note that $z^{(0)} = U^\intercal (y_0-y_0) = 0$ is trivial.
    
    We now construct $U \in \reals^{m\times n}$.
    Note that $S_0 = \emptyset \subseteq S_1 \subseteq \dots \subseteq S_t$.
    $\{u_i\}_{i\in S_t\setminus S_{t-1}}$ is chosen inductively starting from $t=1$.
    Suppose we have already chosen $\{u_i\}_{i\in S_{t-1}}$.
    Choose $\{u_i\}_{i\in S_t \setminus S_{t-1}}$ from the orthogonal complement of 
    \[
        W_t := \Span \left(
            \{
                y_1-y_0, \cdots, y_{t-1}-y_0
            \}
            \cup
            \left\{
                u_i
            \right\}_{i\in S_{t-1}}
        \right)
    \]
    and let them be orthogonal to each other. 
    In case of $S_N \neq \emptyset$, for $i \notin S_N$, choose proper vectors $u_i$ so that $U$ becomes an orthogonal matrix.
    This is possible when the dimension of $W_t^\perp$ is large enough to draw $|S_t \setminus S_{t-1}|$-many orthogonal vectors, or in other words,
    \[
        \dim W_t^\perp \ge |S_t \setminus S_{t-1}|.
    \]
    From the assumption, $m-t+1 \ge m-N+1 \ge n$, so we have a guarantee that
    \[
        \dim W_t^\perp
        =
        m - \dim W_t 
        \ge
        m - \{(t-1) + |S_{t-1}|\}
        \ge
        |S_{t-1}^c| = n - |S_{t-1}|
        \ge
        |S_t \setminus S_{t-1}|.
    \]
    The columns $\{u_i\}_{i=1}^n$ of constructed $U$ satisfies $\langle u_i, y_t-y_0 \rangle = 0$ if $i \notin S_t$, for $t=1,\dots,N$.
    Therefore, 
    \[
        z^{(t)} = U^\intercal (y_t-y_0) 
        \in \Span\{e_i\}_{i\in S_t} 
    \]
    which leads to $\mathrm{supp} \{z^{(t)}\} \subseteq S_t$.
\end{proof}

We now prove the complexity lower bound result for general fixed-point iterations.
\begin{proof}[Proof of \cref{lower-bound:general-halpern-repeat}]
    For any deterministic fixed-point iteration $\al$ and initial point $y_0\in\reals^n$, consider a worst-case operator $\opT\colon\reals^{N+1}\to\reals^{N+1}$ defined in the proof of \cref{lower-bound:oc-halpern}.
    According to \cref{lem:span2}, there exists an orthogonal $U \in \reals^{n\times(N+1)}$ with $n\ge(N+1)+(N-1)=2N$ such that $z^{(k)} = U^\intercal (y_k - y_0)$ for $k=0,\dots,N$,
    \[
        \|z^{(k)} - \opT z^{(k)}\| \le \|y_k-\opT_U y_k\|,
        \quad k=0,\dots,N
    \]
    where the query points $\{y_k\}_{k=0}^N$ are generated from applying $\al$ to $\opT_U$ given initial point $y_0$, and $\{z^{(k)}\}_{k=0}^N$ is a zero-respecting sequence with respect to $\opT$.
    According to \cref{lem:worst-case},
    \[
        \|z^{(N)}-\opT z^{(N)}\|^2 \ge \left(1 + \frac{1}{\gamma}\right)^2 \left(\frac{1}{1+\gamma+\dots+\gamma^N}\right)^2 \|z^{(0)}-z_\star\|^2.
    \]
    According to \cref{lem:invariant}, $y_\star = y_0 + Uz_\star \in \fix\opT_U$ for $z_\star \in \fix\opT$, so
    \[
        \|y_0-y_\star\|^2 = \|U(z^{(0)}-z_\star)\|^2 = \|z^{(0)}-z_\star\|^2
    \]
    where the second identity comes from orthogonality of $U$.
    We may conclude that
    \[
        \|y_N-\opT_U y_N\|^2 \ge \left(1 + \frac{1}{\gamma}\right)^2 \left(\frac{1}{1+\gamma+\dots+\gamma^N}\right)^2 \|y_0-y_\star\|^2
    \]
    and that $\opT_U\colon\reals^n\to\reals^n$ is the desired worst-case $\frac{1}{\gamma}$-contraction with $n\ge 2N$.
\end{proof}

\section{Omitted proofs of Section \ref{sec:restart}} \label{proof:restart}

\subsection{Convergence rate of proximal point method}

\begin{lemma} \label{sharp:lem}
    Let $\{x_k\}_{k\in\mathbb{N}}$ be the iterates generated by applying PPM $x_{k+1} = \resA x_k$ starting from $x_0 \in \reals^n$, given a uniformly monotone operator $\opA$ with parameters $\mu>0$ and $\alpha>1$.
    Now let $A_k := \|x_k-x_\star\|^2$ and $B_k := \|\tilde{\opA} x_{k+1}\|^2$. Then for any $k \in \mathbb{N}\cup\{0\}$,
    \begin{align*}
        A_k &\ge A_{k+1} \left( 1 + \mu A_{k+1}^\frac{\alpha-1}{2} \right)^2 \\
        B_k &\ge B_{k+1}.
    \end{align*}
\end{lemma}
\begin{proof}
    Note that PPM update $x_{k+1} = \resA x_k$ is equivalent to $x_k = x_{k+1} + \tilde{\opA} x_{k+1}$ where $\tilde{\opA} x_{k+1} \in \opA x_{k+1}$.
    \[
        x_k - x_\star = (x_{k+1} + \tilde{\opA} x_{k+1}) - x_\star = (x_{k+1}-x_\star) + \tilde{\opA} x_{k+1}.
    \]
    Then
    \begin{align*}
        A_k
        &=
        A_{k+1} + B_{k} + 2 \langle \tilde{\opA} x_{k+1}, x_{k+1}-x_\star \rangle \\
        &\ge
        A_{k+1} + B_{k} + 2\mu \|x_{k+1}-x_\star\|^{\alpha+1} \\
        &=
        A_{k+1} + B_{k} + 2\mu A_{k+1}^\frac{\alpha+1}{2} \\
        &\ge
        A_{k+1} + \mu^2 A_{k+1}^{\alpha} + 2\mu A_{k+1}^\frac{\alpha+1}{2} \\
        &\ge
        A_{k+1} \left( 1 + \mu A_{k+1}^\frac{\alpha-1}{2} \right)^2
    \end{align*}
    where the second inequality follows from 
    \[
        \|\tilde{\opA} x_{k+1}\| \|x_{k+1}-x_\star\| \ge \langle \tilde{\opA} x_{k+1}, x_{k+1}-x_\star \rangle \ge \mu \|x_{k+1}-x_\star\|^{\alpha+1}.
    \]
    Also, from
    \begin{align*}
        B_{k} - B_{k+1}
        &=
        \|\tilde{\opA} x_k\|^2 - \|\tilde{\opA} x_{k+1}\|^2 \\
        &=
        (\|\tilde{\opA} x_k\|^2 + \|\tilde{\opA} x_{k+1}\|^2) - 2\|\tilde{\opA} x_{k+1}\|^2 \\
        &\ge
        2\langle \tilde{\opA} x_k, \tilde{\opA} x_{k+1} \rangle - 2\|\tilde{\opA} x_{k+1}\|^2
        \tag{Young's inequality} \\
        &=
        - 2\langle \tilde{\opA} x_{k+1} - \tilde{\opA} x_k, \tilde{\opA} x_{k+1} \rangle \\
        &=
        2\langle \tilde{\opA} x_{k+1} - \tilde{\opA} x_k, x_{k+1} - x_k \rangle
        \ge
        0,
        \tag{Monotonicity of $\opA$}
    \end{align*}
    we get $B_k \ge B_{k+1}$.
\end{proof}

\begin{theorem} \label{sharp:ppm-x}
    If $\opA \colon \reals^n \rightrightarrows \reals^n$ is a uniformly monotone operator with parameters $\mu > 0$ and $\alpha > 1$, there exists $C > 0$ such that the iterates $\{x_k\}_{k\in\mathbb{N}}$ generated by PPM exhibits the rate
    \[
        \|x_k-x_\star\|^2 \le \frac{C}{k^\frac{2}{\alpha-1}}
    \]
    for any $k \in \mathbb{N}$.
\end{theorem}

\begin{proof}
    We use the induction on $k$ to show the convergence rate, and find the necessary conditions for $C>0$ to satisfy.
    
    In case of $k=1$, $\|x_1 - x_\star\|^2 \le C$ must be satisfied.
    \cref{sharp:lem} implies the monotonicity of $A_k$, so $C$ with $C \ge \|x_0-x_\star\|^2$ is a suitable choice.
    
    Now, suppose that $A_k \le C k^{-\frac{2}{\alpha-1}}$ and $k \ge 1$.
    We claim that $A_{k+1} \le  C (k+1)^{-\frac{2}{\alpha-1}}$ for the same $C>0$.
    Define $f^\alpha_\mu \colon [0,\infty) \to [0,\infty)$ as
    \[
        f^\alpha_\mu(t) := t \left(1+\mu t^\frac{\alpha-1}{2}\right)^2.
    \]
    Then $f^\alpha_\mu (A_{k+1}) \le A_k$ from \cref{sharp:lem}.
    If $f^\alpha_\mu(A_{k+1}) \le f^\alpha_\mu \left( C (k+1)^{-\frac{2}{\alpha-1}} \right)$, since $f^\alpha_\mu$ is a monotonically increasing function over $[0,\infty)$, we are done.
    Define
    \[
        a_n := (n+1) \left\{ \left(1+\frac{1}{n}\right)^\frac{1}{\alpha-1} - 1 \right\},
    \]
    and function $g\colon (0,\infty) \to \reals$ as
    \[
        g(x) = \left( 1 + \frac{1}{x} \right) \left\{ (1+x)^\frac{1}{\alpha-1} - 1 \right\}
    \]
    so that $a_n = g\left(\frac{1}{n}\right)$ for $n\in\mathbb{N}$.
    Then
    \begin{align*}
        g'(x)
        &=
        -\frac{1}{x^2} \left\{ (1+x)^\frac{1}{\alpha-1} - 1 \right\} + \frac{1}{\alpha-1} \left(1 + \frac{1}{x}\right) (1+x)^{\frac{1}{\alpha-1}-1} \\
        &=
        - \frac{(1+x)^\frac{1}{\alpha-1}}{x^2} + \frac{1}{x^2} + \frac{1}{\alpha-1} \frac{(1+x) (1+x)^{\frac{1}{\alpha-1}-1}}{x} \\
        &=
        \frac{- (1+x)^\frac{1}{\alpha-1} + 1 + {\scriptstyle\frac{x}{\alpha-1}} (1+x)^\frac{1}{\alpha-1}}{x^2} \\
        &=
        \frac{(1+x)^\frac{1}{\alpha-1}}{x^2} \left\{ (1+x)^{-\frac{1}{\alpha-1}} - \left(1 - {\scriptstyle\frac{1}{\alpha-1}}x\right) \right\}.
    \end{align*}
    As $x\mapsto (1+x)^{-\frac{1}{\alpha-1}}$ is a convex function on $[0, \infty)$ and $x\mapsto 1 - \frac{1}{\alpha-1}x$ is a first-order approximation at $0$ of it, $g'(x) \ge 0$ for $x>0$.
    $g$ is a monotonically increasing function, so $g$ obtains its maximum in $(0,1]$ at $x=1$, and we have
    \[
        \sup_{n\in\mathbb{N}} a_n = \sup_{n\in\mathbb{N}} g\left(\frac{1}{n}\right) = g(1) = 2(2^\frac{1}{\alpha-1} - 1) = 2^\frac{\alpha}{\alpha-1} - 2.
    \]
    The boundedness of $a_n$ leads to the equivalency as
    \begin{align*}
        a_k \le 2^\frac{\alpha}{\alpha-1} - 2
        &\iff
        \left(1 + \frac{1}{k}\right)^\frac{2}{\alpha-1} \le \left( 1 + \frac{2^\frac{\alpha}{\alpha-1} - 2}{k+1} \right)^2 \\
        &\iff
        \frac{C}{k^\frac{2}{\alpha-1}} \le \frac{C}{(k+1)^\frac{2}{\alpha-1}} \left\{ 1 + \frac{2^\frac{\alpha}{\alpha-1} - 2}{C^\frac{\alpha-1}{2}} \left(\frac{C}{(k+1)^\frac{2}{\alpha-1}}\right)^\frac{\alpha-1}{2} \right\}^2,
    \end{align*}
    for any choice of $C > 0$.
    Choosing $C \ge \mu^{-\frac{2}{\alpha-1}} (2^\frac{\alpha}{\alpha-1} - 2)^\frac{2}{\alpha-1}$ which is equivalent to
    \[
        \frac{2^\frac{\alpha}{\alpha-1} - 2}{C^{\frac{\alpha-1}{2}}} \le \mu,
    \]
    we get
    \begin{align*}
        \frac{C}{(k+1)^\frac{2}{\alpha-1}} \left\{ 1 + \frac{2^\frac{\alpha}{\alpha-1} - 2}{C^\frac{\alpha-1}{2}} \left(\frac{C}{(k+1)^\frac{2}{\alpha-1}}\right)^\frac{\alpha-1}{2} \right\}^2
        &\le 
        \frac{C}{(k+1)^\frac{2}{\alpha-1}} \left\{ 1 + \mu \left(\frac{C}{(k+1)^\frac{2}{\alpha-1}}\right)^\frac{\alpha-1}{2} \right\}^2 \\
        &=
        f^\alpha_\mu \left( \frac{C}{(k+1)^\frac{2}{\alpha-1}} \right).
    \end{align*}
    Gathering all the inequalities above, if $C \ge \mu^{-\frac{2}{\alpha-1}} (2^\frac{\alpha}{\alpha-1} - 2)^\frac{2}{\alpha-1}$, then
    \[
        f^\alpha_\mu(A_{k+1})
        \le
        A_k 
        \le
        \frac{C}{k^{\frac{2}{\alpha-1}}}
        \le
        f^\alpha_\mu \left( \frac{C}{(k+1)^\frac{2}{\alpha-1}} \right)
    \]
    so we get
    \[
        A_k \le \frac{C}{k^{\frac{2}{\alpha-1}}}
        \implies
        A_{k+1} \le \frac{C}{(k+1)^\frac{2}{\alpha-1}}
    \]
    for $k=1,2,\dots$.
    
    Therefore,
    \[
        \|x_k-x_\star\|^2 \le \frac{C}{k^{\frac{2}{\alpha-1}}}
        =
        \frac{\max\left\{ \left(\frac{2^\frac{\alpha}{\alpha-1}-2}{\mu}\right)^\frac{2}{\alpha-1}, \|x_0-x_\star\|^2 \right\}}{k^\frac{2}{\alpha-1}}.
    \]
\end{proof}

We now prove the convergence rate of PPM in terms of $B_k = \|\tilde{\opA} x_{k+1}\|^2$.

\begin{proof}[Proof of \cref{sharp:ppm}]
    We claim the convergence rate of $B_{k-1} = \|\tilde{\opA} x_k\|^2$ to be as above.
    From the proof of \cref{sharp:lem}, we have
    \[
        B_k \le A_k - A_{k+1} - 2\mu A_{k+1}^\frac{\alpha+1}{2} \le A_k - A_{k+1}.
    \]
    If $N = 1$, then
    \[
        B_0 \le A_0 - A_1 \le A_0 = \|x_0-x_\star\|^2.
    \]
    Suppose $N \ge 2$.
    Let $n := \lfloor \frac{N}{2} \rfloor$ where $\lfloor x \rfloor$ is the largest integer not exceeding $x$.
    Summing up the above inequality from $k=n$ to $k=N-1$ and using the monotonicity of $B_k$, we have
    \[
        \frac{N}{2} B_{N-1} \le \sum_{k=n}^{N-1} B_{N-1} \le \sum_{k=n}^{N-1} B_k \le \sum_{k=n}^{N-1} (A_k - A_{k+1}) = A_{n} - A_N \le A_{n}.
    \]
    Note that from the convergence analysis of $A_k$, or \cref{sharp:ppm-x}, we have
    \[
        A_{n} \le \frac{C}{n^\frac{2}{\alpha-1}}
    \]
    where $C = \max\left\{ \left(\frac{2^\frac{\alpha}{\alpha-1}-2}{\mu}\right)^\frac{2}{\alpha-1}, \|x_0-x_\star\|^2 \right\}$.
    Therefore,
    \[
        \frac{N}{2} B_{N-1} 
        \le \frac{C}{n^\frac{2}{\alpha-1}}
        \le \frac{C}{\left( \frac{N-1}{2} \right)^\frac{2}{\alpha-1}},
    \]
    so we may conclude that, for any $N \ge 2$,
    \begin{align*}
        B_{N-1} 
        &\le 
        \frac{2^\frac{\alpha+1}{\alpha-1} C}{(N-1)^\frac{2}{\alpha-1} N} \\
        &=
        \frac{ {2^\frac{\alpha+1}{\alpha-1} \max\left\{ \left(\frac{2^\frac{\alpha}{\alpha-1}-2}{\mu}\right)^\frac{2}{\alpha-1}, \|x_0-x_\star\|^2 \right\}} }{ {(N-1)^\frac{2}{\alpha-1}N} } \\
        &\le
        \frac{ 2^\frac{\alpha+3}{\alpha-1} \max\left\{ \left(\frac{2^\frac{\alpha}{\alpha-1}-2}{\mu}\right)^\frac{2}{\alpha-1}, \|x_0-x_\star\|^2 \right\} }{ N^\frac{\alpha+1}{\alpha-1} } \\
        &= 
        \mathcal{O} \left( N^{-\frac{\alpha+1}{\alpha-1}} \right)
    \end{align*}
    where the second inequality follows from $2(N-1) \ge N$.
    Since this bound also holds for the case of $N = 1$ from $B_0 \le \|x_0-x_\star\|^2$, we are done.
\end{proof}

\subsection{Convergence rate of restarted OS-PPM \eqref{algo:os-ppm-res}}

\citet{roulet2020sharpness} showed that 
if the objective function $f$ of a smooth convex minimization problem satisfies a \emph{H\"{o}lderian error bound condition}
\[
    \frac{\mu}{r} \|x - x_\star\|^r \le f(x) - f^\star, \qquad \forall\,x\in K \subset \reals^n
\]
where $x_\star\in K$ is a minimizer of $f$ and $K$ is a given set, then the unaccelerated base algorithm can be accelerated with a restarting scheme.
The restarting schedule uses $t_k$ iterations for each $k$\nobreakdash-th outer loop recursively satisfying
\[
    f(x_k) - f^\star \le e^{-\eta k} (f(x_0) - f^\star), \qquad k=1,2,\dots
\]
for some $\eta>0$, where $x_k=\mathcal{A}(x_{k-1},t_k)$ is the output of $k$\nobreakdash-th outer loop, which applies $t_k$ iterations of the base algorithm $\mathcal{A}$ starting from $x_{k-1}$.
If an objective function is strongly convex near the solution ($r=2$), a constant restarting schedule $t_k = \lambda$ provides a faster rate compared to an unaccelerated base algorithm \cite{nemirovski1985optimal}.
If an objective function satisfies a H\"{o}lderian error bound condition but it is not strongly convex ($r>2$), then an exponentially-growing schedule $t_k = \lambda e^{\beta k}$ for some $\lambda > 0$ and $\beta>0$ results in a faster sublinear convergence rate.

As notable prior work, \citet{kim2021accelerated} studied APPM with a constant restarting schedule in the strongly monotone setup but was not able to obtain a rate faster than plain PPM. We show that restarting with an exponentially increasing schedule accelerates \eqref{algo:os-ppm} under uniform monotonicity, as for the case of $r>2$ in \citet{roulet2020sharpness}.

\begin{proof}[Proof of \cref{sharp:appm}]
    Suppose that given an initial point $x_0\in\reals^n$, let $\tilde{x}_0$ be an iterate generated by applying APPM on $x_0$ only once.
    Then
    \[
        \tilde{x}_0 = \frac{1}{2}(2\resA x_0 - x_0) + \frac{1}{2} x_0 = \resA x_0,
    \]
    so we get
    \begin{align*}
        \|x_0 - x_\star\|^2 
        &= 
        \| \tilde{\opA} \tilde{x}_0 + (\tilde{x}_0 - x_\star) \|^2 \\
        &=
        \| \tilde{\opA} \tilde{x}_0 \|^2 + 2 \langle \tilde{\opA} \tilde{x}_0, \tilde{x}_0 - x_\star \rangle + \|\tilde{x}_0 - x_\star\|^2.
    \end{align*}
    From the monotonicity of $\opA$, $\langle \tilde{\opA} \tilde{x}_0, \tilde{x}_0 - x_\star \rangle \ge 0$, so we may conclude that
    \[
        \|\tilde{\opA} \tilde{x}_0\|^2 \le \|x_0-x_\star\|^2.
    \]

    Now we describe the restarting scheme of APPM.
    Let $t_k$ be the number of inner iterations applying APPM for the $k$th outer iteration.
    This iteration starts from $\tilde{x}_{k-1}$ and outputs $\tilde{x}_k$ after applying $t_k$ iterations of APPM.
    Then the $k$th outer iteration results in
    \[
        \|\tilde{\opA} \tilde{x}_k\|^2 
        \le \frac{1}{(t_k+1)^2} \|\tilde{x}_{k-1}-x_\star\|^2
        \le \frac{1}{t_k^2} \|\tilde{x}_{k-1}-x_\star\|^2 
        \le \frac{1}{\mu^{2/\alpha} t_k^2} \|\tilde{\opA} \tilde{x}_{k-1}\|^{2/\alpha},
    \]
    where the last inequality follows from
    \[
        \|\tilde{\opA} \tilde{x}_{k-1}\| \|\tilde{x}_{k-1}-x_\star\| \ge \langle \tilde{\opA} \tilde{x}_{k-1}, \tilde{x}_{k-1}-x_\star \rangle \ge \mu \|\tilde{x}_{k-1}-x_\star\|^{\alpha+1}.
    \]
    
    In order to find a possible choice of restart schedule, we will iteratively find the number $t_k$ of inner iterations for $k$th outer iteration which satisfies
    \[
        \| \tilde{\opA} \tilde{x}_k \|^2 \le e^{-\eta k} \| x_0-x_\star \|^2
    \]
    for some $\eta > 0$.
    The case of $k=0$ holds automatically.
    Suppose $k\ge1$, and $t_1,\dots,t_{k-1}$ are already chosen to satisfy
    \[
        \| \tilde{\opA} \tilde{x}_{k-1} \|^2 \le e^{-\eta (k-1)} \|x_0-x_\star\|^2
    \]
    for $k \ge 1$.
    Then
    \[
        \|\tilde{\opA} \tilde{x}_k\|^2 
        \le 
        \frac{1}{\mu^{2/\alpha} t_k^2} \|\tilde{\opA} \tilde{x}_{k-1}\|^{2/\alpha} 
        \le 
        \frac{1}{\mu^{2/\alpha} t_k^2} e^{-\frac{\eta(k-1)}{\alpha}} \|x_0-x_\star\|^\frac{2}{\alpha},
    \]
    so that the claimed convergence rate is guaranteed if
    \[
        \frac{1}{\mu^{2/\alpha} t_k^2} e^{-\frac{\eta(k-1)}{\alpha}} \|x_0-x_\star\|^\frac{2}{\alpha} 
        \le
        e^{-\eta k} \|x_0-x_\star\|^2.
    \]
    This is equivalent to
    \[
        t_k \ge
        \underbrace{\mu^{-\frac{1}{\alpha}} e^\frac{\eta}{2\alpha} \|x_0-x_\star\|^{\frac{1}{\alpha}-1}}_{:=\lambda} \exp \Big\{ \underbrace{\frac{\eta}{2} \left(1 - \frac{1}{\alpha}\right) }_{:=\beta} k \Big\},
    \]
    so if $t_k\ge \lambda e^{\beta k}$ for $k=1,\dots,R$, then $\|\tilde{\opA}\tilde{x}_k\|^2 \le e^{-\eta k} \|x_0-x_\star\|^2$ for $k=1,\dots,R$.

    Now we prove that the choice of 
    \[
        t_k = \begin{cases}
            \left\lceil \lambda e^{\beta k} \right\rceil & (k=1,\dots,R-1) \\
            N - 1 - \sum_{k=1}^{R-1} t_k & (k=R)
        \end{cases}
    \]
    for integer $R$ satisfying
    \[
        \sum_{k=1}^{R} \lceil \lambda e^{\beta k} \rceil \le N-1 < \sum_{k=1}^{R+1} \lceil \lambda e^{\beta k} \rceil
    \]
    results in $\cO\left(N^{-\frac{2\alpha}{\alpha-1}}\right)$-rate of $\|\tilde{\opA}\hat{x}\|^2$ for restarted OS-PPM \eqref{algo:os-ppm-res}.
    
    For $k=1,\dots,R-1$, $t_k \ge \lambda e^{\beta k}$ by definition of $t_k$.
    If $k=R$, from
    \[
        N-1
        =
        \sum_{k=1}^{R-1} t_k + t_R
        =
        \sum_{k=1}^{R-1} \lceil \lambda e^{\beta k} \rceil + t_R,
    \]
    we have
    \[
        t_R
        =
        N-1 - \sum_{k=1}^{R-1} \lceil \lambda e^{\beta k} \rceil
        \ge
        \lceil \lambda e^{\beta R} \rceil
        \ge
        \lambda e^{\beta R}.
    \]
    Therefore, $t_k \ge \lambda e^{\beta k}$ for $k=1,2,\dots,R$, and we get
    \[
        \|\tilde{\opA}\tilde{x}_R\|^2 \le e^{-\eta R}\|x_0-x_\star\|^2.
    \]
    To find the upper bound to $\|\tilde{\opA}\tilde{x}_R\|^2$ using the inequality above, we obtain a lower bound to $R$.
    From $\lambda e^{\beta k} \le \lceil \lambda e^{\beta k} \rceil \le \lambda e^{\beta k} + 1$
    and $\lceil \lambda e^{\beta R} \rceil \le t_R < \lceil \lambda e^{\beta R} \rceil + \lceil \lambda e^{\beta(R+1)} \rceil$, we have
    \begin{align}
        \sum_{k=1}^R \lambda e^{\beta k} \le N-1 = \sum_{k=1}^R t_k
        = \sum_{k=1}^{R-1} \lceil \lambda e^{\beta k} \rceil + t_R
        \le \sum_{k=1}^{R+1} \lambda e^{\beta k} + R + 1.
        \label{eq:schedule}
    \end{align}
    Using the first inequality in \eqref{eq:schedule}, we have
    \[
        \lambda e^{\beta} \frac{e^{\beta R}-1}{e^\beta -1} \le N-1
    \]
    or equivalently,
    \[
        R \le \frac{1}{\beta} \log \left(\frac{N-1}{\lambda} \frac{e^\beta-1}{e^\beta} + 1 \right).
    \]
    Plugging this upper bound of $R$ to the second inequality of \eqref{eq:schedule}, we get
    \[
        N - 1
        \le 
        \lambda e^{\beta} \frac{e^{\beta(R+1)}-1}{e^\beta-1} + \frac{1}{\beta} \log \left( \frac{N-1}{\lambda} \frac{e^\beta-1}{e^\beta} + 1 \right) + 1.
    \]
    Simplifying this to obtain a lower bound on $R$, we get
    \[
        e^{-\beta} \left\{
            \frac{e^\beta-1}{\lambda e^\beta} \left( N - 2 - \frac{1}{\beta} \log \left( \frac{N-1}{\lambda} \frac{e^\beta-1}{e^\beta} + 1 \right) \right) + 1
        \right\}
        \le e^{\beta R}.
    \]
    Therefore,
    \begin{align*}
        \|\tilde{\opA}\tilde{x}_R\|^2
        &\le
        e^{-\eta R} \|x_0-x_\star\|^2 \\
        &\le
        e^\eta \left\{
            \frac{e^\beta-1}{\lambda e^\beta} \left( N - 2 - \frac{1}{\beta} \log \left( \frac{N-1}{\lambda} \frac{e^\beta-1}{e^\beta} + 1 \right) \right) + 1
        \right\}^{-\frac{\eta}{\beta}} \|x_0-x_\star\|^2 \\
        &=
        \left\{
            {\frac{e^\beta-1}{\lambda e^{2\beta}}} \left(N-2-{\frac{1}{\beta}} 
            \log\left({\frac{e^\beta-1}{\lambda e^\beta}}(N-1) + 1 \right) \right) + {\frac{1}{e^\beta}}
        \right\}^{-\frac{2\alpha}{\alpha-1}} \|x_0-x_\star\|^2 
        \tag{Choose $\eta=2$} \\
        &=
        \cO \left( N^{-\frac{2\alpha}{\alpha-1}} \right)
    \end{align*}
    where $\lambda = \left({\scriptstyle\frac{e}{\mu}}\right)^{\frac{1}{\alpha}} \|x_0-x_\star\|^{-\left(1-\frac{1}{\alpha}\right)}$.
\end{proof}
This is a rate faster than $\mathcal{O}(N^{-\frac{\alpha+1}{\alpha-1}})$-rate of PPM.
Although the monotonicity parameter $\mu > 0$ and $\alpha > 1$ are unknown, one can obtain a suboptimal restart schedule with additional cost for the grid search as in \citet{roulet2020sharpness}, where the total cost for the algorithm is of $\cO\left(N^{-\frac{2\alpha}{\alpha-1}} (\log N)^2 \right)$.

\section{Experiment details} \label{sec:expr-details}
We now describe the experiments of Section~\ref{sec:experiment} in further detail.

\subsection{Experiment details of Section~\ref{exper:toy}}
In the first example, $\opT_\theta$ is constructed with $\theta = 15^\circ$ and $\gamma = \frac{1}{0.95}$, and \eqref{algo:oc-halpern} is applied on $\opT_\theta$ with the same $\gamma=\frac{1}{0.95}$.
In the second example, $\opM$ is constructed with $\mu=0.035$, and \eqref{algo:os-ppm} is applied on $\opM$ with the same $\mu=0.035$.
For both experiments, we use $N=101$ total iterations.
The plots of both experiments display the position of every iterate with markers, when methods started from initial point $y_0 = \begin{bmatrix}1&0\end{bmatrix}^\intercal \in \reals^2$.

\subsection{Experiment details of Section \ref{exper:ct}}

X-ray CT reconstructs the image from the received from a number of detectors.
Reconstruction of the original image is often formulated as a least-squares problem with total variation regularization
\begin{equation}
    \begin{array}{ll}
        \underset{x\in \reals^n}{\mbox{minimize}} & \frac{1}{2}\|Ex-b\|^2 + \lambda\|Dx\|_1,
    \end{array}
    \label{prob:ct}
\end{equation}
where $x\in\reals^n$ is a vectorized image, $E\in\reals^{m\times n}$ is the discrete Radon transform, $b=Ex$ is the measurement, and $D$ is the finite difference operator.
This regularized least-squares problem can be solved using PDHG, also known as the Chambolle--Pock method \cite{chambolle2011first}.
PDHG can be interpreted as an instance of variable metric PPM \cite{he2012convergence}; it is a nonexpansive fixed-point iteration $(x^{k+1}, u^{k+1}, v^{k+1}) = \opT(x^k, u^k, v^k)$ defined as
\begin{align*}
    x^{k+1}
    &=
    x^k - \alpha E^\intercal u^k - \beta D^\intercal v^k \\
    u^{k+1}
    &=
    \frac{1}{1+\alpha} \left( u^k + \alpha E(2x^{k+1} - x^k) - \alpha b \right) \\
    v^{k+1}
    &=
    \Pi_{[-\lambda\alpha/\beta, \lambda\alpha/\beta]} \left( v^k + \beta D(2x^{k+1} - x^k) \right)
\end{align*}
with respect to the metric matrix
\[
    M = \begin{bmatrix}
        (1/\alpha) I & -E^\intercal & -(\beta/\alpha) D^\intercal \\
        -E & (1/\beta)I & 0 \\
        -(\beta/\alpha) D & 0 & (1/\beta)I
    \end{bmatrix}.
\]
Therefore, we apply OHM on $\opT$ as 
\[
    (x^{k+1},u^{k+1},v^{k+1}) = \left(1 - \frac{1}{k+2}\right) \opT(x^k,u^k,v^k) + \frac{1}{k+2}(x^0,u^0,v^0)
    \tag{PDHG with OHM}
\]
and use additional restarting strategy to yield a faster convergence.

In our experiment, we use the a Modified Shepp-Logan phantom image.
We applied PDHG, PDHG combined with OHM, and PDHG combined with restarted OC-Halpern \eqref{algo:os-ppm-res}, where the parameters are given as $\alpha = 0.01$, $\beta = 0.03$ and $\lambda = 1.0$.
We applied restarting with the schedule illustrated in \cref{sharp:appm}, with properly chosen $\lambda>0$ and $\beta>0$.

\begin{figure}[H]
    \centering
    \includegraphics[height=5cm]{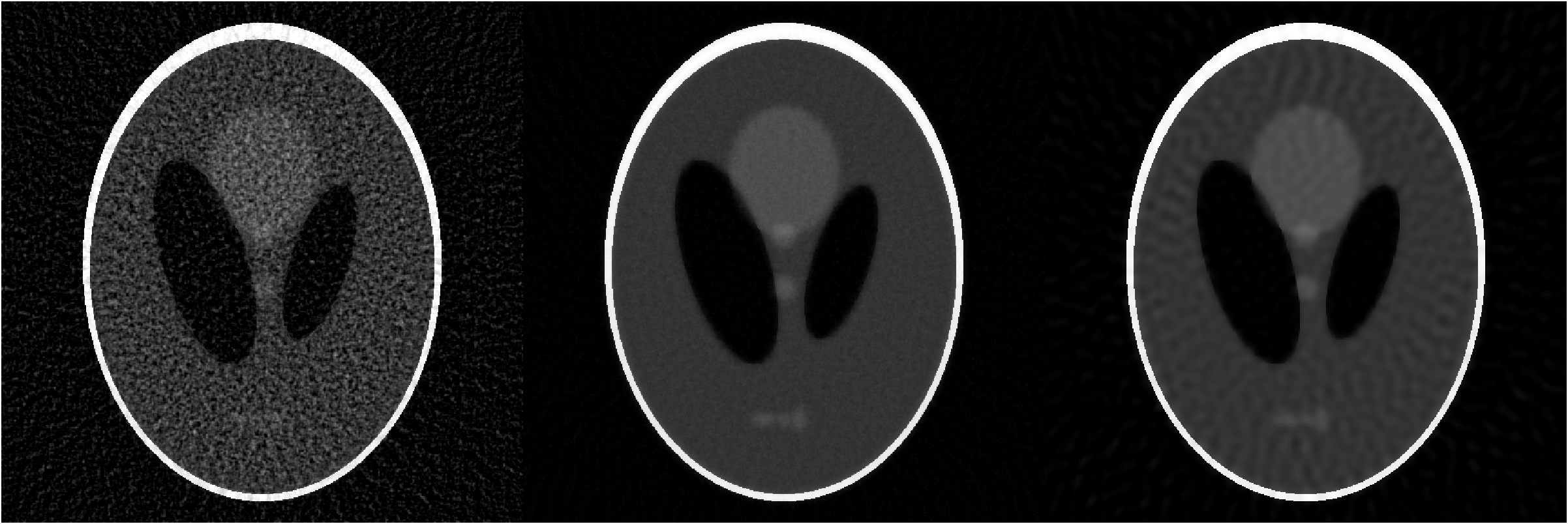}
    \caption{Images reconstructed by applying PDHG, PDHG with OHM, and PDHG with restarted OC-Halpern for 1000 iterations.}
    \label{fig:ct-real-image}
\end{figure}

\begin{figure}[H]
    \centering
    \includegraphics[width=0.5\textwidth]{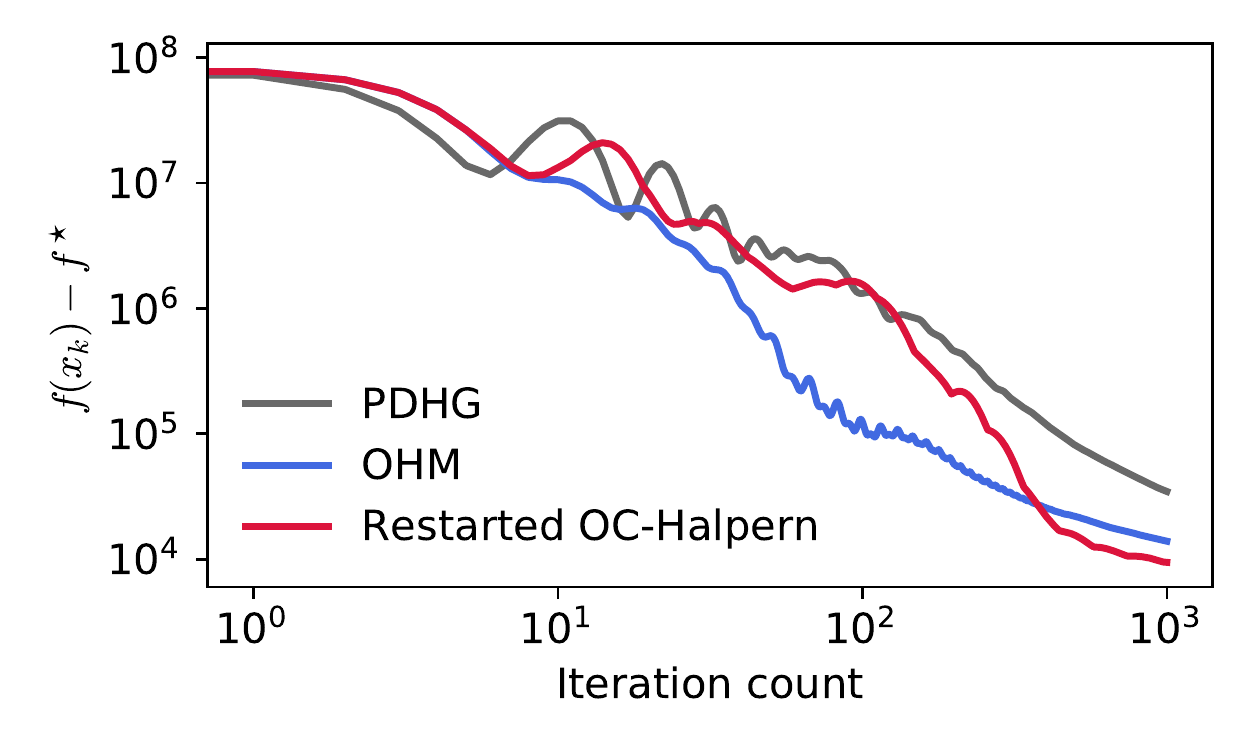}
    \caption{Function value suboptimality $f(x_k)-f^\star$ plot of PDHG, PDHG with OHM, and PDHG with restarted OC-Halpern \eqref{algo:os-ppm-res} in CT image reconstruction.}
    \label{fig:ct-function-value}
\end{figure}

\cref{fig:ct-real-image} shows the reconstructed images after 1000 iterations.
Restarted OC-Halpern \eqref{algo:os-ppm-res} can effectively recover the original image, in a faster rate.
\cref{fig:ct-function-value} shows that even without theoretical guarantee, the function value suboptimality decreases in a faster rate for OHM and restarted OC-Halpern.

\subsection{Experiment details of Section \ref{exper:emd}}
In this section, we approximated the Wasserstein distance (or Earth mover's distance) of two different probability distributions by solving the following discretized problem
\[
    \begin{array}{ll}
        \underset{m_x, m_y}{\mbox{minimize}} & \|\mathbf{m}\|_{1,1} = \sum_{i=1}^n \sum_{j=1}^n |m_{x,ij}| + |m_{y,ij}| \\
        \mbox{subject to} & \mathrm{div}(m) + \rho_1 - \rho_0 = 0.
    \end{array}
\]
To solve this problem \citet{li2018parallel} used PDHG \cite{chambolle2011first}
\begin{align*}
    \tilde{m}_{x,ij}^{k+1} &= \frac{1}{1+\varepsilon \mu} \mathrm{shrink}_1 \left( \tilde{m}_{x,ij}^k + \mu (\nabla \Phi^k)_{x,ij}, \mu \right) \\
    \tilde{m}_{y,ij}^{k+1} &= \frac{1}{1+\varepsilon \mu} \mathrm{shrink}_1 \left( \tilde{m}_{y,ij}^k + \mu (\nabla \Phi^k)_{y,ij}, \mu \right) \\        \Phi_{ij}^{k+1} &= \Phi_{ij}^k + \tau \Big( (\mathrm{div}(2\mathbf{m}^{k+1} - \mathbf{m}^k))_{ij} + \rho_{ij}^1 - \rho_{ij}^0 \Big)
    \tag{Primal-dual method for EMD-$L_1$}
\end{align*}
for $k=1,2,\dots$ where $\tilde{\mathbf{m}} = (\tilde{m}_x, \tilde{m}_y)$ is $\mathbf{m} = (m_x, m_y)$ with zero padding on their last row and last column, respectively, hence making $\tilde{m}_x, \tilde{m}_y \in \reals^{n\times n}$.
We denote this fixed-point iteration by $\opT$, so that  $(\tilde{m}_{x}^{k+1}, \tilde{m}_{y}^{k+1}, \Phi^{k+1}) = \opT (\tilde{m}_{x}^{k}, \tilde{m}_{y}^{k}, \Phi^{k})$.
Combining OHM on this fixed-point iteration yields the iteration
\[
    (\tilde{m}_{x}^{k+1}, \tilde{m}_{y}^{k+1}, \Phi^{k+1}) 
    = 
    \left(1 - \frac{1}{k+2}\right) \opT (\tilde{m}_{x}^{k}, \tilde{m}_{y}^{k}, \Phi^{k})
    + \frac{1}{k+2} (\tilde{m}_{x}^{0}, \tilde{m}_{y}^{0}, \Phi^{0})
\]
for $k=1,2,\dots$, and we also combine restarting technique with exponential schedule to hope for further acceleration.

\begin{figure}[ht]
    \centering
    \subfigure[Probabilistic distribution $\rho_0$]{
        \centering
        \includegraphics[width=0.25\columnwidth]{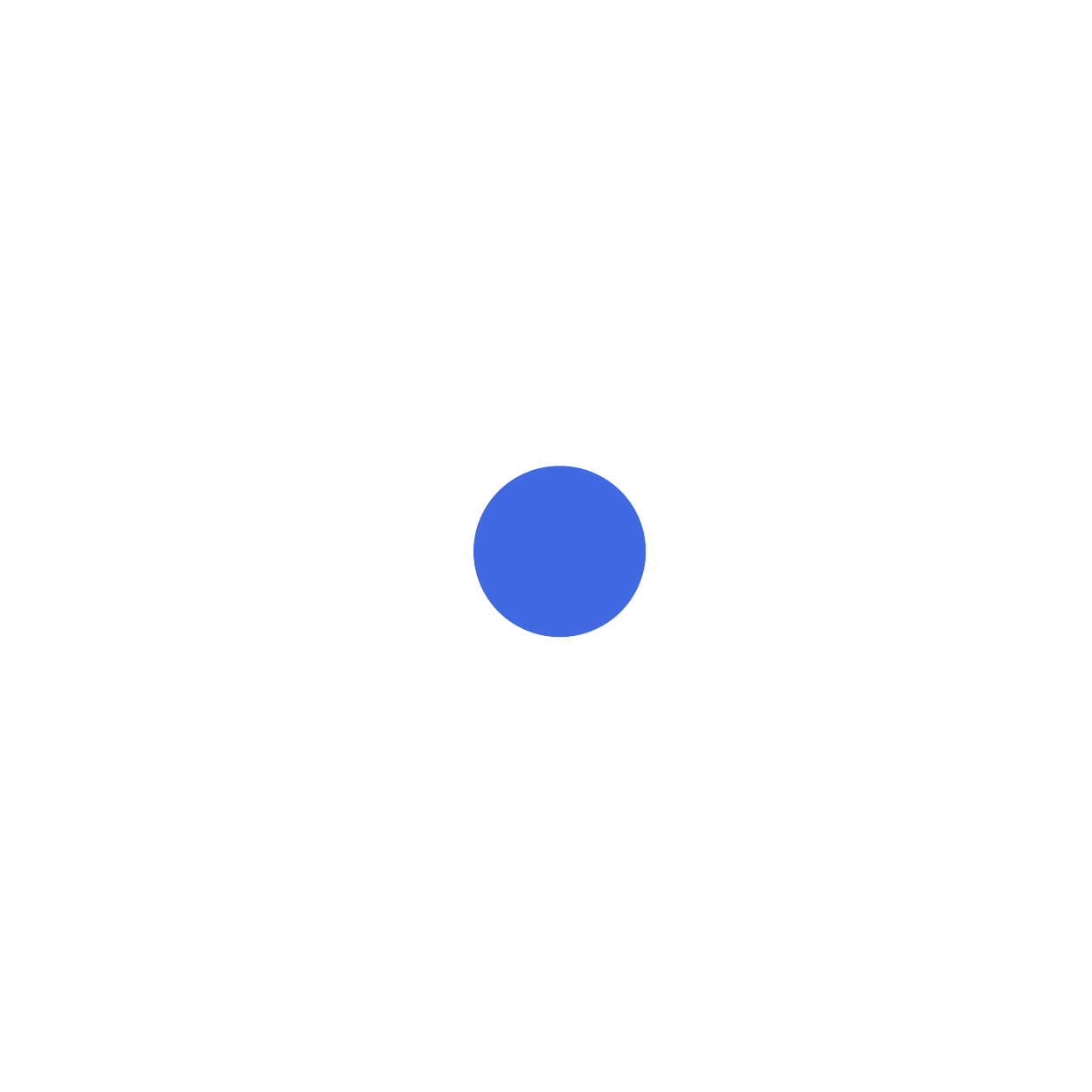}
    }
    \hfill
    \subfigure[Probabilistic distribution of $\rho_1$]{
        \centering
        \includegraphics[width=0.25\columnwidth]{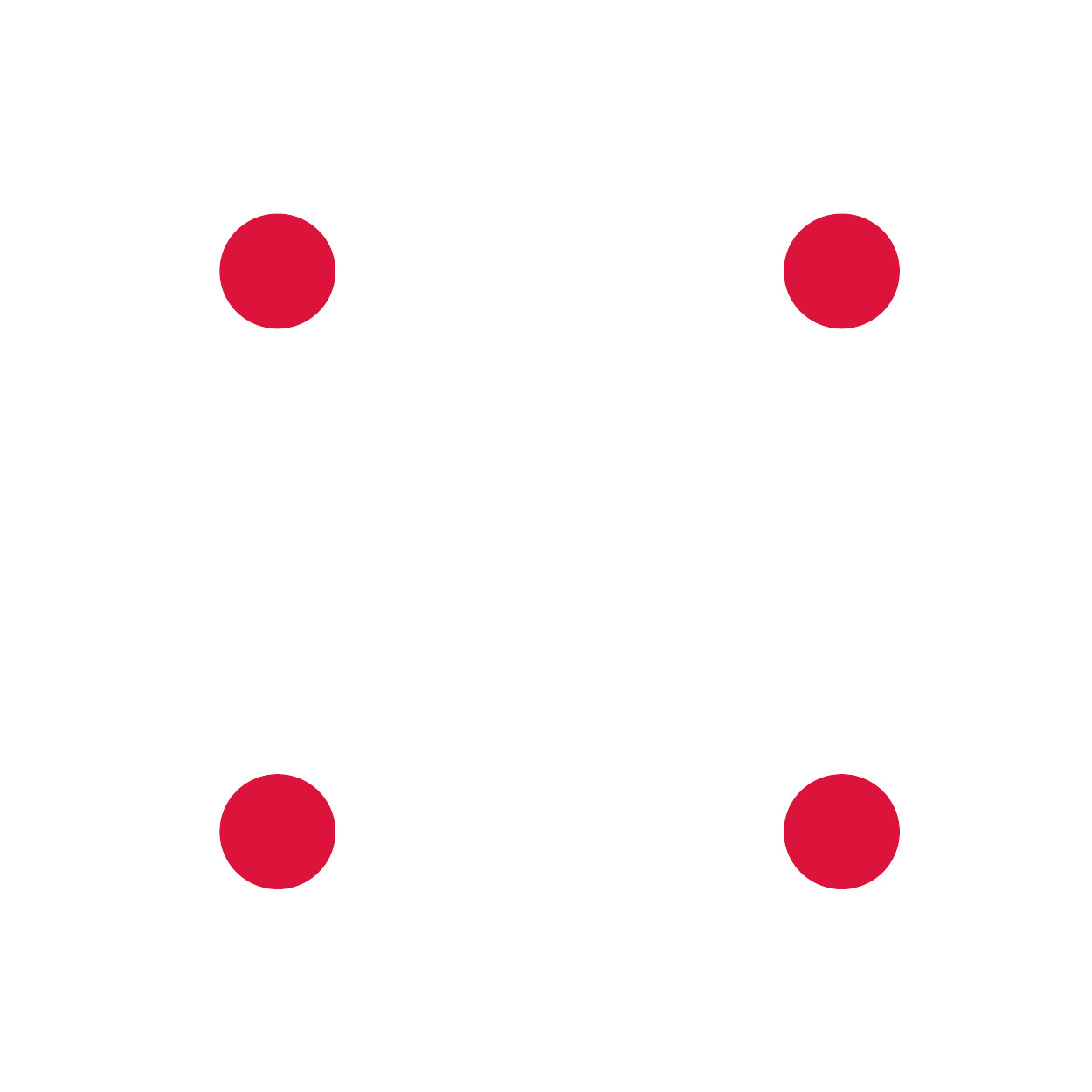}
    }
    \hfill
    \subfigure[Solution of discretized problem in \cref{exper:emd}]{
        \centering
        \includegraphics[width=0.25\columnwidth]{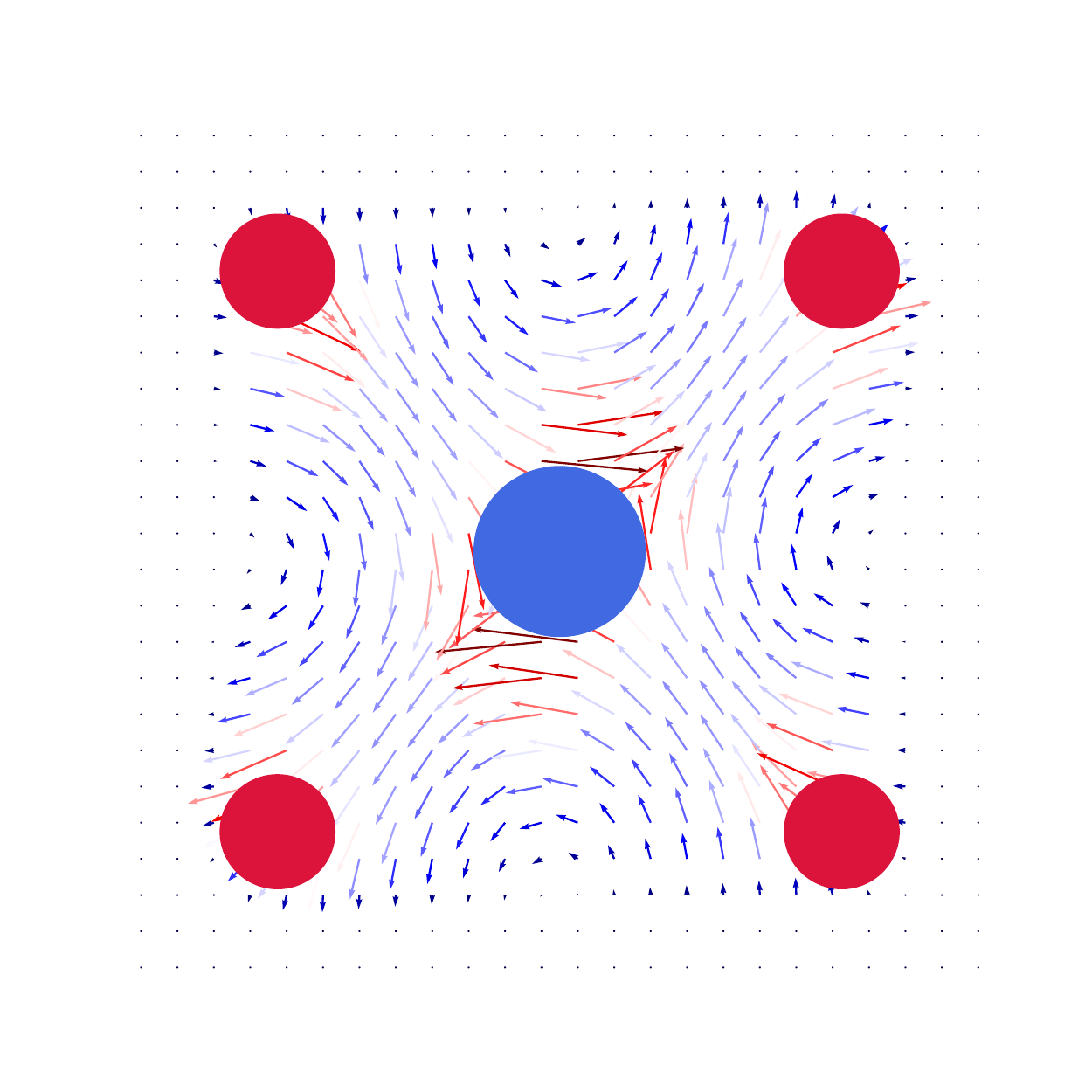}
    }
    \caption{Probabilistic distribution of $\rho_0$ and $\rho_1$.
    This distribution is expressed as the colored parts in $256 \times 256$ grid, there $\rho_0$ contains the part of $x^2 + y^2 \le (0.3)^2$ and $\rho_1$ contains $4$ identical circles with radius $0.2$, centered at $(\pm 1, \pm 1)$.}
    \label{fig:emd-distributions}
\end{figure}

\begin{figure}[ht]
    \centering
    \includegraphics[width=.5\textwidth]{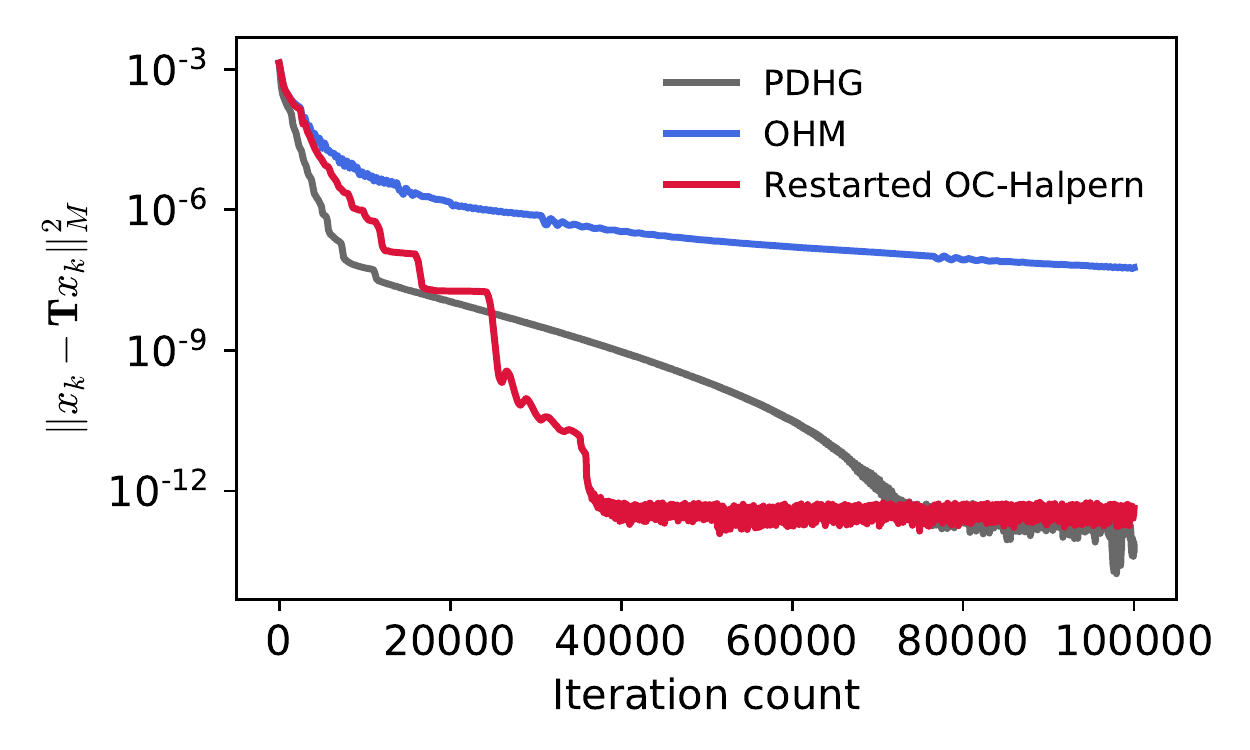}
    \caption{Fixed-point residual of $\opT$ versus iteration count plot for approximating Wasserstein distance.}
    \label{fig:emd-fpr2}
\end{figure}

This experiment calculated an approximation of Wasserstein distance between the two probability distributions as in \cref{fig:emd-distributions}.
We applied 3 different algorithms for $N=100,000$ iterations with algorithm parameters $\mu = 1.0\times10^{-6}$ and $\varepsilon = 1.0$.
Restarting the algorithm with Halpern scheme every $10,000$ iterations provided the accelerated rate, but we chose better exponential schedule for the plot in \cref{fig:emd-fpr2}.

\begin{figure}[ht]
    \centering
    \includegraphics[width=.5\textwidth]{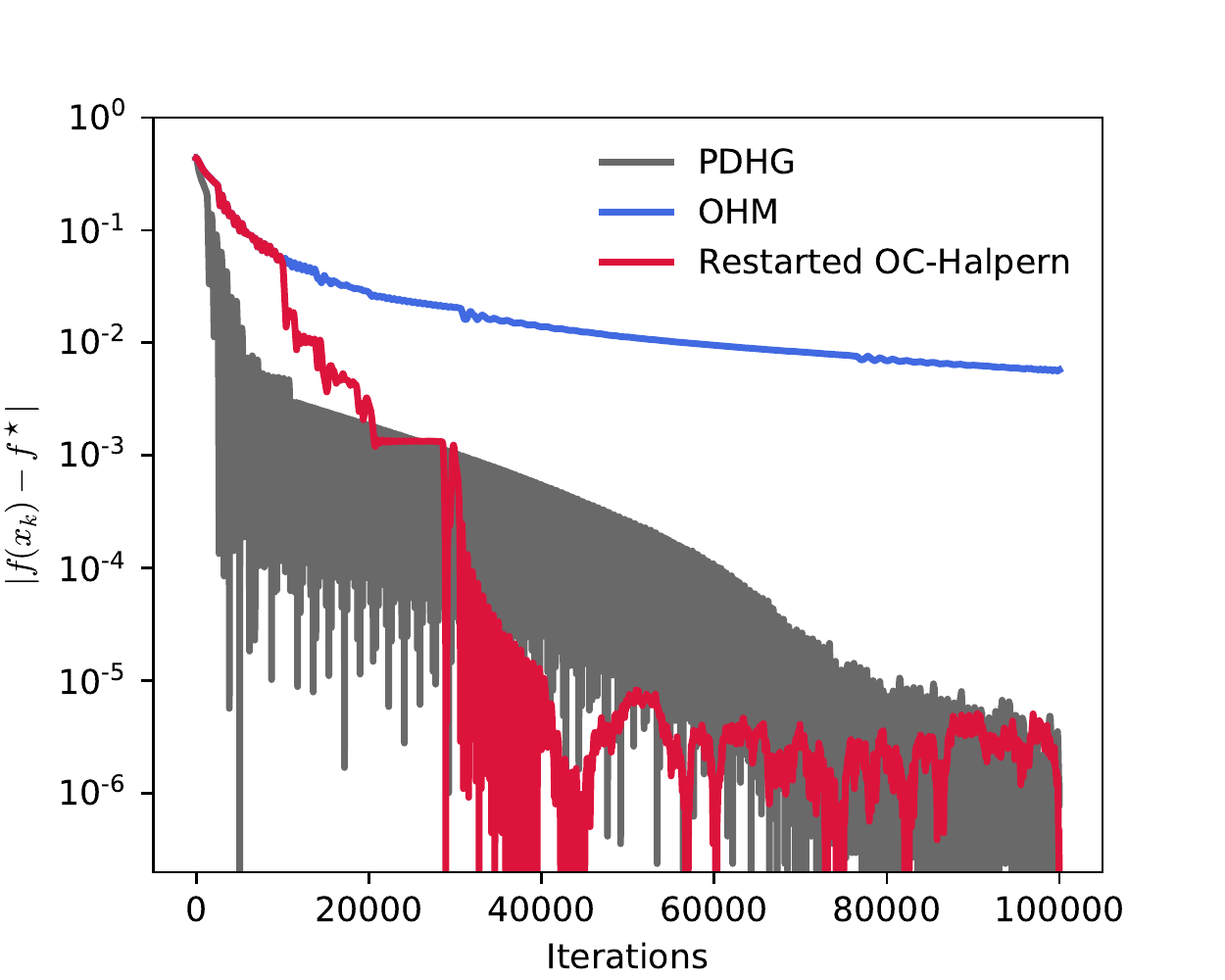}
    \caption{Absolute function-value suboptimality $|f(x_k)-f^\star|$ versus iteration count plot for approximating Wasserstein distance.}
    \label{fig:emd-f-value}
\end{figure}

\subsection{Experiment details of Section \ref{exper:decentralized}}
\begin{figure}[ht]
    \centering
    \includegraphics[width=0.4\textwidth]{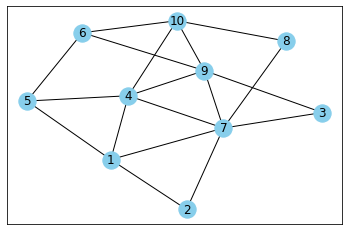}
    \caption{The network underlying the setting of Section \ref{exper:decentralized}.}
    \label{fig:decentralized-graph}
\end{figure}

We follow the settings of decentralized compressed sensing experiment in section IV of \citep{shi2015proximal}.
The underlying network has 10 nodes and 18 edges, and these edges connect the nodes as in \cref{fig:decentralized-graph}.

Experiment considered the regularized least-squares problem on $\reals^{50}$, where the sparse signal $x_\star$ has 10 nonzero entries.
\[
    \begin{array}{ll}
        \underset{x\in \reals^n}{\mbox{minimize}} & \displaystyle \frac{1}{n} \sum_{i=1}^n \|A_{(i)} x - b_{(i)}\|^2 + \lambda \|x\|_1.
    \end{array}
\]
Each node $i$ maintains its local estimate $x_i$ of $x\in\reals^n$, and have access to sensing matrix $A_{(i)}\in\reals^{m_i\times n}$, where $m_i$ is the number of accessible sensors.
Here, we assume to have $m_i=3$ many sensors for each node and has total $m=30$ sensors.

We applied PG-EXTRA, PG-EXTRA combined with OHM, PG-EXTRA with \eqref{algo:oc-halpern}, and PG-EXTRA with Restarted OC-Halpern \eqref{algo:os-ppm-res}, since PG-EXTRA can be understood as a fixed-point iteration \cite{wu2017decentralized}.
Let $\vx^k\in\reals^{n\times 10}$ be a vertical stack of $\reals^n$ vectors, where each $i$\nobreakdash-th row vector $\vx^k_i$ is a local copy of $x$ stored in node $i$.
The vectors in node $i$ only interact with other vectors in close neighborhood of node $i$.
The fixed-point iteration $(\vx^{k+1}, \vw^{k+1}) = \opT (\vx^k, \vw^k)$ is
\begin{align*}
    \vx^{k+1}_i &= \prox_{\alpha\lambda\|\cdot\|_1} \left(
        \sum_{j} W_{i,j} \vx^k_j - \alpha A_{(i)}^\intercal (A_{(i)} \vx^k_i - b_{(i)}) - \vw^k_i
    \right) \\
    \vw^{k+1} &= \vw^{k} + \frac{1}{2}(I - W) \vx^k
\end{align*}
and PG-EXTRA combined with OHM is
\[
    (\vx^{k+1}, \vw^{k+1}) = \left(1 - \frac{1}{k+2}\right) \opT(\vx^k, \vw^k) + \frac{1}{k+2}(\vx^0, \vw^0)
    \tag{PG-EXTRA with OHM}
\]
for $k=0,1,\dots$.
For all these methods, we chose the mixing matrix $W\in\reals^{10\times10}$ to be Metropolis-Hastings weight with each $(i,j)$-entry $W_{i,j}$ being
\[
    W_{i,j} =
    \begin{cases}
        \frac{1}{\max\{\mathrm{deg}(i), \mathrm{deg}(j)\}} & (i\neq j) \\
        1 - \sum_{j\neq i} W_{i,j} & (i=j)
    \end{cases}
\]
where $\mathrm{deg}(i)$ is the number of edges connected to node $i$.
We applied each methods (PG-EXTRA, PG-EXTRA with OHM, PG-EXTRA with \ref{algo:oc-halpern}, and PG-EXTRA with restarted OC-Halpern \eqref{algo:os-ppm-res}) with stepsize $\alpha=0.005$ and regularization parameter $\lambda=0.002$ for 100 iterations.

\begin{figure}[H]
    \centering
    \includegraphics[width=0.6\textwidth]{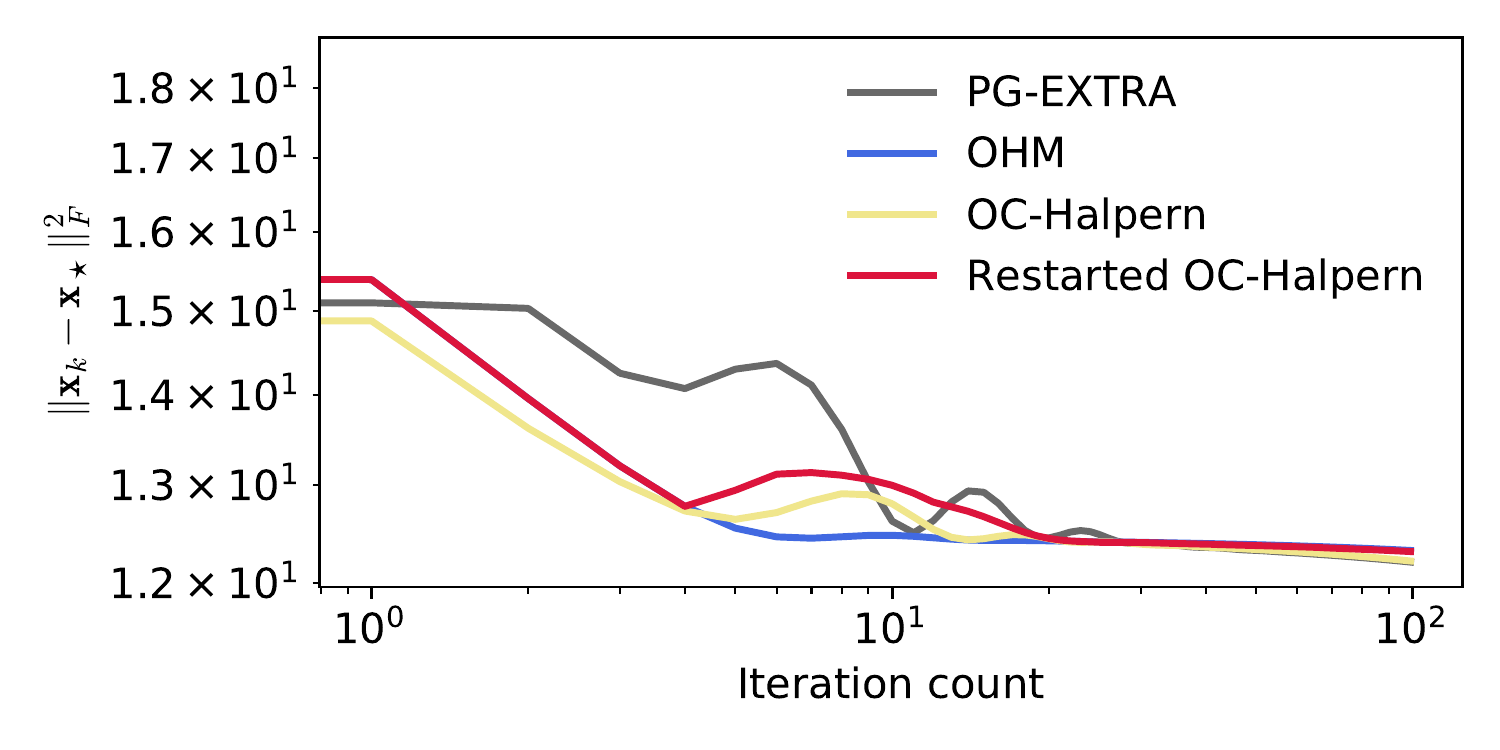}
    \caption{Distance to solution $\|\vx_k-\vx_\star\|^2_F$ versus iteration count plot for PG-EXTRA, PG-EXTRA with OHM, PG-EXTRA with \eqref{algo:oc-halpern}, and PG-EXTRA with Restarted OC-Halpern \eqref{algo:os-ppm-res}.}
    \label{fig:pg-extra-dist}
\end{figure}

\end{document}